\documentclass[12pt,oneside,openany,article]{memoir}
\usepackage{mempatch}

\nouppercaseheads 
\usepackage[amsthm,thmmarks,hyperref]{ntheorem}
\usepackage[noTeX]{mmap}
\usepackage{etex}

\usepackage[english]{babel}
\usepackage[leqno]{mathtools}

\usepackage{mathrsfs}
\usepackage{upgreek}
\usepackage[compress,square,comma,numbers]{natbib}
\usepackage{hypernat} 

\usepackage{sfmath}

\usepackage{microtype}

\usepackage{subfig}
\usepackage{wrapfig}
\usepackage{array}

\usepackage{graphicx,color}
\definecolor{gray}{gray}{0}
\pagecolor{white}

\usepackage{hyperxmp}
\usepackage{xr-hyper}
\usepackage{nameref}
\usepackage[pdftex,bookmarks,pdfnewwindow,plainpages=false,unicode]{hyperref}

\usepackage{bookmark}

\usepackage{enumitem}

\usepackage{amssymb} 

\hypersetup{
colorlinks=true,
linkcolor=black,
citecolor=black,
urlcolor=blue,
pdfauthor={Victor Ivrii},
pdftitle={Multidimensional magnetic Schr\"odinger operator. II.\\
Non-Full-rank case},
pdfsubject={Sharp Spectral Asymptotics},
pdfkeywords={Microlocal Analysis,  Sharp Spectral Asymptotics, Schort loops},
bookmarksdepth={4}
}

\numberwithin{equation}{chapter}

\theoremstyle{plain}
\newtheorem{theorem}{Theorem}[section]

\newtheorem{proposition}[theorem]{Proposition}
\newtheorem{corollary}[theorem]{Corollary}

\theoremstyle{definition}

\newtheorem{definition}[theorem]{Definition}
\newtheorem{Problem}[theorem]{Problem}
\theoremstyle{remark}
\newtheorem{remark}[theorem]{Remark}
\newtheorem{example}[theorem]{Example}

\numberwithin{equation}{section}

\DeclareMathAlphabet{\mathpzc}{OT1}{pzc}{m}{it}

\newcommand{\cA}{\mathcal{A}}
 \newcommand{\cB}{\mathcal{B}}

 \newcommand{\cE}{\mathcal{E}}
 
 \newcommand{\cG}{\mathcal{G}}

 \newcommand{\cK}{\mathcal{K}}
 \newcommand{\cL}{\mathcal{L}}
 
 \newcommand{\cN}{\mathcal{N}}

 \newcommand{\cQ}{\mathcal{Q}}
 
 \newcommand{\cT}{\mathcal{T}}

 \newcommand{\cW}{\mathcal{W}}

 \newcommand{\cZ}{\mathcal{Z}}

 \newcommand{\sC}{\mathscr{C}}

 \newcommand{\sF}{\mathscr{F}}

 \newcommand{\sL}{\mathscr{L}}

\newcommand{\interm}{{\mathsf{intm}}}
\newcommand{\inn}{{\mathsf{inn}}}

\newcommand{\out}{{\mathsf{out}}}

\newcommand{\corr}{{\mathsf{corr}}}

\newcommand{\D}{{\mathsf{D}}}

\newcommand{\MW}{{\mathsf{MW}}}

\newcommand{\R}{{\mathsf{R}}}

\newcommand{\T}{{\mathsf{T}}}

\newcommand{\sU}{{\mathsf{U}}}
\newcommand{\W}{{{\mathsf{W}}}}
\newcommand{\w}{{\mathsf{w}}}

\newcommand{\const}{{\mathsf{const}}}
\newcommand{\dist}{{{\mathsf{dist}}}}

\newcommand{\new}{{{\mathsf{new}}}}

\newcommand{\bC}{{\mathbb{C}}}

\newcommand{\bR}{{\mathbb{R}}}

\newcommand{\bZ}{{\mathbb{Z}}}

\newcommand{\fA}{{\mathfrak{A}}}

\newcommand{\fJ}{{\mathfrak{J}}}

\newcommand{\fm}{{\mathfrak{m}}}
\newcommand{\fM}{{\mathfrak{M}}}
\newcommand{\fn}{{\mathfrak{n}}}
\newcommand{\fN}{{\mathfrak{N}}}

\def\1{\boldsymbol {|}}
\newcommand{\boldtau}{{\boldsymbol{\tau}}}

\newcommand{\Def}{\mathrel{\mathop:}=}

\newcommand{\Hess}{\operatorname{Hess}}

\newcommand{\Ker}{\operatorname{Ker}}
\newcommand{\mes}{\operatorname{mes}}

\newcommand{\Ran}{\operatorname{Ran}}
\newcommand{\rank}{\operatorname{rank}}
\newcommand{\supp}{\operatorname{supp}}

\newcommand{\tr}{\operatorname{tr}}

\newenvironment{claim}[1][{\textup{(\theequation)}}]{\refstepcounter{equation}\vglue10pt
\begin{trivlist}
\item[{\hskip\labelsep#1}]}{\vglue10pt\end{trivlist}}

\newenvironment{claim*}[1][{}]{\vglue10pt
\begin{trivlist}
\item[{\hskip\labelsep#1}]}{\vglue10pt\end{trivlist}}

\newenvironment{phantomequation}[1][]{\refstepcounter{equation}}{}
\newcounter{note}

\DeclareTextCommand{\textinfty}{PU}{\9042\036}

\DeclareTextCommand{\textge}{PU}{\9042\145}
\DeclareTextCommand{\textle}{PU}{\9042\144}
\DeclareTextCommand{\texthat}{PD1}{\136}

\begin{document}
\title{Multidimensional magnetic Schr\"odinger operator. II.\\
Non-Full-rank case}
\author{Victor Ivrii}

\maketitle
{\abstract%
This paper is a continuation of Sections~\ref{book_new-sect-4-5}, \ref{book_new-sect-7-5} and Chapters~\ref{book_new-sect-18}, \ref{book_new-sect-19} of \cite{futurebook}. I derive sharp spectral asymptotics (with the remainder estimate
$O(h^{1-d}+\mu ^r h^{1-r-q})$ for $d$-dimensional Schr\"odinger operator with a strong magnetic field; here $h$ and $\mu$ are Plank and binding constants respectively and magnetic intensity matrix has constant but not full rank $2r$ at each point and $q=d-2r$.

In comparison with version 1 of 5.5 year ago this version contains more results (we also study some degenerations), improvements and some minor corrections.
\endabstract}

\setcounter{chapter}{-1}

\chapter{Introduction}
\label{sect-20-1}
\section{Preface}
\label{sect-20-1-1}

In this Chapter we consider multidimensional Schr\"odinger operator (\ref{book_new-19-1-1}) of \cite{futurebook}

\begin{multline}
A=A_0+V(x),\qquad A_0=\sum_{j,k\le d}P_jg ^{jk}(x)P_k, \\
P_j=hD_j- \mu V_j(x),
\quad h\in (0,1],\ \mu \ge 1.
\label{20-1-1}
\end{multline}
Recall that it is characterized by \emph{magnetic field intensity tensors\/}
$(F_{jk})$ with 
\begin{equation}
F_{jk}=\partial _kV_j-\partial _jV_k,
\label{20-1-2}%
\end{equation}
which is skew-symmetric $d\times d$-matrix, and $(F^j_p)=(g^{jk})(F_{kp})$ which is unitarily equivalent to the skew-symmetric matrix
$(g^{jk})^{\frac{1}{2}}(F_{jk})(g^{jk})^{\frac{1}{2}}$. Then all the eigenvalues of
$(F^j_k)$ (with multiplicities) are $\pm i f_p$ ($f_p>0$, $p=1,\ldots, r$) and $0$ of multiplicity $q=d-2r$ where $2r=\rank (F^j_k)$.

In this Chapter we assume first that the magnetic field intensity matrix has a constant but not full rank:
\begin{gather}
\rank (F_{jk} )=2r= d-q,\quad q\ge 1, \label{20-1-3}\\
\shortintertext{and}
|\bigl((F_{jk})\bigr|_{\cK^\perp} \bigr)^{-1}|\le c_0
\label{20-1-4}
\end{gather}
where $\cK=\Ker (F_{jk})$ and $\cK^\perp =\Ran (F_{jk})$ and (under certain conditions) we derive sharp spectral asymptotics (with the remainder estimate $O(h^{1-d}+\mu ^rh^{1-r-q})$.
The typical example (already studied in Chapters~\ref{book_new-sect-13} and~\ref{book_new-sect-18}) is $3\D$ magnetic Schr\"odinger operator. 

\begin{remark}\label{rem-20-1-1}
Obviously the most interesting case is $q=1$ which is the generic case for skew-symmetric real $d\times d$-matrices. However it does not exclude such generic matrix-valued functions with extra degenerations at some points. 
\end{remark}

Condition (\ref{20-1-3})--(\ref{20-1-4}) is equivalent to (\ref{20-1-3}) and
\begin{equation}
|f_p|\ge \epsilon_0\qquad p=1,\ldots,r;
\label{20-1-5}
\end{equation}
however later we replace it with a \emph{multiscale assumption\/} according to which $f_1,\ldots,f_r$ are broken into groups of different magnitudes which will allow us to analyze generations.

Recall that we consider operator in some domain or on some manifold $X$ with some boundary conditions, assuming that it is self-adjoint in $\sL^2$ and  denote by $e(x,y,\tau)$ Schwartz' kernel of it's spectral projector.

There is a profound difference between full-rank case considered in Chapter~\ref{book_new-sect-19} and non-full-rank case considered here: in its canonical form the full-rank magnetic Schr\"odinger operator is reduced to the family of $r$-dimensional $\mu^{-1}h$-pseudo-differential operators and therefore we can expect asymptotics with the remainder estimate (in the smooth case) as good as 
\begin{equation*}
C(\mu h^{-1})^{r-1} \times \bigl(1+(\mu h)^{-r}\bigr)=
C\bigl(\mu^{-1}h^{1-d}+ \mu^{r-1}h^{1-r}\bigr)
\end{equation*} 
where the first factor is the remainder estimate for a single  scalar $\mu^{-1}h$-pseudo-differential operator and the second factor is the number of these operators to be taken into account; however one needs some non-degeneracy (microhyperbolicity) assumption; otherwise the remainder estimate could be as bad as $C\bigl(\mu h^{1-d}+ \mu^{r}h^{-r}\bigr)$ which is the magnitude of the principal part as $\mu \ge h^{-1}$.

On the other hand, in its canonical form non-full-rank the magnetic Schr\"odinger operator is reduced to the family of $q$-dimensional Schr\"odinger operators which are also $r$-dimensional $\mu^{-1}h$-pseudo-differential operators and in the best case we can expect asymptotics with the remainder estimate 
\begin{equation*}
Ch^{1-q}\times (\mu h^{-1})^r \times \bigl(1+(\mu h)^{-r}\bigr)\asymp+ Ch^{1-d}C\mu^rh^{1-d+r}
\end{equation*}
with the principal part of magnitude 
$h^{-q}+\times (\mu h^{-1})^r \times \bigl(1+(\mu h)^{-r}\bigr)$; however, non-degeneracy condition is not needed at all (if smoothness allows) as 
$q\ge 2$; case $q=1$ is slightly more complicated.

We see that there is a more subtle distinction between $q=1$ and $q\ge 2$; actually $q=2$ might stay apart from $q\ge 3$ due to the lack of the smoothness. There is even more subtle distinction between $r=1$ and $r\ge 2$ but it makes the difference only as $q=0,1$ setting $d=2,3$ a bit apart (and in fact it plays role mainly in the non-smooth case).

\section{Assumptions}
\label{sect-20-1-2}

Note first that $\cK=\Ker F= \Ker \omega_F$ with 
$\omega_F=\sum_{j,k} F_{jk}dx_j\wedge dx_k=\frac{1}{2} d \sum_j V_j dx_j$ and therefore $d\omega_F=0$ and in virtue of Frobenius theorem

\begin{claim}\label{20-1-6}
If $\rank F =d-q$ at each point then $\cK (x) =\Ker F(x)$ is an integrable foliation of $q$-dimensional spaces, which means that for any $y$ there exists a $q$-dimensional manifold $K=K_y\ni y$ such that $T_x K=\cK(x)$ at any point $x\in K$. 
\end{claim}
Then there exists a coordinate system $x=(x''';x^\perp)=(x_1,\ldots,x_q;x_{q+1},\dots,x_d)$ with 
$x^\perp=\const$ along leafs $K_y$ of the above foliation. This (almost) fixes 
$x^\perp$ and also implies that in this coordinate system 
\begin{equation}
F_{jk}=0\qquad \text{as\ \ } j=1,\ldots,q \;\text{or\ \ } k=1,\ldots,q.
\label{20-1-7}
\end{equation}
Then there exists function $\phi$ such that $V_j=\partial_j\phi$ for $j=1,\ldots,q$ and one can eliminate $V_j$ for all $j=1,\ldots,q$ by a gauge transformation\footnote{\label{foot-20-1} Really, we always can eliminate $V_1$; but then $V_{j,1}=V_{1,j}=0$ and $V_j=V_j(x_2,\ldots, x_d)$. Then we can eliminate $V_2$ etc repeating previous arguments for $V_1,\ldots, V_q$. Actually we can eliminate also $V_{q+1}$}. So, we can assume that 
\begin{equation}
V_j=0\; \text{as\ \ } j=1,\ldots,q \quad(\implies V_j=V_j(x^\perp) \; \text{as\ \ } j=q+1,\ldots, d. 
\label{20-1-8}
\end{equation}
Now we assume that these reductions are already made and after these reductions
\begin{phantomequation}\label{20-1-9}\end{phantomequation}
\begin{gather}
g^{jk} \in \sC^{\bar{l},\bar{\sigma}}, \qquad V\in \sC^{l,\sigma}, \tag*{$\textup{(\ref*{20-1-9})}_{1,3}$}\label{20-1-9-1}
\shortintertext{and} 
V_j = \partial_j \phi _j, \qquad \phi_j \in \sC^{\bar{l}+2,\bar{\sigma}}
\tag*{$\textup{(\ref*{20-1-9})}_{2}$}\label{20-1-9-2}
\end{gather}
where the last assumption is a bit stronger than more natural 
$F_{jk}\in  \sC^{\bar{l},\bar{\sigma}}$ or 
$V_\in \sC^{\bar{l}+1,\bar{\sigma}}$; cf. $\textup{(\ref{book_new-19-1-25})}_{1-3}$ of \cite{futurebook}; 
again $(\bar{l},\bar{\sigma})\succeq (l,\sigma)\succeq (1,1)$. Assume also that
\begin{equation}
q=1\implies g^{jd}=\updelta_{jd}, \label{20-1-10}.
\end{equation}
We also assume that
\begin{equation}
|V|\ge \epsilon
\label{20-1-11}
\end{equation}

\section{Canonical Form}
\label{sect-20-1-3}

Recall (see Section~\ref{book_new-sect-13-1-2} of \cite{futurebook}) that if $X={\bR}^d$ and $g^{jk},F_{jk},V$ are constant then the magnetic Schr\"odinger operator is unitarily equivalent to operator
\begin{equation}
\sum_{1\le j\le r} f_j(h^2D_j^2+\mu ^2 x_j^2) +
\sum_{r+1\le j\le d-2r} h^2D_j^2 +V
\label{20-1-12}
\end{equation}
and that   $e(x,x,\tau) =h^{-d}\cN_d^\MW (\tau)$ where 
\begin{multline}
\cN_d^\MW (\tau)\Def\\ \varpi_{d-2r}(2\pi )^{-d+r} \mu ^r h^{r}
\sum _{\alpha \in \bZ^{+r}}
\Bigl(\tau - \sum_j (2\alpha_j +1) f_j\mu h -V\Bigr)_+^{\frac q 2}
f_1\cdots f_r\sqrt g
\label{20-1-13}
\end{multline}
and $\varpi_k$ is a volume of unit ball, $g=\det(g^{jk})^{-1}$.

Now in the general case (for variable $V$ and, may be, $g^{jk},F_{jk}$) we are  interested in the asymptotics of $\int  e(x,x,0)\psi (x)\,dx $ as
$h\to +0,\mu \to +\infty$ where $\psi $ is a fixed function, smooth and compactly supported in $X$.

It is known from Section~\ref{book_new-sect-13-2} of \cite{futurebook} that in the smooth case for $d=3$ one can reduce magnetic Schr\"odinger operator  to its canonical form
\begin{multline}
\mu ^2 \sum_{m+p+j+k\ge 1} a_{m,p,j,k}(x_2,x_3, \mu ^{-1}hD_2) \times
\\
\qquad \qquad\times \bigl(x_1 ^2+\mu ^{-2} h ^2D_1 ^2\bigr)^m (hD_3)^{2p} 
\mu ^{-2m-2j-2p-k}
h^k\label{20-1-14}
\end{multline}
with the ``main part''
\begin{multline}
a_{1,0,0}(x_2,x_3,\mu^{-1}hD_2)\bigl(\mu ^2x_1 ^2+h ^2D_1 ^2\bigr)+\\
h ^2D_3 ^2+a_{0,1,0}(x_2,x_3,\mu ^{-1}hD_2). 
\tag*{$\textup{(\ref{20-1-14})}_0$}\label{20-1-4-0}
\end{multline}
with $a_{1,0,0}=F\circ \Psi_0$, $a_{0,1,0}= V\circ \Psi_0$ and certain diffeomorphism $\Psi_0:T\bR^1\times \bR \to \bR^3$; we ignore terms with $k\ge1$ or $m+p+j\ge 2$. Then one can replace harmonic oscillator $\bigl(h^2D_1^2+\mu^2x_1^2\bigr)$ by one of its eigenvalues $(2\alpha +1)\mu h$.

Multidimensional case is much more tricky. The main part of the canonical form will be 
\begin{multline}
\sum_{j,k\le r} Z_k^* a_{jk}^\w(x'',x''',\mu^{-1}hD'')Z_j +\\
\sum _{2r+1\le j,k\le d}h^2 D_j b_{jk}^\w(x'',x''',\mu^{-1}hD'') D_k +
a_0^\w(x'',x''',\mu^{-1}hD'')
\label{20-1-15}
\end{multline}
with $Z_j=(hD_j -i\mu x_j)$, $Z^*_j=(hD_j +i\mu x_j)$ $j=1,\dots,r$,
$x''=(x_{r+1},\dots, x_{2r})$, $x'''=(2r+1,\dots,d)$. Here $(a_{jk})_{j,k=1,\dots,r}$ is a Hermitian positive definite matrix with eigenvalues $f_1,\dots,f_r$ and $(b_{jk})_{j,k=2r+1,\dots,d}$ is a real symmetric positive definite matrix. We would like to have $(a_{jk})$ in the diagonal form; then instead of (\ref{20-1-15}) we would have
\begin{multline}
\sum_{j} f_j^\w(x'',x''',\mu^{-1}hD'' )\bigl(H_j-\mu h\bigr)+\\
\sum _{2r+1\le j,k\le d}h^2 D_j b_{jk}^\w(x'',x''',\mu^{-1}hD'') D_k+
a_0(x'',x''',\mu^{-1}hD'')^\w
\label{20-1-16}
\end{multline}
with $1\D$-harmonic oscillators $H_j=(h^2D_j^2+\mu^2x_j^2)$.  However second-order resonances $f_j=f_k$ with $j\ne k$ prevent us; in the general case we can assume only that locally
\begin{multline}
\sum_{\fm \in \fM} \sum _{j,k\in \fm } a_{jk}(x'',x''', \mu ^{-1}hD'') Z_j^* Z_k + \\
\sum _{2r+1\le j,k\le d}h^2 D_j b_{jk}^\w(x'',x''',\mu^{-1}hD'') D_k + a_0(x'',\mu ^{-1}hD'',x'')
\label{20-1-17}
\end{multline}
where $\fm \in \fM  $ are disjoint subsets of $\{1,\ldots, r\}$ and eigenvalues of each of matrices $(a_{jk})_{j,k\in \fm }$ are close to one another. This leads to the necessity of the matrix rather than the scalar non-degeneracy (microhyperbolicity) condition.

The reduction of non-principal terms of operator are prevented by higher-order resonances (see Chapter~\ref{book_new-sect-19}); we discuss it later.

\section{Microhyperbolicity}
\label{sect-20-1-4}

Thus we arrive to the final \emph{microhyperbolicity condition\/}\index{microhyperbolicity} at point $\bar{z}=(\bar{x}'',\bar{x}''',\bar{\xi}'')$ and thus at point 
$\bar{x} = \Psi_0 ^{-1}(z)$:

\begin{definition}\label{def-20-1-2} Assume that

\begin{claim}\label{20-1-18}
Set $\{1,\ldots,r\}$ is split into disjoint groups $\fm \in \fM (\bar{z})$ and matrix $(a_{ij})$ is block-diagonal and in a small vicinity of $\bar{z}$ its blocks are close to scalar matrices $f_\fm I_{\#\fm }$ where $f_\fm $ are disjoint.
\end{claim}

Then we call operator \emph{microhyperbolic\/} if  for each real vector
$\boldtau=(\tau_\fm )_{\fm \in \fM }$ such that 
$|\sum _{\fm \in \fM } \tau _\fm  +V|\le \epsilon$ there exists vector 
$\ell= \ell (\bar{z}, \boldtau)\in \bR^{2r+q}$ such that
\begin{gather}
\mu^2 \sum_{j,k} \bigl(\ell (a_{jk}a_0^{-1}) \bigr)\zeta _j^\dag \zeta _k\ge \epsilon_1\qquad \forall \zeta\in \bC^r
\label{20-1-19}\\
\shortintertext{as long as}
a_\fm =\sum _{j,k\in \fm  }a_{jk}\zeta_j^\dag \zeta_k = \mu^2 \tau_\fm 
\qquad  \forall \fm \in \fM .
\label{20-1-20}
\end{gather}
\end{definition}

\begin{remark}\label{rem-20-1-3}
\begin{enumerate}[label=(\roman*),fullwidth]
\item
To impose this condition we need to assume first that symbol $a_0$ does not vanish; thus we need assumption (\ref{20-1-11});

\item
Exactly as in Chapters~\ref{book_new-sect-13} and~\ref{book_new-sect-18} we do not need microhyperbolicity assumption if magnetic field is weak enough;

\item In contrast to Chapters~\ref{book_new-sect-13} and~\ref{book_new-sect-18} we do not need microhyperbolicity condition as $q\ge 3$, and also in $q=2$ if there are no $3$-rd order resonances;

\item  In contrast to Chapter~\ref{book_new-sect-19} we no longer consider $\fN$ partition of $\{1,\ldots,r\}$ joining in one group all possible $3$-rd order resonances as well; the reason is that now we need to consider evolution with 
$|t|\le \epsilon_0$ instead of $|t|\le \epsilon_0\mu $ as it was in Chapter~\ref{book_new-sect-19} and therefore $3$-rd order resonances are not a problem anymore  here.
\end{enumerate}
\end{remark}

In the case of $\mu h\ge \epsilon$ we will need to reformulate microhyperbolicity assumption

\section{Plan of the Chapter}
\label{sect-20-1-5}

In Section~\ref{sect-20-2} we develop weak magnetic field approach to cover both weak magnetic field case and the outer zone in the intermediate magnetic field case. Definition of the ``weak'' magnetic field depends on $q$.

In Section~\ref{sect-20-3} we reduce operator to its canonical form which basically is what we got in the previous Chapter~\ref{book_new-sect-19} albeit with the coefficient depending on all variable $x$ plus a ``free'' Hamiltonian (i.e. without magnetic field) in $q$-dimensional space.

Sections~\ref{sect-20-4} and~\ref{sect-20-5} are  devoted to the remainder estimate with the Tauberian main part of asymptotics in the cases of the intermediate  and a stronger\footnote{\label{foot-20-2} I.e. either strong, or very strong, or superstrong} magnetic field respectively. In all these two and the next two sections presence of $3$-rd order resonances and corresponding non-removable terms poses one of the main obstacles. As $q\ge 3$ ($q\ge 2$ provided there are no such terms) we recover the remainder estimate $O\bigl(h^{1-d}(1+\mu h)\bigr)$ but in the remaining cases the remainder estimate is not as good unless we impose microhyperbolicity or non-degenerateness condition.

Calculations and simplifications of the Tauberian main part  are done in  Sections~\ref{sect-20-6} and~\ref{sect-20-7}  in the cases of the intermediate  and a stronger\footref{foot-20-2}  magnetic field respectively. 

Finally, in sections~\ref{sect-20-8} we consider the case when one pair of the non-zero eigenvalues of $F^j_K$ vanishes on the manifold of codimension $3$. It is more complicated than in $3\D$-case because different eigenvalues have different magnitudes and a simple rescaling does not work.

\chapter{Weak magnetic field case}
\label{sect-20-2}

Analysis in this case follows one of Section~\ref{book_new-sect-18-6} of \cite{futurebook} without serious modifications. We assume that $\mu \le h^{\delta-1}$. Let us  assume that f operator conditions (\ref{20-1-3})--(\ref{20-1-4}), (\ref{20-1-9}) are fulfilled.

\section{Preliminary remarks}
\label{sect-20-2-1}

Let us start from the case (further restrictions to follow)
\begin{gather}
\varepsilon \ge Ch |\log h|,\quad \mu \le \epsilon (h |\log h|)^{-1}.
\label{20-2-1}\\
\shortintertext{Obviously}
[P_j,P_k]=i\mu h F_{jk}.
\label{20-2-2}
\end{gather}
Recall that $\varepsilon$ is a mollification parameter. Consider in the framework of (\ref{20-1-7}) commutators
\begin{align}
h ^{-1}[A,x_m]=&\sum_k \bigl(g ^{mk}P_k+P_kg ^{mk}\bigr),\qquad &&m=1,\ldots, d,
\label{20-2-3}\\
h ^{-1}[A,P_m]=& h ^{-1}\sum_{j,k}P_j[g ^{jk},P_m] P_k+h^{-1}[V,P_m]
&&m=1,\ldots,q,
\label{20-2-4}
\end{align}
and
\begin{multline}
\mu ^{-1}h ^{-1}[A,P_m]=
\sum_{j,k} \bigl(g^{jk}F_{jm}P_k + P_j g^{jk}F_{km}\bigr)+\\
\mu ^{-1}h ^{-1}\sum_{jk}P_j[g ^{jk},P_m] P_k
+ \mu ^{-1}h ^{-1}[V,P_m]\qquad  m=q+1,\ldots,d.  \label{20-2-5}
\end{multline}
These commutators have symbols of the class $\sF ^{l,\sigma }$ in any domain
$\{(x,\xi):\ a(x,\xi )\le c\}$ of the bounded energy where
\begin{multline}
a(x,\xi )= a_0(x,\xi )+V, \qquad
a_0(x,\xi)=\mu^2 \sum_{j,k}g ^{jk}p_j(x,\xi )p_k(x,\xi ),\\
\label{20-2-6}
p_j(x,\xi )\Def(\xi _j- V_j);
\end{multline}
recall that classes $\sF^{l,\sigma}$ are introduced in Chapter~\ref{book_new-sect-18}. To finish the preliminary remarks let us notice that the following statement holds

\begin{proposition}
\footnote{\label{foot-20-3} Microlocal boundaries for $P_j$. Compare with proposition~\ref{book_new-prop-18-6-1} of \cite{futurebook}.}\label{prop-20-2-1} 
Let $f\in \sC_0 ^\infty(\bR ^d)$, $f=1$ in $B(0,C_0)$ with large enough constant $C_0$. Let $T\ge h ^{1-\delta }$.
Then under condition \textup{(\ref{20-2-1})} with large enough constant $C=C_s$
\begin{multline}
|F_{t\to h ^{-1}\tau }
\Bigl(1-f\bigl(p_1(x,\xi ),p_2(x,\xi ),\dots, p_d(x,\xi )\bigr)\Bigr)^\w
\chi _T(t)u |\le Ch ^s\\
\forall x,y\in B(0,{1\over 2})\qquad \forall \tau \le c.\label{20-2-7}
\end{multline}
\end{proposition}

\begin{proof}
One can see easily that symbol 
$f\bigl(p_1(x,\xi ),p_2(x,\xi ),\dots, p_d(x,\xi )\bigr)$ is quantizable under condition (\ref{20-2-1}) and the proof of proposition \ref{prop-20-2-1} is rather obvious. Further standard details we leave to the reader.
\end{proof}

\section{Tauberian estimate}
\label{sect-20-2-2}

\subsection{Outer zone: general case}
\label{sect-20-2-2-1}

Let us prove several results about propagation: 

\begin{proposition}\label{prop-20-2-2}\footnote {\label{foot-20-4} Finite speed of propagation with respect to $x$; cf. propositions~\ref{book_new-prop-18-2-4} and~\ref{book_new-prop-19-2-8} of \cite{futurebook}.}
Let $\mu \le \epsilon (h |\log h|)^{-1}$ and let
\begin{equation}
M \ge \sup _{\sum _{jk}g^{jk}\xi _j \xi _k +V = 0} 2\sum _k g^{jk}\xi _k
+\epsilon
\label{20-2-8}
\end{equation}
with arbitrarily small constant $\epsilon >0$.

Let $\phi _1$ be supported in $B(0,1)$, $\phi _2=1$ in $B(0,2)$, $\chi $ be
supported in $[-1,1]$\,\footnote {\label{foot-20-5} Recall  that we pick up all such auxiliary functions to be admissible in the sense of Section~\ref{book_new-sect-2-3} of \cite{futurebook}.}. Finally, let 
$T_* = Ch |\log h| \le T \le T^*=\epsilon _0$.

Then
\begin{gather}
|F_{t\to h^{-1}\tau } \bar{\chi}_T(t)   \psi_{2x} U(x,y,t)\,\psi_{1y} | \le Ch^s\qquad
\forall x,y\ \forall \tau \le \epsilon_1
\label{20-2-9}
\\
\shortintertext{for}  
\psi_{2x}=\bigl(1- \phi _{2, MT} (x-\bar{x})\bigr), \quad
\psi_{1y}=\phi_{1, MT} (y-\bar{x})\notag
\end{gather}
where here and below $\epsilon_1>0 $ is a small enough constant.
\end{proposition}

\begin{proposition}\label{prop-20-2-3}\footnote {\label{foot-20-6} Finite speed of propagation with respect to $P'''$; cf. proposition~\ref{book_new-prop-18-6-2} of \cite{futurebook}.}
Let $\mu \le \epsilon h^{\delta -1}$ and let
\begin{equation}
T_* = C\varepsilon ^{-1} h |\log h| \le T \le T^*=\epsilon _0.
\label{20-2-10}
\end{equation}
Then
\begin{gather}
|F_{t\to h^{-1}\tau } \bar{\chi} _T(t) Q_{2x}
\psi_{2x} U(x,y,t)\psi_{1y} \,^t\!Q_{1y} |  \le Ch^s\qquad
\forall x,y\ \forall \tau \le \epsilon_1
\label{20-2-11}\\
\shortintertext{for} 
Q_{2x}=\bigl(1- \phi _{2, MT} (hD_x '''-\bar{\xi}''')\bigr), \quad
Q_{1y}= \phi _{1, MT} (hD_y''' -\bar{\xi}'''),\notag
\end{gather}
where here and below  $\psi _{1,2}$ are admissible functions supported in 
$B(0, \frac{1}{2})$ and $\psi_{1x}=\psi_1(x)$, $\psi_{2y}=\psi_2(y)$.
\end{proposition}

\begin{proposition}\label{prop-20-2-4}\footnote {\label{foot-20-7} Singularities leave diagonal; cf. proposition~\ref{book_new-prop-18-6-3} of \cite{futurebook}.} Let $\mu \le \epsilon h^{\delta -1}$ and condition \textup{(\ref{20-2-10})} be fulfilled, and let
\begin{gather}
|\bar{\xi}'''|\Def  \rho\ge C\varepsilon ^{-1} h|\log h| + C\mu^{-1}.
\label{20-2-12}\\
\shortintertext{Then for}
T_*\Def C(\mu^{-1}h|\log h|)^{\frac{1}{2}}+ C\rho ^{-2} h|\log h| \le T \le T^*\Def \epsilon \rho
\label{20-2-13}\\
\shortintertext{estimate}
|F_{t\to h^{-1}\tau } \chi _T(t) \Gamma  U(x,y,t)\psi_y \,^t\! Q_y | \le Ch^s\qquad 
\forall \tau \le \epsilon_1
\label{20-2-14}
\end{gather}
holds for $Q_y=\phi_{1, \frac{1}{2}\rho} (hD_y'''-\bar{\xi }''')$. 
\end{proposition}

\begin{proof}[Proofs of propositions~\ref{prop-20-2-2}--\ref{prop-20-2-4}]
Proofs of these propositions repeat proofs of the corresponding propositions of Chapters~\ref{book_new-sect-18} and~\ref{book_new-sect-19} and are left to the reader.
\end{proof}

Now let us pick up 
\begin{equation}
\varepsilon = C\rho^{-1}h|\log h|.
\label{20-2-15}
\end{equation}
Then any  $T\in [\epsilon_0 \mu^{-1}, \epsilon \rho]$ satisfies condition (\ref{20-2-13}) as long as
\begin{equation}
\rho \ge \bar{\rho}_1 \Def C\max\bigl(\mu^{-1},(\mu h|\log h|)^{\frac{1}{2}}\bigr).
\label{20-2-16}
\end{equation}
This is a definition of the \emph{outer zone\/}\index{zone!outer} and $\varepsilon$ is defined by (\ref{20-2-15}) there.

Therefore contribution of zone $\{|\xi''|\asymp \rho\}$ with given 
$\rho\ge \bar{\rho}_1\}$ to the Tauberian remainder does not exceed $Ch^{1-d}T^{*\,-1}\rho^{q-1}$ where $\rho^{q-1}$ is the measure of $\{\xi''':\ a(x,\xi)=0,\ |\xi'''|\asymp \rho\}$; then the total contribution of the outer zone to the Tauberian remainder does not exceed
\begin{equation}
Ch^{1-d}\int \frac{1}{T^*(\rho)} \rho^{q-1}\, d\rho
\label{20-2-17}
\end{equation}
with all integrals here and until the end of the section taken from $\bar{\rho}_1$ to 1, and plugging $T^*(\rho)\asymp \rho$ we get $Ch^{1-d}\int \rho^{q-2}\, d\rho \asymp h^{1-d}$ as $q\ge 2$.

So, we have proven

\begin{proposition}\label{prop-20-2-5}
Let $q\ge 2$, $\varepsilon$ and $\bar{\rho}_1$  be defined by \textup{(\ref{20-2-15})} and \textup{(\ref{20-2-16})}, and $(l,\sigma)=(1,1)$. Then  the contribution of the outer zone
$\{|\xi '''|\ge  \bar{\rho}_1\}$to the Tauberian remainder with $T=\epsilon \mu^{-1}$ is  $O(h^{1-d})$.
\end{proposition}

\subsection{Outer zone: special case $q=1$}
\label{sect-20-2-2-2}

As $q=1$ expression (\ref{20-2-17}) with $T^*(\rho) =\epsilon \rho$ results in  
$Ch^{1-d}|\log \mu |$  and to get rid of this factor $|\log \mu |$ we need to increase $T^*(\rho)$. We consider few cases in the increased generality and complexity.

\medskip\noindent
(i) If $g^{jk}$ and $F_{jk}$ are constant it is easy: we can take 
\begin{equation}
T^*(\rho) = C\rho |\log \rho|^2
\label{20-2-18}
\end{equation}
as in section Section~\ref{book_new-sect-18-6} of \cite{futurebook} provided $(l,\sigma)\succeq (1,2)$ and we can consider propagation in an appropriate time direction (in which $\xi_1$ increases); see Section~\ref{book_new-sect-18-6} of \cite{futurebook} for details. Namely, let us cover $B(0,1)$ by balls of radii 
$\gamma=\rho ^{\delta_1}$ (with a small exponent $\delta_1>0$). In each element one of the following cases holds:
\begin{align}
&|\partial _{x_1} V(x,\xi)|\asymp \varrho \ge \bar{\varrho}
\Def C_0|\log\rho|^{-2},\label{20-2-19}\\[2pt]
&|\partial _{x_1} V(x,\xi)|\le \bar{\varrho}.
\label{20-2-20}
\end{align}
During propagation the classical trajectory does not leave the initial ball. Furthermore, in the  case (\ref{20-2-20}) \ $| {d\xi_1} / {dt}|=|\partial _{x_1}V|$ does not exceed $\bar{\varrho}$ and 
$|\xi_1|\ge \frac{1}{2}\rho - \bar{\varrho}|t|\ge  \frac{1}{3}\rho$ as 
$|t|\le T^*(\rho)$ defined by (\ref{20-2-18}). 

On the other hand, in the case (\ref{20-2-19}) $|  {d\xi_1}/{dt}|$ is larger than $\varrho$ and does not change sign for $|t|\le T^*(\rho)$. 
Then $|\xi_1|\ge  \frac{1}{2}\rho + \varrho |t|$ for 
$0\le \pm t \le T(\rho) $ with an appropriate sign. 

In both cases $x_1$ shifts away from its initial position and for time 
$t$ (of the appropriate sign in the case (\ref{20-2-19})) this shift is of magnitude $\rho |t|+ \varrho |t|^2$ as 
$\epsilon\mu^{-1}= T_* \le \pm t\le T^*(\rho)$ (with an appropriate sign) and we can take a prescribed $T^*(\rho) =\epsilon \rho|\log \rho|^2$ and recover $O(h^{1-d})$ in (\ref{20-2-17}).

\medskip\noindent
(ii) Similar arguments work if $f_j$ have constant multiplicities. Really, we can rewrite then the symbol as 
\begin{equation}
a(x,\xi)=\xi_1^2+\mu^2 \sum_{1\le j\le r} f_j(\eta_{2j-1}^2+\eta_{2j}^2) +V
\label{20-2-21}
\end{equation}
with $\eta_k$ being linear combinations of $p_1,\dots, p_{2r}$ (with $\sF^{1,2}$ coefficients) such that
\begin{equation}
\mu \{\eta_k,\eta_m\}\equiv \Lambda_{km}
\label{20-2-22}
\end{equation}
modulo linear combinations of $\eta_1,\dots, \eta_{2r}$ with $\sF^{0,2}$ coefficients where $\Lambda_{km}=1$ if $k=2j-1, m=2j$ for some $j$, 
$\Lambda_{km}=-1$ if $k=2j, m=2j-1$ for some $j$, $\Lambda_{km}=0$ otherwise.

Then the derivative of $\mu^2 (\eta_{2j-1}+i\eta_{2j})(\eta_{2k-1}- i\eta_{2k})$ along classical trajectories 
\begin{equation}
|\mu^2 \bigl \{a ,(\eta_{2j-1}+i\eta_{2j})(\eta_{2k-1}- i\eta_{2k})\bigr \}|\le M (1+\mu |f_j-f_k|) \qquad \forall j,k 
\label{20-2-23}
\end{equation}
is bounded as $f_j=f_k$ and during time $T$  its  variation  does not exceed $CT$ which is way less than  $\epsilon \bar{\varrho}$. This could be justified on the quantum level as well because one can quantize symbol 
\begin{equation*}
f\bigl(\bar{\varrho}^{-1}\bigl( (\eta_{2j-1}+i\eta_{2j})(\eta_{2k-1}- i\eta_{2k})-\lambda \bigr)\bigr)
\end{equation*}
with function $f$ supported in $[-c,c]$ and $|\lambda|\le c$.

Consider now $\{a, \xi_1\}$; modulo $O(\rho)$ it is equal to
$-\partial _{x_1}V + \mu^2\sum _{j,k} \alpha_{jk} \eta_j\eta_k$ with $\sF^{0,2}$ coefficients. Then we can find $\alpha'_{jk}$ such that
\begin{multline}
\bigl\{a, \xi_1- \mu \sum_{j,k}\alpha_{jk}\eta_j\eta_k \bigr\}\equiv\\
-\partial _{x_d}V +\mu^2\sum _{j,k:f_j=f_k} \alpha'' _{jk}(\eta_{2j-1}+i\eta_{2j})(\eta_{2k-1}- i\eta_{2k})
\label{20-2-24}
\end{multline}
modulo $O(\bar{\varrho})$: we remove all non-resonance terms but previously 
we make $\rho^\delta$-mollification of $\alpha_{jk}$ thus making $O(\bar{\varrho})$ error. But then all the above arguments work since correction $ \mu \sum_{j,k}\alpha_{jk}\eta_j\eta_k$ does not exceed $C\mu^{-1}\le \epsilon \rho$. The role of $-\partial _{x_1}V$ is played now by the right-hand expression of (\ref{20-2-24}).

\medskip\noindent
(iii) Consider now the general case. First of all let us make again $\gamma=\rho^{\delta_1}$-covering and consider one element. Then the symbol
could be rewritten in the form 
\begin{gather}
a(x,\xi)=\xi_1^2+ \sum_\fm  a_\fm  +V,\label{20-2-25} \\ 
\shortintertext{where}
a_\fm = \mu^2 \sum_{j,k\in \fm }
a_{jk}(\eta_{2j-1}-i\eta_{2j})(\eta_{2k-1}+i\eta_{2k})
\label{20-2-26}
\end{gather}
with Hermitian matrices $a_\fm =(a_{jk})_{j,k\in \fm }$ whose eigenvalues 
are $f_j(x)$, $f_j$ and $f_k$ differ by less than $\epsilon$ (more than $2\epsilon$) as $j,k $ belong to the same group $\fm $ (different groups, respectively). 

Consider traceless matrices $b_\fm  = a_\fm  - {\frac{1} {(\# \fm )}}\tr a_\fm  \cdot I_\fm $ with the corresponding unit matrices $I_\fm $ and matrices
$c_\fm = \partial_{x_1} b_\fm $ calculated at some fixed point of this partition element. Without any loss of the generality one can assume that each $c_\fm $ is a diagonal matrix; let classify  its eigenvalues  into groups: $\lambda_j$ and $\lambda_k$ differ less than $\epsilon \bar{\varrho}$ (more than $2\epsilon \bar{\varrho}$) as $j,k $ belong to the same subgroup $\fn \subset \fm $ (different subgroups, respectively). 

Then one can cover each element of partition by $\gamma_1=\rho^{\delta_1+\delta_2}$ subelements of two kinds:

\smallskip\noindent
(a) With $f_j$, $f_k$ different by less than $\gamma_1$ for some $j,k$ belonging to the different subgroups. One can prove easily that the total measure of this type elements (i.e. of all elements together) does not exceed $C\bar{\varrho}$ and thus their total contribution to the Tauberian remainder does not exceed $Ch^{1-d}|\log \rho|^{-2}$ which after integration over $ {d\rho}/{\rho}$ results in  $O(h^{1-d})$.

\smallskip\noindent
(b) With $f_j$, $f_k$ different by more than $\gamma_1$ for all $j,k$ belonging to the different subgroups. Then we can rewrite (\ref{20-2-25}) with $\fn$ (which indicates the finer partition) instead of $\fm $  and to prove that 
$|\{a,a_\fn \}|\le C\gamma^{-\delta_3}$ where exponent $\delta_3$ could be made arbitrarily small (by taking all the previous exponents small). Then $a_\fn $ are conserved modulo $O(\bar{\varrho})$ in the classical evolution. 

Moreover, we can make even finer subpartition depending on subelement such that 
$|f_j-f_k|\le C\bar{\varrho}$ if $j,k$ belong to the same element of it. 

But then 
$f_\fn \mu^2 \sum_{j\in \fn }(\eta_{2j-1}^2+\eta_{2j}^2)$
are conserved modulo $O(\bar{\varrho})$ and then it is also true for 
$a_\fn $ and $a_\fm $ as well and also
\begin{equation*}
\{a, \xi_1- \mu \sum_{j,k}\alpha_{jk}\eta_j\eta_k \}\equiv 
\sum_j \lambda_j (\eta_{2j-1}^2+\eta_{2j}^2)\mod O(\bar{\varrho})
\end{equation*}
for appropriate $\alpha_{jk}$.

But then arguments of (i),(ii) work and the total contribution of elements of type (b) does not exceed $Ch^{1-d}|\log \rho|^{-2}$ as well and after integration over $ {d\rho}/{\rho}$ we get  $O(h^{1-d})$ again. So we have proven

\begin{proposition}\label{prop-20-2-6}
Let $q=1$, $\varepsilon$ and $\bar{\rho}_1$  be defined by \textup{(\ref{20-2-15})} and \textup{(\ref{20-2-16})}, and $(l,\sigma)=(1,2)$.  Then  the contribution of the outer zone $\{|\xi '''|\ge  \bar{\rho}_1\}$ to the Tauberian remainder with $T=\epsilon \mu^{-1}$ is $O(h^{1-d})$.
\end{proposition}

\begin{remark}\label{rem-20-2-7}
Similar arguments were applied in Section~\ref{book_new-sect-7-5} of \cite{futurebook} to study propagation near the boundary.
\end{remark}

\subsection{Conclusion}
\label{sect-20-2-2-3}

We have not specified $\varepsilon$ for $\rho=|\xi'''|\le \bar{\rho}_1$ but it is reasonable to take there $\varepsilon=C\bar{\rho}_1^{-1}h|\log h|$ and therefore
\begin{equation}
\varepsilon = C\left\{\begin{aligned} 
&\rho^{-1} h|\log h| \qquad  &&\rho \gtrsim \bar{\rho}_1,\\
&\bar{\rho}_1^{-1} h|\log h|\qquad &&\rho \lesssim \bar{\rho}_1.
\end{aligned}
\right.
\label{20-2-27}
\end{equation}

\begin{remark}\label{rem-20-2-8}
Note that the threshold for $\bar{\rho}_1$ is 
$\mu =\bar{\mu}= C(h|\log h|)^{-\frac{1}{3}}$: $\bar{\rho}_1=C\mu^{-1}$ as $\mu \le \bar{\mu}$, and $\bar{\rho}_1=C(\mu h |\log h|)^{\frac{1}{2}}$ as $\mu \ge \bar{\mu}$.
\end{remark}

On the other hand,  the standard estimate in the complementary 
\emph{inner zone\/}\index{zone!inner}
$\{|\xi '''|\le \bar{\rho}_1\}$ brings the Tauberian remainder estimate 
\begin{equation}
Ch^{1-d}\times \mu \bar{\rho}_1^q \asymp
Ch^{1-d}\times \mu (\mu h|\log h|)^{\frac{q}{2}}+Ch^{1-d}\mu ^{1-q}.
\label{20-2-28}
\end{equation}
Combining with propositions~\ref{prop-20-2-5} and~\ref{prop-20-2-6} we arrive to 

\begin{proposition}\label{prop-20-2-9}
\begin{enumerate}[label=(\roman*),fullwidth]
\item In the framework of either proposition~\ref{prop-20-2-5} or~\ref{prop-20-2-6}
\begin{multline}
\R^\T \Def |\Gamma (\tilde{e}\psi)(\tau) - h^{-1}\int_{-\infty}^\tau 
F_{t\to h^{-1}\tau'} \bar{\chi}_T(t)\Gamma (U \psi) | \le \\[2pt]
Ch^{1-d}+ C \mu (\mu h|\log h|)^{\frac{q}{2}}h^{1-d}\qquad \forall \tau: |\tau| \le \epsilon_1
\label{20-2-29}
\end{multline}
where the left-hand expression is the total Tauberian remainder with $T=\epsilon \mu^{-1}$.

\item  Furthermore, this estimate holds with any $T\ge \bar{T}$ with
$\bar{T}=Ch|\log h|$.
\end{enumerate}
\end{proposition}

\section{Main theorems}
\label{sect-20-2-3}

\subsection{General case}
\label{sect-20-2-3-1}

Now in virtue of the standard results rescaled we can replace Tauberian expression by multiterm Weyl expression $h^{-d}\tilde{\cN}^\W_{(\infty)}$ defined by (\ref{book_new-19-2-83} of \cite{futurebook}) for approximate operators then we can replace $\R^\T$ by $\tilde{\R}^\W_{(\infty)}$ also for approximate operators (where tilde means that we consider approximate operators). 

Note that the contribution of the outer zone to an approximation error does not exceed
\begin{equation}
C\int \bigl( \frac {h|\log h|} {\rho}\bigr)^l
\bigl(1+ \frac {\mu h|\log h|} {\rho}\bigr)
|\log (\rho^{-1}h|\log h|)|^{-\sigma} \rho ^{q-1}\,d\rho \times h^{-d}
\label{20-2-30}
\end{equation}
and this expression  is $O(h^{1-d})$ provided
\begin{claim}\label{20-2-31}
\underline{Either} \ $q\ge 2$,
$(l,\sigma) \succeq (1,1)$ \ \underline{or} \ $q=1$, $(l,\sigma) \succeq (1,2)$.
\end{claim}

Also one can prove easily that replacing $h^{-d}\cN^\W_{(\infty)}$ by the magnetic Weyl expression $h^{-d}\cN^\MW$ does not increase an error and therefore we arrive to the following

\begin{theorem}\label{thm-20-2-10}
Let conditions \textup{(\ref{20-1-1})}, \textup{(\ref{20-1-5})}, \textup{(\ref{20-1-7})}, \textup{(\ref{20-1-8})}--\textup{(\ref{20-1-11})}, and  \textup{(\ref{20-2-31})} be fulfilled, and $\mu \le h^{\delta -1}$.

Then there exist two framing approximations (see footnote \footref{book_new-foot-18-16} of Chapter~\ref{book_new-sect-18}) of \cite{futurebook} such that for both of them
\begin{multline}
\R^\MW = |\int \Bigl(\tilde{e} (x,x,\tau)-
h^{-d}\cN^\MW (x,\tau )\Bigr)\psi(x)\,dx|\le\\
Ch^{1-d}+C\mu h^{1-d}(\mu h |\log h|)^{\frac{q}{2}};
\label{20-2-32}
\end{multline}
the same estimate holds for $\R^\W_{(\infty)}$.
\end{theorem}

\begin{remark}\label{rem-20-2-11} 
\begin{enumerate}[label=(\roman*), fullwidth]
\item Theorem \ref{thm-20-2-10} uses no non-degeneracy condition and it does not benefit from stronger smoothness assumptions; 

\item For $\mu \le \bar{\mu}_q$,
\begin{gather}
\bar{\mu}_{(q)} \Def  (h|\log h|)^{-q/(q+2)}
\label{20-2-33}\\
\shortintertext{estimate}
\R^\MW \le Ch^{1-d}
\label{20-2-34}
\end{gather}
holds. In particular, $\bar{\mu}_{(1)}=C(h|\log h|)^{-\frac{1}{3}}=\bar{\mu}_1$, 
$\bar{\mu}_{(2)}=C(h|\log h|)^{-\frac{1}{2}}$ and
$\bar{\mu}_{(q)}\ge C(h|\log h|)^{-\frac{3}{5}}$ as $q\ge 3$.
\end{enumerate}
\end{remark}

\begin{Problem}\label{Problem-20-2-12} 
Using rescaling technique get rid off logarithms in the definitions of $\bar{\rho}_1$ and thus in estimates (\ref{20-2-29}), (\ref{20-2-32}) and in (\ref{20-2-33}):
\begin{gather}
\R^\T  \le  Ch^{1-d}+ C  (\mu h)^{1+\frac{q}{2}}h^{-d}
\tag*{$\textup{(\ref*{20-2-29})}^*$}\label{20-2-29-*}\\
\shortintertext{and}
\R^\MW \le
Ch^{1-d}+ C  (\mu h)^{1+\frac{q}{2}}h^{-d}.
\tag*{$\textup{(\ref*{20-2-32})}^*$}\label{20-2-32-*}
\end{gather}
\end{Problem}

\subsection{Microhyperbolic case}
\label{sect-20-2-3-2}

To push $\mu$ in theorem \ref{thm-20-2-10} up\footnote{\label{foot-20-8} Actually we need to do this for $q=1$ only.} we need a microhyperbolicity condition and some extra smoothness.

First of all, let us improve propagation of singularities:

\begin{proposition}\label{prop-20-2-13}  
Let  $\bar{\mu}_{(q)}\le \mu \le \epsilon h^{-1}|\log h|^{-1}$ and microhyperbolicity condition (see definition~\ref{def-20-1-2}) be fulfilled. Then for  $C\varepsilon ^{-1}h|\log h|\le T \le T^* =\epsilon_1$ estimate \textup{(\ref{20-2-14})} holds.
\end{proposition}

\begin{proof}
An easy proof based on technique of subsections~\ref{book_new-sect-18-6-6} and~\ref{book_new-sect-19-2-2} of \cite{futurebook} we leave to the reader.
\end{proof}

In the inner  zone $\{\rho \le \bar{\rho}_1\}$ we can select 
$\varepsilon = C\mu h|\log h|$ (as  $\mu \ge (h|\log h|)^{-\frac{1}{3}}$) and then for $Q$ with the symbol supported there estimate (\ref{20-2-14})  holds for 
$T_*=\epsilon \mu^{-1} \le T\le T^*=\epsilon_1$.

Redefining in the outer zone $\varepsilon$, so finally 
\begin{equation}
\varepsilon = C\left\{\begin{aligned}
&\max\bigl(\mu h|\log h|\bigl(\frac {\bar{\rho}_1} {\rho}\bigr)^s , 
\rho ^{-1}h|\log h|\bigr)
\qquad && \rho\ge \bar{\rho}_1,\\
&\bar{\varepsilon}\Def \mu h|\log h| \qquad && \rho\le \bar{\rho}_1
\end{aligned}\right.
\label{20-2-35}
\end{equation}
and combining with proposition \ref{prop-20-2-4} we get then that for $Q$ with the symbol supported in the outer zone estimate also (\ref{20-2-14}) holds for 
$T_* \le T\le T^*$.

Thus we arrive to

\begin{proposition}\label{prop-20-2-14} 
Let  $(h|\log h|)^{-\frac{1}{3}}\lesssim \mu \le  \epsilon (h|\log h|)^{-1}$ and microhyperbolicity condition (see definition~\ref{def-20-1-2}) be fulfilled. Let $\varepsilon $ be defined by \textup{(\ref{20-2-35})}. 

Then estimate \textup{(\ref{20-2-14})} holds with $Q=I$ and  any 
$T\in [T_*, T^*]$, $T_*=\epsilon \mu^{-1}$, $T^*=\epsilon_1$.
\end{proposition}

Combining with the standard results rescaled and applying the standard Tauberian technique we arrive to estimate (\ref{20-2-29}) again.

Note that an approximation error in the operator does not exceed 
\begin{equation}
C\varepsilon^q |\log \varepsilon|^{-\sigma}| +
C\mu \varepsilon^{\bar{l}+1}  |\log \varepsilon|^{-\bar{\sigma}} 
\label{20-2-36}
\end{equation}
where one can skip the second term as $\varepsilon \lesssim \mu^{-1}$ which is the case in the inner zone iff $\mu \lesssim (h|\log h|)^{-\frac{1}{3}}$; then it is true in the outer zone as well.

Then one can prove easily  that an approximation error in both $\R^\MW$ and $\R^\W_{(\infty)}$   does not exceed (\ref{20-2-35}) calculated as 
$\varepsilon= \bar{\varepsilon}$, multiplied by $\bar{\rho}_1^{\frac{q}{2}}h^{-d}$ i.e. 
\begin{equation}
C(\mu h|\log h|)^{l+\frac{q}{2}}|\log h|^{-\sigma}h^{-d} 
+C (\mu h|\log h|)^{2+\bar{l}+\frac{q}{2}} |\log h|^{-1-\bar{\sigma}}h^{-1-d}
\label{20-2-37}
\end{equation}
where one can skip the second term as  
$\mu \lesssim (h|\log h|)^{-\frac{1}{2}}$. This proof uses microhyperbolicity assumption as $q=1$ and we consider $\R^\MW$.

However if $\ell_\perp=0$ we can take mollification with 
$\varepsilon = Ch\rho^{-1}|\log h|$ with respect to all variables $x$ and also an extra mollification with $\varepsilon$ defined by (\ref{20-2-35}) with respect to $x'''$; then an approximation error will be given  by (\ref{20-2-37}) but without the second term.

One can prove easily that the difference between $h^{-d}\cN^\MW$ and 
$h^{-d}\cN^\W_{(\infty)}$ does not exceed (\ref{20-2-37}) without second term.

This implies 

\begin{theorem}\label{thm-20-2-15}
Let conditions \textup{(\ref{20-1-1})}, \textup{(\ref{20-1-5})}, \textup{(\ref{20-1-7})}, \textup{(\ref{20-1-8})}--\textup{(\ref{20-1-11})}, and  \textup{(\ref{20-2-31})} be fulfilled, and $\mu \le h^{\delta -1}$.

Let the microhyperbolicity condition (see definition~\ref{def-20-1-2}) be fulfilled. Then there exist two framing approximations (see footnote \footref{book_new-foot-18-16} of Chapter~\ref{book_new-sect-18}) of \cite{futurebook} such that for both of them

\begin{enumerate}[label=(\roman*), fullwidth]
\item  The following estimate holds
\begin{multline}
\R^\MW \le Ch^{1-d}+\\
C(\mu h|\log h|)^{l+\frac{q}{2}}|\log h|^{-\sigma}h^{-d} 
+C  (\mu h|\log h|)^{2+\bar{l}+\frac{q}{2}} |\log h|^{-1-\bar{\sigma}}h^{-1-d};
\label{20-2-38}
\end{multline}
\item  Furthermore, if \underline{either} 
$\mu \le \bar{\mu}_2=C(h|\log h|)^{-\frac{1}{2}}$ \underline{or} $\ell_\perp=0$ in the microhyperbolicity condition, then one can skip the last term in the right-hand expression of \textup{(\ref{20-2-38})}: 
\begin{equation}
\R^\MW \le Ch^{1-d}+
C (\mu h |\log h|)^{l+ \frac{q}{2}}h^{-d}|\log h|^{-\sigma};
\label{20-2-39}
\end{equation}
\item  In particular, $\R^\MW\le Ch^{1-d}$ as $q=1$, 
$(l,\sigma)=(\frac{3}{2}, 1)$ and $\mu \le \bar{\mu}_2$.
\end{enumerate}
\end{theorem}

Recall that all these estimates hold for $\R^\W_{(\infty)}$ as well.

\begin{Problem}\label{Problem-20-2-16} 
Using rescaling technique replace $\mu h|\log h|$ by $\mu h$ in the definitions of $\bar{\rho}_1$ and $\varepsilon$ and thus in estimates (\ref{20-2-38}) and  (\ref{20-2-39}):
\begin{multline}
\R^\MW \le Ch^{1-d}+\\
C(\mu h)^{l+\frac{q}{2}}|\log h|^{-\sigma}h^{-d} 
+C  (\mu h)^{2+\bar{l}+\frac{q}{2}} |\log h|^{-\bar{\sigma}}h^{-1-d}
\tag*{$\textup{(\ref*{20-2-38})}^*$}\label{20-2-38-*}
\end{multline}
and
\begin{equation}
\R^\MW \le Ch^{1-d}+
C (\mu h )^{l+ \frac{q}{2}}h^{-d}|\log h|^{-\sigma};
\tag*{$\textup{(\ref*{20-2-39})}^*$}\label{20-2-39-*}
\end{equation}
\end{Problem}

\subsection{Case of constant $g^{jk}$, $F_{jk}$}
\label{sect-20-2-3-3}

Consider now the case of constant $g^{jk}$ and $F_{jk}$; then without any loss of the generality one can assume that 
\begin{claim}\label{20-2-40}
$g^{jk}=\updelta_{jk}$ and $F_{jk} = \pm f_i$ for $j=q+i, k=q+i+r$ and $j=q+i+r, k=q+i$ with $i=1,\ldots, r$ and $F_{jk}=0$ otherwise.
\end{claim}
Then the microhyperbolicity assumption is equivalent to 
\begin{gather}
|\nabla V|\ge \epsilon_0.
\label{20-2-41}\\
\intertext{We replace it by the non-degeneracy condition}
|\nabla V|\le \epsilon_0\implies |\det \Hess V|\ge \epsilon_0.
\label{20-2-42}
\end{gather}
We leave to the reader the following

\begin{Problem}\label{Problem-20-2-17} 
Let conditions \textup{(\ref{20-1-1})}, \textup{(\ref{20-1-5})}, \textup{(\ref{20-1-7})}, \textup{(\ref{20-1-8})}, $\textup{(\ref{20-1-9})}_3$, \textup{(\ref{20-1-11})}, and  \textup{(\ref{20-2-31})} be fulfilled, and $\mu \le h^{\delta -1}$.

Assume that $g^{jk}$, $F_{jk}$ are constant and (\ref{20-2-42}) is fulfilled.  

\medskip\noindent
\begin{enumerate}[label=(\roman*), fullwidth]
\item  Using arguments similar to those of subsection~\ref{book_new-sect-19-2-5} of \cite{futurebook} prove that for appropriate framing approximations
\begin{equation}
\R^\T \le Ch^{1-d}+ C\mu (\mu h|\log h|)^{q+r} h^{1-d}
\label{20-2-43}
\end{equation}
where $q+r = \frac{1}{2}q+\frac{1}{2}d$ and
\begin{multline}
\R^\MW \le Ch^{1-d}+\\
C\mu (\mu h|\log h|)^{q+r} h^{1-d}+
C(\mu h|\log h|)^{l+\frac{q}{2}}|\log h|^{-\sigma}h^{-d} ;
\label{20-2-44}
\end{multline}
\item Adding rescaling arguments prove that for appropriate framing approximations
\begin{gather}
\R^\T \le Ch^{1-d}+ C (\mu h)^{q+r+1} h^{-d}
\tag*{$\textup{(\ref*{20-2-43})}^*$}\label{20-2-43-*}\\
\shortintertext{and}
\R^\MW \le Ch^{1-d}+C(\mu h)^{q+r+1} h^{-d}+
C(\mu h)^{l+\frac{q}{2}}|\log h|^{-\sigma}h^{-d}. 
\tag*{$\textup{(\ref*{20-2-44})}^*$}\label{20-2-44-*}
\end{gather}
\end{enumerate}
\end{Problem}

\chapter{Canonical form}
\label{sect-20-3}

\section{Pilot-model}
\label{sect-20-3-1}

Assume temporarily that $g^{jk}$ and $F_{jk}$ are constant. Then without any loss of the generality we can assume  $V_j(x)$  are linear functions. 

Then $A^0$ is transformed into exactly
\begin{equation}
\sum_{1\le k\le q} h^2 D_{k}^2+
\sum_{1\le j\le r} f_j\Bigl(h^2 D_{r+j}^2+\mu ^2 x_{r+j}^2\Bigr)
\label{20-3-1}
\end{equation}
by $\mu ^{-1}h$-metaplectic transformation which consists of the following steps:

\begin{enumerate}[label=(\roman*), fullwidth]
\item  Change of variables $(x,\mu^{-1}hD) \mapsto (Qx, \,^t\!Q^{-1}\mu^{-1}hD)$ transforming $g^{jk}$ into $\updelta_{jk}$ and $F$ into matrix satisfying (\ref{20-2-40}). It transforms $V(x)$ into $V(Qx)$. 

\item  Gauge transformation (multiplication by $e^{i\mu h^{-1}S(x)}$ with quadratic form $S(x)$), transforming $V_j(x)$ into $0$ for $j=1,\ldots,q+r$ and $V_{j+r+q}(x)$ into $f_j x_{j+q}$ for $j=1,\ldots,r$. Then $A^0$ is transformed into
\begin{equation}
\sum_{1\le k\le q} h^2 D_{k}^2+ h^2 \sum_{1\le j\le r} \Bigl( D_{q+j}^2 + \bigl(D_{q+j+r}- f_j\mu h^{-1}x_{q+j}\bigr)^2\Bigr)
\label{20-3-2}
\end{equation}
and $V(x)$ is preserved.

\item
[(iii)--(iv)] Partial $\mu^{-1}h$-Fourier transform, change of variables like in subsection~\ref{book_new-sect-19-3-1} of \cite{futurebook}  $A^0$ into
\begin{equation*}
\sum_{1\le k\le q} h^2 D_{k}^2 + \sum_{1\le j\le r} \Bigl(h^2 D_j^2 +  f_j^2\mu^2x_j^2\Bigr)
\end{equation*}
and $V$ into 
$\tilde{V}\Def V(x''',x'-Kx'', \mu ^{-1}h D''- \,^t\!K \mu ^{-1}hD')$ with Weyl quantization. Finally, $x_j\mapsto f_j^{\frac{1}{2}}x_j$, 
$D_j\mapsto f_j^{-\frac{1}{2}}D_j$ reduces operator $A^0$ into (\ref{20-3-1}) 
and transforms $\tilde{V}$ accordingly.
\end{enumerate}

\section{General case: framework}
\label{sect-20-3-2}

We can assume now that
\begin{equation}
\bar{\mu}_{(q)}=C(h|\log h|)^{-\frac{q}{q+2}} \le \mu \le 
\epsilon (h|\log h|)^{-1}
\label{20-3-3}
\end{equation}
(a case of larger $\mu $ we'll consider later; we also can redefine $\bar{\mu}_{(q)}=Ch^{-\frac{q}{q+2}}$).

Anyway we need to consider an inner zone 
\begin{equation}
\cZ_\inn\Def \{|\xi '''|\le \bar{\rho}_1=C(\mu h|\log h|)^{\frac{1}{2}}\}
\label{20-3-4}
\end{equation}
and in this zone in the general case we already set 
$\varepsilon = C(\mu^{-1}h|\log h|)^{\frac{1}{2}}$. However we reset it to a far larger value $\varepsilon = C\mu^{-1}$ thus defining 
\begin{equation}
\varepsilon = \left\{\begin{aligned}
&C\mu ^{-1} \bigl(\frac {\bar{\rho}_1}{\rho}  \bigr)^s + 
C\rho ^{-1}h|\log h|\qquad &&\text{as\ \ } \rho\ge \bar{\rho}_1,\\
&C\mu^{-1} &&\text{as\ \ } \rho\le \bar{\rho}_1
\end{aligned}\right.
\label{20-3-5}
\end{equation}
with a possible increase later.

From this section point of view there are two different cases $q=1$ and 
$q\ge 2$; if \underline{either} microhyperbolicity condition is fulfilled \underline{or} $g^{jk}$, $F_{jk}$  are constant and non-degeneracy condition (\ref{20-2-42}) is fulfilled, then $\bar{\mu}_{(1)}$ is pushed up and $q=1$ falls in the latter case as well.  

Really, for $q\ge 2$, $\mu \ge \bar{\mu}_{(2)}= (h|\log h|)^{-\frac{1}{2}}$ (and thus $\mu^{-l}|\log h|^{-\sigma}\le \mu h$ even as $(l,\sigma)=(1,1)$) this will bring an approximation error not exceeding
\begin{equation*}
C \mu^r h^{-r} (\mu h|\log h|)^{\frac{q}{2}-1} \mu^{-l}|\log h|^{-\sigma}\times (\mu h)^{1-r} |\log h| 
\end{equation*}
where the first factor estimates an error in each term in the sum defining $h^{-d}\cN^\MW$ and the second factor estimates the number of the terms effected; the result is $O\bigl(h^{1-d}\bigr)$ for sure.

As  $q=1$ we need to assume that\footnote{\label{foot-20-9} This is not necessarily true as $q=1$; however without this assumption no improvement of the previous section results is possible.} 
\begin{multline}
\bar{\rho}_2\Def\bigl(\mu^{-l}|\log h|^{-\sigma}\bigr)^{\frac{1}{2} }\le \
\rho_1^*\Def (\mu h)^{\frac{1}{2}} \\
\iff \mu \ge (h|\log h|^{\sigma})^{1/(l+1)}) .
\label{20-3-6}
\end{multline}
Consider strip 
$\{\rho \le |\xi '''|\le 2\rho \}$ with 
$\bar{\rho}_2\le \rho \le \bar{\rho}_1$ where
an approximation error can be estimated similarly by
\begin{equation*}
Ch^{-d}\rho \bar{\rho}_2^2+C\mu h^{1-d} \rho^{-1}\bar{\rho}_2^2
\end{equation*}
and after summation over $\rho$ we conclude that the contribution of of zone 
$\{\bar{\rho}_2\le |\xi'''|\le \bar{\rho}_1\}$ to an approximation error does not exceed
\begin{equation}
C h^{-d} (\mu h|\log h|)^{\frac{1}{2}}\mu^{-l}|\log h|^{-\sigma} +
C h^{1-d} \mu^{1- \frac{l}{2}}|\log h|^{-\frac{\sigma}{2}}.
\label{20-3-7}
\end{equation}
Meanwhile, in the zone $\{|\xi '''|\le C\bar{\rho}_2 \}$ no more than 
$C(\mu h)^{1-r}$ terms in the sum defining $h^{-d}\cN^\MW$ are affected and an error in each does not exceed $C\mu^r h^{-r} \bar{\rho}_2$; this brings an approximation error estimate $C\mu h^{1-d} \bar{\rho}_2$ which is the exactly the second term in (\ref{20-3-7}).

\begin{remark}\label{rem-2-3-1}
Expression (\ref{20-3-7}) is $O\bigl(h^{1-d}\bigr)$ for all 
$\mu \ge C(h|\log h|)^{-\frac{1}{3}}$ iff $(l,\sigma)\succeq (2,1)$.
\end{remark}

On the other hand, if $q=1$ and microhyperbolicity condition \ref{def-20-1-2} is fulfilled, then an  the approximation error does not exceed
\begin{equation}
C h^{-d} (\mu h|\log h|)^{\frac{1}{2}}\mu^{-l}|\log h|^{-\sigma}
\label{20-3-8}
\end{equation}
and it should be weighted against theorem~\ref{thm-20-2-15}; so in this case for
$\mu \le \bar{\mu}_{(2)}= C(h|\log h|)^{-\frac{1}{2}}$ we use theorem \ref{thm-20-2-15} and for $\mu \ge \bar{\mu}_{(2)}$ we will use the intermediate magnetic field construction; an approximation error estimate in both cases is $O\bigl(h^{1-d}\bigr)$ for all corresponding $\mu$ iff 
$(l,\sigma)\succeq (\frac{3}{2},1)$.

Similar arguments work if $q=1$, $g^{jk}$, $F_{jk}$  are constant and non-degeneracy condition (\ref{20-2-42}) is fulfilled; we refer in this case to problem~\ref{Problem-20-2-17} as $\mu \le \bar{\mu}_{(2)}$.

Thus we have proven

\begin{proposition}\label{prop-20-3-2} 
Let $\varepsilon$ be defined by \textup{(\ref{20-3-5})}. Then

\begin{enumerate}[label=(\roman*), fullwidth]
\item 
If $q\ge 2$, $(l,\sigma)\succeq (1,1)$ and
$\mu \ge (h|\log h|)^{-\frac{1}{2}}$, then a mollification error is $O(h^{1-d})$;

\item
If $q=1$, $(l,\sigma)\succeq (1,2)$, $\mu \ge (h|\log h|)^{-\frac{1}{3}}$ and condition \textup{(\ref{20-3-6})} is fulfilled, then a mollification error does not exceed \textup{(\ref{20-3-7})};

\item
If $q=1$, $(l,\sigma)\succeq (1,2)$, $\mu \ge (h|\log h|)^{-\frac{1}{2}}$ and \underline{either} microhyperbolicity condition  is fulfilled \underline{or} $g^{jk}$, $F_{jk}$  are constant and non-degeneracy condition \textup{(\ref{20-2-42})} is fulfilled, then  a mollification error does not exceed \textup{(\ref{20-3-8})}.
\end{enumerate}
\end{proposition}

\section{Reduction. Main part}
\label{sect-20-3-3}

So we can start from the operator
\begin{gather}
A=\sum _{1\le j,k\le d} P_j^* g^{jk}P_k +V, \qquad P_j= hD_j-\mu V_j \label{20-3-9}\\
\shortintertext{with}
[P_j,P_k]=i\mu hF_{jk} \label{20-3-10}\\
\shortintertext{where}
F_{jk}=0\quad \text{as\ \ }j\le q\  \text{or\ \ }\ k\le q;\;
\rank (F_{jk})_{j,k=1,\dots,2r}=2r,\label{20-3-11}
\end{gather}
and we can assume without any loss of the generality at some point $\bar{x}$ matrix $(F_{jk})$ is in the canonical form; namely 
\begin{equation*}
F_{jk}(\bar{x})=\left\{ \begin{aligned}
f_i\qquad &\text{as\ \ } j=q+i, k=q+r+i,\ i=1,\ldots, r,\\
-f_i\qquad &\text{as\ \ } k=q+i, j=q+r+i,\ i=1,\ldots, r,\\
0 \qquad &\text{otherwise.}
\end{aligned}\right.
\end{equation*}
Then without any loss of the generality we can assume that in its vicinity 
\begin{equation}
F_{jk}=0\ne 0 \implies j\ge q+1,\ k\ge q+1\quad\text{and\ \ }
\exists \fm \in \fM : j,k \in \fm .
\label{20-3-12}
\end{equation}

\begin{proposition}\label{prop-20-3-3} 
Let 
\begin{equation}
\cL = \mu^2\bigl(\sum _{j,k\le d}\beta ^{jk} p_jp_k\bigr)^\w
\label{20-3-13}
\end{equation}
with $\beta^{jk}=\beta^{kj}\in \sF^{\bar{l},\bar{\sigma}}$. Then modulo operators with symbols belonging to $\mu^{-2}\sF^{\bar{l}-1,\bar{\sigma}}$
\begin{equation}
e^{\frac{i}{2} \mu^{-1}h^{-1}t\cL}(p_j)^\w e^{-\frac{i}{2}  \mu^{-1}h^{-1}t\cL}
\equiv \Bigl(\sum_k \bigl(e^{\Lambda t}\bigr)_j^k p_k\Bigr)^\w
\label{20-3-14}
\end{equation}
with matrix $\Lambda= (\Lambda _{jk})_{j,k}=(F_{jp})(\beta^{pk})$.
\end{proposition}

\begin{proof}
An easy proof  following arguments of Sections~\ref{book_new-sect-18-7} and~\ref{book_new-sect-19-3} of \cite{futurebook} is left to the reader.
\end{proof}

One can see easily that due to (\ref{20-1-7})--(\ref{20-1-8}) 
\begin{equation*}
e^{\Lambda t} = \begin{pmatrix} 
I& 0\\[2pt]
\int_0^t e^{\Lambda 't'} \Lambda''\, dt' &e^{\Lambda 't}
\end{pmatrix}
\end{equation*}
with $\Lambda ' = (\Lambda _{jk})_{j,k\ge q+1}$, 
$\Lambda'' = (\Lambda _{jk})_{j\le q+1, k\le q}$. It provides us with the first step of the reduction:

\medskip\noindent
\emph{Step 1}. 
Applying transformation $\exp (\frac{i}{2} \mu^{-1}h^{-1}\cL)$ with operator $\cL$ defined by (\ref{20-3-13}) with appropriate coefficients $\beta^{jk}$ one can transform operator $A$ to the block form
\begin{gather}
\mu^2 \sum _{q+1\le j,k\le q+r} \zeta_j^{\dag\,\w}a^{jk\,\w}\zeta_k^\w + 
\mu ^2\sum _{j,k\le q} p_j^\w a^{jk\,\w}p_k^\w + a_0^\w+ 
\frac{1}{2} \bigl(B+B^*\bigr),\label{20-3-15}\\ 
B= \Bigl(\mu^2\sum_{j,k,m} b^{jkm}p_jp _kp _m + \sum_j b^jp_j\Bigr)^\w \label{20-3-16} \\
\zeta_j = p_j+ ip_{j+r} \; \text{with\ \ }j=q+1,\ldots,q+r, \label{20-3-17}\\[2pt]
a^{jk}\ne 0 \implies \exists \fm : j,k\in \fm; \label{20-3-18}
\end{gather}
let us recall that $\fM $ is the resonance partition at point $\bar{x}$. 

\begin{claim}\label{20-3-19}
Here and below $a^{jk}\in \sF^{\bar{l},\bar{\sigma}}$, $a_0\in\sF^{l,\sigma}$, $b^{jkl}\in \sF^{\bar{l}-1,\bar{\sigma}}$, $b^j\in \sF^{l-1,\sigma}$
and while $a^{jk}=a^{jk}(x)$ we have $b^{jkl}=b^{jkl}(x,\mu ^{-1}\xi)$,
$b^j = b^j(x,\mu^{-1}\xi)$ and $c^\w$ means $h$- rather than $\mu^{-1}h$-quantization of symbol $c$. 
\end{claim}
\emph{Step 2}. Applying an appropriate gauge transformation we can achieve 

\begin{claim}\label{20-3-20}
$V_k=0$ as $k\le q+r$ and
$V_k=V_k(x_{q+1},\ldots,x_{d})$ for $k=q+r+1,\ldots,d$
\end{claim}
where condition \ref{20-1-9-2} provides  $\sF^{\bar{l}+1,\bar{\sigma}}$-regularity of new $V_k$. This is the only reason for this condition. Then 
\begin{multline}
p_j=0\ \ \forall j=q+r+1,\dots, d \iff \\
x_j =\lambda_j(x'',\xi')\ \ \forall j=q+1,\ldots,q+r
\label{20-3-21}
\end{multline}
with $x'=(x_{q+1},\ldots, x_{q+r})$, $\xi'=(\xi_{q+1},\ldots, \xi_{q+r})$, 
$x''=(x_{q+r+1},\ldots, x_{q+2r})$, $\xi''=(\xi_{q+r+1},\dots, \xi_{q+2r})$, and 
$\lambda_j\in \sF^{\bar{l}+1,\bar{\sigma}}$.

After this, following arguments of Sections~\ref{book_new-sect-18-3},~\ref{book_new-sect-18-7}, and~\ref{book_new-sect-19-3} of \cite{futurebook} and applying transformation  $e^{i\mu^{-1} \cL'}$ with 
$\cL'=\sum _{q+1\le j\le q+r}\lambda_j ^\w D_{j+r}$ we can transform operator (\ref{20-3-15}) to the form (\ref{20-3-15})--(\ref{20-3-16}) but now with
\begin{align}
\zeta_j = \eta_j+ i\eta_{j+r},\; &\eta_j= \mu^{-1}\xi_j,\; \eta_{j+r}=x_j\qquad &&j=q+1,\ldots, q+r, \label{20-3-22}\\[2pt]
&\eta_j= \mu^{-1}\xi_j \quad &&j=1,\ldots, q.
\label{20-3-23}
\end{align}
However in contrast to the previous step of the reduction now
\begin{equation*}
a^{jk}=a^{jk}(x'',\mu^{-1}\xi''; x'''),\quad a_0=a_0(x'',\mu^{-1}\xi''; x''')
\end{equation*}
rather than $a^{jk}=a^{jk}(x)$, $a_0=a_0(x)$.

\begin{remark}\label{rem-20-3-4}
\begin{enumerate}[label=(\roman*), fullwidth]
\item  One can see easily that if $f_j$ have constant multiplicities, then the operator can be reduced on Step 1 to the form
\begin{equation}
\mu^2 \sum _{q+1\le j\le q+r} \zeta_j^{\dag\,\w}f_j^\w\zeta_j^\w + 
\mu ^2\sum _{j,k\le q} p_j^\w a^{jk\,\w}p_k^\w + a_0^\w+ 
\frac{1}{2} \bigl(B+B^*\bigr)
\label{20-3-24}
\end{equation}
and this form will be preserved after Step 2.

\item  One can see easily that if $q=1$, then after Step 2 by means of transformation $e^{ih^{-1}\cL''}$ with 
$\cL''= \frac{1}{2}(hD_1 k^\w + k^\w hD_1)$, $k=k(x'',x_d, \mu^{-1}\xi'')$ one can reduce operator to the same form (\ref{20-3-15}) but with $a^{11}=1$.

\item 
In the general case we can achieve $a^{jk}=\updelta_{jk}$ for $j,k=1,\ldots, q$ only assuming that $g^{jk}$, $F_{jk}$ are constant.
\end{enumerate}
\end{remark}

\section{Reduction. Junior terms}
\label{sect-20-3-}

Recall that junior terms are of the form $\frac{1}{2}(B+B^*)$ with $B$ defined
by (\ref{20-3-16}) with $b^{jkm}=b^{jkm}(x,\mu^{-1}\xi)$ and $b^j=b^j(x,\mu^{-1}\xi)$. 

\begin{remark}\label{rem-20-3-5}
As $(1,1)\preceq (l,\sigma)\preceq (2,0)$ one can rewrite $B$ in the same form but with 
\begin{equation}
b^{jkm}=b^{jkm}(x'',\mu^{-1}\xi''; x'''),\qquad b^j=b^j(x'',\mu^{-1}\xi''; x''')
\label{20-3-25}
\end{equation}
modulo operator $B'$ with the symbol belonging to 
$\mu^{-l}|\log \mu|^{-\sigma}\sF ^{0,0}$. Our reduction is done modulo operators of this type.
\end{remark}

Then applying operator $e^{i\mu^{-2}h^{-1}\cL_1}$ with
\begin{equation}
\cL_1 = \Bigl(\sum _{\alpha: 1\le |\alpha|\le 3}
\beta ^{(\alpha)} \mu ^{|\alpha|} \eta^\alpha \Bigr)^\w, 
\label{20-3-26}
\end{equation}
where $\eta^\alpha = \eta_1^{\alpha_1}\cdots \eta_d^{\alpha_d}$, one can reduce operator to the same form (\ref{20-3-15}) with 
\begin{multline}
B\equiv \mu ^2 \Bigl(\sum_{q+1\le j,k,m\le q+r} b^{jkm}
\zeta_j^\dag \zeta_k\zeta_m + \\
\sum_{q+1\le j,k\le q+r,\ 1\le m\le q} b^{jkm}\eta_m\zeta_j^\dag \zeta_k+
\sum_{1\le j,k, m\le q} b^{jkm}\eta_m\eta_j\eta_k+ 
\sum_{1\le j\le q} b^j\eta_j\Bigr)^\w.
\label{20-3-27}
\end{multline}

Here in the first term $b^{jkm}=0$ unless $|f_j-f_k-f_m|\le \epsilon $ and  in the second term $b^{jkm}=0$ unless $|f_j-f_k|\le \epsilon$ (so there are the $3$-rd and $2$-nd order resonances respectively). Further, $b^{jkm}=0$ always for constant $g^{jk}, F_{jk}$.

Also, in (\ref{20-3-26}) $\beta^{(\alpha)}$ belongs to $\sF^{\bar{l}-1,\bar{\sigma}}$, $\sF^{l-1,\sigma}$ as $|\alpha |\ge 2$, 
$|\alpha |=1$ respectively;  also 
$b^{jkm}\in \sF^{\bar{l}-1,\bar{\sigma}}$, $b^j\in \sF^{l-1,\sigma}$.

\begin{remark}\label{rem-20-3-6} 
\begin{enumerate}[label=(\roman*), fullwidth]
\item  In the intermediate magnetic field case $|\eta_j|\le C\mu^{-1}$ in the microlocal sense (on the energy levels below $c$) and therefore junior terms are of magnitude $\mu^{-1}$ and skipping them we would arrive to an approximation error $Ch^{-d}(\mu h|\log h|)^{\frac q 2}\mu^{-1}$ (as either $q\ge 2$, or $q=1$ and microhyperbolicity condition is fulfilled, or $(l,\sigma)\succeq (2,0)$,  and non-degeneracy condition is fulfilled) which is $O(h^{1-d})$ as $q\ge 3$ and $O(h^{1-d}|\log h|)$ as $q=2$. 

Thus, in the case of an intermediate magnetic field only for $q=1,2$ we will need to take junior terms into account at all. This will be a minor annoyance for $q=2$ but a major obstacle for $q=1$, when without either microhyperbolicity or non-degeneracy assumptions this error would be 
$O\bigl(h^{-d}(\mu h|\log h|)^{\frac{1}{2}}\mu^{-1}+ 
\mu^{\frac{1}{2}}h^{1-d}\bigr)$
but even under the microhyperbolicity or non-degeneracy assumptions these terms cannot be simply ignored.

\item  Actually, not all of these terms are equally bad: terms containing factor(s) $\eta_m$ with $1\le m \le q$ are of magnitude $O(\mu^{-1}|\xi'''|)$ and can be simply ignored as $q=2$ but even as $q=1$ they are minor annoyances. Only the first term in (\ref{20-3-27}) are really important. 

\item  One can see easily that if $q=1$ and $f_j$ have constant multiplicities, then by means of transformation $e^{i\mu^{-1}h^{-1}\cL'_1}$ with 
\begin{equation*}
\cL'_1 = \mu^2 \Bigl(\sum_\fm  \sum_{j,k\in \fm } 
c^{jk} \zeta_j^{\dag}\zeta_k + c_0 \eta_1^2 +c_1\Bigr)^\w 
\end{equation*}
with $c^{jk}=c^{jk}(x'',x_d,\mu^{-1}\xi'')$, $c_j=c_j(x'',x_d,\mu^{-1}\xi'')$ one can remove from $B$ all terms containing factor $\eta_1$ thus arriving to 
\begin{equation}
B\equiv \mu ^2 \Bigl(\sum_{q+1\le j,k,m\le q+r} b^{jkm}
\zeta_j^\dag \zeta_k\zeta_m\Bigr)^\w.
\label{20-3-28}
\end{equation}
\end{enumerate}
\end{remark}

Therefore, we constructed operator $\cT$ with the properties described below:

\begin{proposition}\label{prop-20-3-7} 
Let $(\bar{l},\bar{\sigma})=(2,1)$, $(1,1)\preceq (l,\sigma)\preceq (2,0)$. Then there exists an operator $\cT$ such that
\begin{align}
&\cT^* \cT \bar{Q}\bar{\psi} \equiv \bar{Q}{\bar\psi},\label{20-3-29}\\[2pt]
&\cT^* Q_1\cT (I-\bar{Q})\bar{\psi} \equiv 0,\label{20-3-30}\\[2pt]
&\cT^* (I-Q_2)\cT \bar{Q}\bar{\psi} \equiv 0&& \mod \mu^{-\infty} (\sF^{0,0})^\w,\label{20-3-31}\\[2pt]
&\cT^* \cA \cT \bar{Q}\bar{\psi} \simeq A\bar{Q}\bar{\psi},
&& \mod \mu^{-l}|\log h|^{-\sigma} (\sF^{0,0})^\w \label{20-3-32}
\end{align}
as $\bar{\psi}\in \sC_0^\infty (B(0,1))$, $\bar{Q}=\bar{Q}(hD''')$ has a symbol which is supported in $\{|\xi'''|\le 2C\bar{\rho}_1\}$ and equal to $1$ in 
$\{|\xi'''|\le 2C\bar{\rho}_1\}$, $Q_1$ and $Q_2$ are the same type operators but with symbols which are supported in $\{|\xi'''|\le 2\bar{\rho}_1\}$ and 
$\{|\xi'''|\le 2C^2\bar{\rho}_1\}$ and equal to $1$ in 
$\{|\xi'''|\le \bar{\rho}_1\}$ and $\{|\xi'''|\le C^2\bar{\rho}_1\}$ respectively with large enough constant $C$, $\cA$ is the reduced operator in the form \textup{(\ref{20-3-15})} with junior terms defined by \textup{(\ref{20-3-27})}. 

Further, as $q=1$ and $f_j$ have constant multiplicities, junior are terms defined by \textup{(\ref{20-3-28})}. 
\end{proposition}

\begin{remark}\label{rem-20-3-8}
\begin{enumerate}[label=(\roman*), fullwidth]
\item  Note that while (\ref{20-3-29})--(\ref{20-3-31}) hold modulo negligible operator, (\ref{20-3-32}) is satisfied modulo $\mu^{-1}h$-pseudo-differential operators with symbols of 
$\mu^{-l}|\log h|^{-\sigma} \sF^{0,0}$, and to emphasize the difference we use ``$\simeq$'' instead of ``$\equiv$''. 

\item
In contrast to almost all this Chapter in propositions \ref{prop-20-3-7}, \ref{prop-20-3-9} $b^\w$ means $\mu^{-1}h$-quantization of symbol $b$.
\end{enumerate}
\end{remark}

We also need the following properties of  operator constructed above:

\begin{proposition}\label{prop-20-3-9} 
Furthermore the constructed operator has the following properties: 
\begin{enumerate}[label=(\roman*), fullwidth]
\item For functions $ \bar{\psi}=1$ in $B(0,\frac{3}{4})$ and 
$\psi\in \sC_0^\infty (B(0,\frac{1}{2})$ 
\begin{align}
& \psi\bar{Q}\bar{\psi} \equiv \cT^* \tilde{\psi}^\w\cT \bar{Q} 
\bar{\psi},\label{20-3-33}\\
&\tilde{\psi}\simeq \psi _0(x'',x''',\xi '') +
\sum _{1\le j\le d} \psi _j (x'',x''',\xi '')\eta_j \label{20-3-34}\\
&\qquad\qquad\qquad\qquad\qquad\qquad\quad \mod \mu^{-l}|\log h|^{-\sigma}\sF^{0,0};\notag
\end{align}
\item For operator $Q=Q^\w(hD''')$ as above 
\begin{align}
&\cT^*Q \cT \bar{\psi}\simeq \bigl(Q(f(x,\xi))\bigr)^\w \bar{\psi}\qquad &&\mod \mu^{-l} |\log h|^{-\sigma}(\sF^{0,0})^\w ,\label{20-3-35}\\
&f(x,\xi)\simeq \xi''' -\sum_{1\le j\le 2r} \kappa_j\eta_j -\mu^{-1}\kappa_0&&\mod \mu^{-l}|\log h|^{-\sigma}\sF^{0,0}\label{20-3-36}
\end{align}
with coefficients $\kappa_j=\kappa_j(x'',x''',\mu^{-1} \xi'')$ and all estimates holding for $|\eta_j|\le C\mu^{-1}$, $\psi_j,\alpha_j \in \sF^{l-1,\sigma}$.
\end{enumerate}
\end{proposition}

\begin{proof}
We prove more difficult statement (ii) leaving easier (i) to the reader.  

\begin{enumerate}[label=(\roman*), fullwidth]
\item On Step 1 of the main part reduction we considered operator 
$e^{\frac{i}{2} \mu^{-1}h^{-1}t\cL}$
with the corresponding symplectomorphism which is equivalent to the map
$f_t: (x,\eta_1,\dots,\eta_d)\to (x,\eta_1,\dots,\eta_d)$ which is of $\sF^{\bar{l},\bar{\sigma}}$ class. Let $Q_t= \bigl(Q(f_t)\bigr)^\w$; then the principal symbol of operator
\begin{equation}
\frac{\partial\ }{\partial t}Q_t + 
\frac{i}{2} \mu ^{-1}h^{-1}\bigl[ \cL, Q_t\bigr]
\label{20-3-37}
\end{equation}
vanishes and then as long as $(\bar{l},\bar{\sigma})\succeq (2,1)$ and 
$\varepsilon \ge C(\mu^{-1}h|\log h|)^{\frac{1}{2}}$ the norm of this operator does not exceed
\begin{equation*}
C\mu^{-1}h^{-1} \times \bigl( {\frac h \rho} \bigr)^2 \le 
C \mu^{-1} h (\mu h|\log h)^{-1} = C\mu^{-2}|\log h|^{-1}
\end{equation*}
for $\rho \ge C(\mu h|\log h|)^{\frac{1}{2}}$. Then the norm of the operator
\begin{equation}
\frac{d\ }{dt} 
\Bigl(e^{-\frac{i}{2} \mu^{-1}h^{-1}t\cL}Q_t e^{\frac{i}{2} \mu^{-1}h^{-1}t\cL}
\Bigr) 
\label{20-3-38}
\end{equation}
does not exceed $C\mu^{-2}|\log h|^{-1}$ and then the same is true for operators
\begin{equation}
e^{-\frac{i}{2} \mu^{-1}h^{-1}t\cL}Q_te^{\frac{i}{2} \mu^{-1}h^{-1}t\cL} - Q \quad\text{and}\quad
e^{\frac{i}{2} \mu^{-1}h^{-1}t\cL}Qe^{-\frac{i}{2} \mu^{-1}h^{-1}t\cL} - Q_t.
\label{20-3-39}
\end{equation}
\item
The same arguments work perfectly on Step 2 of the main part reduction.

\item
Now consider the junior terms reduction. First of all, there are removable terms 
$\sum_{q+1\le j\le q+2r} \beta_j \eta_j$ with 
$\beta_j=\beta_j(x'',\mu^{-1}\xi'')$, $\beta_j\in \sF^{l-1,\sigma}$ which are removed by a transformation of the same type as above with 
$\cL= \mu^{-1}\sum_{1\le j\le 2r} \beta'_j \eta_j$ and $\beta'_j=\beta'_j(x'',\mu^{1}\xi'')$, $\beta'_j\in \sF^{l-1,\sigma}$ as well.

Then the corresponding Hamiltonian map is of class $\sF^{l,\sigma}$ and then the principal symbol of operator (\ref{20-3-34}) vanishes and the norm of this operator does not exceed 
\begin{multline*}
C\mu^{2-l}|\log h|^{-\sigma}\times\mu^{-1}h^{-1}\times 
\Bigl(\bigl( \frac{h}{\rho}\bigr)^2 + 
\mu^{-1}h \bigl(\frac{h}{\rho}\bigr) \times \mu^{2-l}|\log h|^{-\sigma}\Bigr)\le \\
C\mu^{-l}|\log h|^{-\sigma}\times \frac {\mu h}{\rho^2}
\end{multline*}
which is less than $C\mu^{-l}|\log h|^{-\sigma}$ and then the same estimate holds for operators (\ref{20-3-38})--(\ref{20-3-39}).

With the removable third and second order terms the construction is the same but simpler since these terms are coming from the principal part which belongs to $\sF^{2,1}$ for sure.
\end{enumerate}
\end{proof}

\begin{remark}\label{rem-20-3-10}  
In Chapter~\ref{book_new-sect-18} in the case $d=3$ we took $\varepsilon = C(\mu^{-1}h|\log h|)^{\frac{1}{2}}$ instead of $C\mu^{-1}$. The above construction would not work. In this case note first that the rescaling 
$x \mapsto \varepsilon^{-1}(x-\bar{x})$ and  taking 
$\hbar = \mu^{-1}h\varepsilon^{-2}= C^{-1}|\log h|^{-1}$ we will
see that the original operator $\bar{\chi}(hD_1/\rho )$ becomes 
$\bar{\chi}(\hbar D_1)$ and $\mu^{-1}h^{-1}\cL$ becomes 
$\nu \ell (x,\hbar D)$ with  $\nu =\mu^{-2}\varepsilon^{-1} =
C^{-1}(\mu^3 h|\log h|)^{-\frac{1}{2}}\le \epsilon_1$ and therefore $e^{i\mu^{-1}h^{-1}\cL}$ becomes $\hbar$-pseudo-differential operator with ``analytic'' symbol (more precisely its symbol and those of $Q$ satisfy assumptions of Subsection~\ref{book_new-sect-1-1-3} of \cite{futurebook}). 

Then $e^{-i\mu^{-1}h^{-1}\cL}Qe^{i\mu^{-1}h^{-1}\cL}$ becomes $\hbar$-pseudo-differential operator which means that in the original scale it is $\mu^{-1}h$-pseudo-differential operator with $\sF^{0,0}$
symbol supported in $\{|\xi_3|\le 2C\rho\}$ and equal $1$ in 
$\{|\xi_3|\le C\rho\}$ with $\rho =\bar{\rho}_1=C(\mu h|\log h|)^{\frac{1}{2}}$.
\end{remark}

\section{Strong and Stronger Magnetic Field}
\label{sect-20-3-5}

In this case we have 
\begin{gather}
\mu\ge \epsilon (h|\log h|)^{-1} \label{20-3-40}\\
\shortintertext{and}
\varepsilon = C(\mu^{-1}h|\log \mu |)^{\frac{1}{2}}.\label{20-3-41}
\end{gather}
Then construction of Subsections~\ref{book_new-sect-3-3}--\ref{book_new-sect-3-4} of \cite{futurebook} works on every its step. Now in a microlocal sense  
$|\eta_j|\le C(\mu^{-1}h|\log h|)^{\frac{1}{2}}$ (on the energy levels below $c$),  and therefore junior terms seem to be of magnitude 
$ C\mu^2 (\mu^{-1}h|\log \mu|)^{\frac{3}{2}}$ and cannot be skipped; even next junior terms are of magnitude $C\mu^2 (\mu^{-1}h|\log \mu|)^2$ which is
$O(h)$ as $h|\log \mu|^2\le C$ only. 

In fact, however, as we will see, junior terms will have no impact because the microlocal estimate will be replaced by an operator one
$|\eta_j^\w|\le C(\mu^{-1}h)^{\frac{1}{2}}$ which is sufficient to estimate junior terms by $Ch$ as long as $\mu \le C_0h^{-1}$ (strong and very strong magnetic field). 

Further,  for $\mu h\ge C_0h^{-1}$ (\emph{superstrong\/} magnetic field, we consider Schr\"odinger-Pauli operator then)
only the lowest Landau level is important where irreducible terms due to the $3$-rd order resonances vanish and irreducible terms due to the $4$-th order resonances  are of magnitude $C\mu^2(\mu^{-1}h)^2=O(h^2)$ in the operator sense.

\chapter{Intermediate magnetic field: estimates}
\label{sect-20-4}

In this section we assume that the magnetic field is intermediate i.e.  condition (\ref{20-3-3}) is fulfilled. In this case we know that the contribution of the outer zone 
$\Omega_\out\Def \{|\xi'''|\ge \bar{\rho}_1\Def 
C(\mu h|\log h|)^{\frac{1}{2}}\}$ to the Tauberian estimate with 
$T=\epsilon \mu ^{-1}$ is $O(h^{-2})$ and therefore we need to consider a complementary zone 
$\Omega_\out^c =\{|\xi'''|\le \bar{\rho}_1\Def C(\mu h|\log h|)^{\frac{1}{2}}\} = \Omega_\interm \cup\Omega_\inn $ with an \emph{intermediate zone\/}\index{zone!intermediate}
\begin{equation}
\Omega_\interm \Def
\bigl\{\bar{\rho}^*_1\Def C\max\bigl(\mu^{-1}, (\mu h)^{\frac{1}{2}}\bigr) \le |\xi'''|\le \bar{\rho}_1=(\mu h|\log h|)^{\frac{1}{2}}\bigr\}
\label{20-4-1}
\end{equation}
and an inner zone $\Omega_\inn = \{|\xi'''|\le \bar{\rho}^*_1\}$.

\section{Intermediate zone}
\label{sect-20-4-1}

\subsection{General arguments}
\label{sect-20-4-1-1}

Let us consider an intermediate zone (\ref{20-4-1}). We split it into strips $\cZ_\rho$ with $\rho$ in the indicated frames:
\begin{equation}
\cZ_\rho\Def \{(x,\xi):\ \rho \le |\xi'''|\le 2\rho\},\qquad 
\bar{\rho}^*_1\le \rho \le \bar{\rho}_1.
\label{20-4-2}
\end{equation}

\begin{proposition}\label{prop-20-4-1} 
Let condition \textup{(\ref{20-3-3})} be fulfilled and let $Q=Q^\w$ be an operator with the symbol supported in the strip $\cZ_\rho$, and such that
\begin{equation}
|\partial_{x,\xi}^\alpha Q|\le C_\alpha \rho^{-|\alpha|}\qquad \forall \alpha,
\label{20-4-3}
\end{equation}
and let $\psi \in \sC_0^\infty (B(0,1))$. 

Let $\chi \in \sC_0^\infty ([-1,-\frac{1}{2}]\cup [\frac{1}{2}, 1])$. Then estimate
\begin{gather}
|F_{t\to h^{-1}\tau} \chi_T(t)\Gamma Q\psi U |\le 
C h^{-d} T \rho^q \bigl(\rho^2 + \rho T+\mu^{-1}\bigr)
\bigl( \frac {h} {\rho^2 T} \bigr)^s
\label{20-4-4}\\
\shortintertext{holds as}
T_*=\frac{h} {\rho^2} \le T\le T^*=\epsilon \rho
\label{20-4-5}
\end{gather}
with an arbitrarily large exponent $s$.
\end{proposition}

\begin{proof} Proof repeats those of analysis in Section~\ref{book_new-sect-18-8} of \cite{futurebook}. 
First of all, let us make $\gamma$-partition with respect to $(x'',\xi'',x''')$ with 
\begin{equation}
\gamma \ge C\mu^{-1}+ C(\mu^{-1}h|\log h|)^{\frac{1}{2}}
\label{20-4-6}
\end{equation}
and consider some element $\psi'(x)$ of this partition supported in $B(\bar{x},\gamma)$. Let us reduce operator to its canonical form in $B(\bar{x},2\gamma)$; without any loss of the generality one can assume that this canonical form is
\begin{equation}
\sum_{1\le j\le r} \bar{f}_j (h^2D_{q+j}^2 + \mu^2 x_{q+j}^2) + 
\sum_{1\le j \le q}h^2\bar{g}^{jk}D_j D_k+ \bar{V}+O(\gamma)
\label{20-4-7}
\end{equation}
there where $\bar{f}_j$, $\bar{g}^{jk}$ and $\bar{V}$ are constant (and we can assume that $\bar{g}^{jk}=\updelta_{jk}$, $O(\gamma)$ covers variations of $f_j$, $V$ as well as unremovable terms and we are discussing only 
$B(\bar{x},2\gamma)$.

Let us decompose $\sU \Def \cT U$ into the sum
\begin{equation}
(\cT U) (x,y,t)= \sum_{\alpha,\beta\in \bZ^{+\,r }} \sU_{\alpha,\beta}(x'',x''',y'',y''',t)
\Upsilon_\alpha (x')\Upsilon_\beta (y')
\label{20-4-8}
\end{equation}
with $\Upsilon_\alpha (x')=
\hslash^{-\frac{r}{2}}\prod_j \upsilon_{\alpha_j}(x_{q+j}/ \hslash )$, 
$\hslash =(\mu^{-1}h)^{\frac{1}{2}}$ and Hermite functions $\upsilon_*$.

Then modulo $O(\gamma)$ in $\sL^2$-norm operator $\cA$ applied to the terms containing $\Upsilon_\alpha(x')$ in this sum becomes 
\begin{equation}
\bar{\cA}_\alpha = \sum_{1\le j\le r} \bar{f}_j (2\alpha_j+1)\mu h + 
\sum_{1\le j\le q}h^2D_j^2+\bar{V}
\label{20-4-9}
\end{equation}
which is a temporary notation.

Let us strengthen condition (\ref{20-4-6}) assuming that 
\begin{equation}
\gamma \ge \rho^2 + \mu h+\mu^{-1}.
\label{20-4-10}
\end{equation}

Then using our standard methods one can prove easily that for $Q'=Q \cT \psi' \cT^*$ the following inequality holds
\begin{multline}
|F_{t\to h^{-1}\tau} \chi_T(t)\Gamma ' Q'\psi \sU_{\alpha,\beta}|\le \\[3pt]
C\mu^r h^{-r-q} \rho^q \gamma^d T \times
\Bigl(\bigl( \frac{h} {\rho^2 T} \bigr)^s + 
\bigl( \frac{h}{\rho \gamma} \bigr)^s + 
\bigl( \frac {\mu^{-1}h} { \gamma^2} \bigr)^s \Bigr) \times \\[3pt]
\Bigl(\frac{\gamma} 
{\gamma+|\mu h\sum_j \bar{f}_j(\alpha_j -\bar{\alpha}_j)|}\Bigr)^s \times 
\Bigl(\frac{\gamma} {\gamma+|\mu h\sum_j \bar{f}_j(\beta_j-\bar{\alpha}_j)|} \Bigr)^s
\label{20-4-11}
\end{multline}
with some $\bar{\alpha}=\bar{\alpha}(\tau)\in \bZ^{+\, r}$ where 
\begin{equation}
\Gamma 'v\Def \int v|_{x''=y'',x'''=y'''}\,dy''\,dy'''.
\label{20-4-12}
\end{equation}

Here factor $\mu h^{1-d}\rho ^q\gamma^d $ in the right-hand expression comes as a trace norm of $Q'\sU_{\alpha,\beta}$; factor $T$ is a result of integration with respect to $t$ in the Fourier transform,  the second factor is the sum of three terms: 
the first one $\bigl(  h / (\rho^2 T) \bigr)^s$  is due to the microhyperbolicity with respect to $\xi'''$ and rescaling while  terms 
$\bigl( h/ (\rho\gamma)\bigr)^s$ and 
$\bigl( {\mu^{-1}h} /( \gamma^2)\bigr)^s$ are operator calculus errors.

Finally,  two last factors in the right-hand expression of (\ref{20-4-11}) are due to the ellipticity of operators $\bigl(hD_t -\bar{\cA}_\alpha \bigr)$ and  
$\bigl(hD_t -\bar{\cA}_\beta \bigr)$ respectively.

More precisely, let us note that operator $\bigl(\bar{\cA}_\alpha+O(\gamma)-\tau \bigr)$ is elliptic in this vicinity as long as $\alpha \notin \fA$  where
\begin{equation}
\fA\Def \{\alpha\in \bZ^{+\,r}: 
|\sum_j (2\alpha_j+1)\bar{f}_j\mu h +\bar{V}-\tau |\le \gamma\}.
\label{20-4-13}
\end{equation}
Really, we can express there $\sU_{\alpha,*}$ with $\alpha\notin \fA$ via $\sU_{\alpha',*}$ with $\alpha'\in \fA$.  Namely, for $\alpha \notin \fA$ we have
\begin{multline}
|F_{t\mapsto h^{-1}\tau} \chi_T(t) Q' \sU_{\alpha,\beta}|\le\\
\max_{\alpha'\in  \fA } |F_{t\mapsto h^{-1}\tau} \chi_T(t) Q' \sU_{\alpha',\beta}|\times
\Bigl(\gamma 
\bigl(|\sum_j (2\alpha_j+1)\bar{f}_j\mu h +\bar{V} |+\gamma \bigr)^{-1}\Bigr)^s
\label{20-4-14}
\end{multline}
The similar inequality holds for $\sU_{*,\beta}$ with $\beta\notin \fA$.
Finally, for  $\sU_{\alpha,\beta}$ with arbitrary  $\alpha,\beta \notin \fA$ 
\begin{multline}
|F_{t\mapsto h^{-1}\tau} \chi_T(t) Q' U{\alpha,\beta}|\le
\max_{\alpha',\beta'\in  \fA } |F_{t\mapsto h^{-1}\tau} \chi_T(t) Q' \sU_{\alpha',\beta'}|\times\\[3pt] 
\Bigl(
\gamma\bigl(|\sum_j (2\alpha_j+1)\bar{f}_j\mu h+\bar{V} |+\gamma \bigr)^{-1} \Bigr)^s\times \Bigl(
\gamma \bigl(|\sum_j (2\beta_j+1)\bar{f}_j\mu h +\bar{V} |+\gamma \bigr)^{-1} \Bigr)^s
\label{20-4-15}
\end{multline}
which in combination with the microhyperbolicity arguments justifies (\ref{20-4-11}). 

Note also that  $\#\fA \le C\gamma (\mu h)^{-r}$. Then (\ref{20-4-11}) implies
\begin{multline}
|F_{t\to h^{-1}\tau} \chi_T(t)\Gamma Q'\psi U |\le \\
C\mu^r h^{-r-q}  \rho^q \gamma^d T\times 
\Bigl(\bigl( \frac{h} {\rho^2 T} \bigr)^s + 
\bigl(\frac{h} {\rho \gamma} \bigr)^s + 
\bigl( \frac {\mu^{-1}h} { \gamma^2} \bigr)^s\Bigr) \times 
\gamma (\mu h)^{-r}=\\
C h^{-d} T \rho^q \gamma^{d +1}
\Bigl(\bigl( \frac {h} {\rho^2 T}  \bigr)^s + 
\bigl( \frac {h} {\rho \gamma} \bigr)^s + 
\bigl(\frac {\mu^{-1}h} { \gamma^2} \bigr)^s\Bigr) .
\label{20-4-16}
\end{multline}
Restoring to original $Q$, $\psi$ (by summation over $x$-partition) we lose factor $\gamma^d$. As exponent $s$ is large enough, the optimal value of $\gamma$ honoring (\ref{20-4-10}) is  $\gamma= \rho T + \rho^2+\mu ^{-1}$ because for 
$\rho \ge \mu^{-1}$ (which we assume) $\gamma ^2 \ge \mu^{-1}\rho^2 T$.

Then we arrive to estimate (\ref{20-4-4}).
\end{proof}

Proposition \ref{prop-20-4-1} implies immediately estimate
\begin{equation}
|F_{t\to h^{-1}\tau} \bigl(\bar{\chi}_T(t)-
\bar{\chi}_{\bar{T}}(t)\bigr)\Gamma Q\psi U|\le 
C h^{-d} \mu^{-1} \rho^q (\rho^2+\mu^{-1}) \times 
\bigl( \frac {\mu h} {\rho^2 } \bigr)^s 
\label{20-4-17}
\end{equation}
with $ \bar{\chi}\in \sC_0^\infty ([-1,1])$ equal $1$ on 
$[-\frac{1}{2}, \frac{1}{2}]$, $T_*=\epsilon \mu^{-1}$ as $T$ satisfies
\begin{equation}
\bar{T}=\epsilon\mu^{-1} \le T\le T^*=\epsilon \rho.
\label{20-4-18}
\end{equation}

In turn (\ref{20-4-17}) implies 
\begin{equation}
|F_{t\to h^{-1}\tau}\bar{\chi}_T(t)\Gamma Q\psi U|\le 
C h^{-d} \mu^{-1}\rho^q (\rho^2+\mu^{-1}) \times 
\bigl( \frac {\mu h} {\rho^2 } \bigr)^s + Ch^{1-d}\rho^q
\label{20-4-19}
\end{equation}
where $C\rho^q h^{1-d}$ estimates 
$|F_{t\to h^{-1}\tau}\bar{\chi}_{T_*}(t)\Gamma Q\psi U|$ due to the standard results rescaled.

Therefore the standard Tauberian arguments immediately imply estimate
\begin{multline}
\R_{Q}^\T \Def |\Gamma (Q\psi \tilde{e}) - h^{-1} \int_{-\infty}^0 
\Bigl(F_{t\to h^{-1}\tau}\bar{\chi}_T(t)\Gamma Q\psi u\Bigr)\,d\tau |\le\\
CT^{-1}\Bigl(C h^{-d} \mu^{-1}\rho^q (\rho^2+\mu^{-1}) \times 
\bigl( \frac {\mu h} {\rho^2} \bigr)^s + Ch^{1-d}\rho^q\Bigr)=\\
C h^{-d} \mu^{-1} \rho^{q-1} (\rho^2+\mu^{-1})\times 
\bigl( \frac {\mu h} {\rho^2}  \bigr)^s + Ch^{1-d}\rho^{q-1}
\label{20-4-20}
\end{multline}
with any
\begin{equation}
T_*\Def Ch |\log h|/\rho^2 \le T\le \epsilon \rho.
\label{20-4-21}
\end{equation}

\subsection{Case $q\ge 2$}
\label{sect-20-4-1-2}

Consider $q\ge 2$ (case $q=1$ we analyze later). Integrating the first term in the right-hand expression of (\ref{20-4-20}) by  $d\rho /\rho$ from $\bar{\rho}^*_1=C(\mu h)^{\frac{1}{2}}$ to $\bar{\rho}_1$ we get its value as $\rho=\bar{\rho}^*_1$; one can see easily that the result will be then less than $Ch^{1-d}$. The same integration  by $d\rho/\rho$ applied to the second term in the right-hand expression of (\ref{20-4-20}) results in $O(h^{1-d})$ for sure.

Thus for $q\ge 2$ in the whole intermediate zone 
$\{\bar{\rho}^*_1\le |\xi'''|\le \bar{\rho}_1\}$ we get a proper remainder estimate
\begin{multline}
\R^\T_{1\,Q}\Def 
|\Gamma (Q\psi \tilde{e}) - h^{-1} \int_{-\infty}^0 \sum_{0\le n\le \bar{n}}
\Bigl(F_{t\to h^{-1}\tau}\bar{\chi}_{T_n}(t)\Gamma Q_n\psi U\Bigr)\,d\tau |\le\\
Ch^{1-d}
\label{20-4-22}
\end{multline}
as $Q=\sum_{0\le n\le \bar{n}} Q_n $, where  $Q_n$ are operators with the symbols supported in  $\{\rho_n \le |\xi '''|\le \rho_{n+1}\}$, 
$\rho_n= 2^n\bar{\rho}^*_1$,  $\rho_{\bar{n}}\asymp \bar{\rho}_1$, and satisfying (\ref{20-4-3}), $Ch|\log h|/\rho_n^2 \le T_n\le \epsilon \rho_n$. Therefore we arrive to

\begin{proposition}\label{prop-20-4-2} 
As $q\ge2$  in the framework of proposition \ref{prop-20-4-1} estimate \textup{(\ref{20-4-22})} holds for  operator $Q=Q^\w$ with the symbol supported in  $\Omega_\interm= \{\bar{\rho}^*_1\le |\xi '''|\le \bar{\rho}_1\}$ and satisfying \textup{(\ref{20-4-3})}.
\end{proposition}

Now let us consider an expression 
\begin{equation}
h^{-1} \int_{-\infty}^0 \Bigl(F_{t\to h^{-1}\tau}\chi_T(t)\Gamma Q\psi U\Bigr)\,d\tau
\label{20-4-23}
\end{equation}
where again $Q$ is an operator with the symbol supported in $\cZ_\rho$, $T\in[T_*,T^*]$, $T^*=\epsilon\rho$, $\rho\in [\bar{\rho}^*_1,\bar{\rho}_1]$. Recall that $\chi \in \sC_0^\infty ([-1,-\frac{1}{2}]\cup [\frac{1}{2}, 1])$.
Also recall that one can rewrite (\ref{20-4-23}) as
\begin{equation}
T^{-1} \Bigl(F_{t\to h^{-1}\tau}\check{\chi} _T(t)\Gamma Q\psi U \Bigr)\Bigr|_{\tau =0}
\label{20-4-24}
\end{equation}
with $\check{\chi}(t)=i t^{-1}\chi(t)$ and due to proposition~\ref{prop-20-4-1} expression (\ref{20-4-24}) does not exceed 
\begin{equation}
C h^{-d} \rho^{q+2} \times \bigl( \frac{h} {T\rho^2} \bigr)^s. 
\label{20-4-25}
\end{equation}
Making a summation with respect to $T\in[T_*,T^*]$
we get its value at $T=\bar{T}=\epsilon\mu^{-1}$ i.e. 
$Ch^{-d}\rho^{q+2}(\mu h/\rho^2)^s$. Then after summation with respect to 
$\rho\in [\bar{\rho}^*_1,\bar{\rho}_1]$ we get its value at $\rho=\bar{\rho}^*_1$ i.e. $Ch^{-d}(\mu h)^{\frac{q}{2}+1}$. Estimating roughly contribution of the inner zone $\Omega_\inn= \{|\xi'''|\le C\bar{\rho}^*_1\}$ by 
$C\mu h^{1-d}(\bar{\rho}^*_1)^q$ and estimating contribution of the outer zone by $Ch^{1-d}$ due to results of Section~\ref{sect-20-2} we arrive to the estimate
\begin{multline}
\R^\T= |\Gamma (\psi \tilde{e}) - h^{-1} \int_{-\infty}^0 
\Bigl(F_{t\to h^{-1}\tau}\bar{\chi}_T(t)\Gamma \psi U\Bigr)\,d\tau |\le\\
C h^{1-d} + C(\mu h)^{\frac{q}{2}+1} h^{-d}
\label{20-4-26}
\end{multline}
with $T=\bar{T}=\epsilon \mu^{-1}$. 

On the other hand, due to the standard results rescaled for $q\ge 2$ 
\begin{multline}
|h^{-1} \int_{-\infty}^0 \Bigl(F_{t\to h^{-1}\tau}\bar{\chi}_T(t)\Gamma \psi U \Bigr)\,d\tau - h^{-d}\int \cN^\MW(x,\tau) \psi (x)\,dx | \le \\
Ch^{1-d} + C (\mu h)^{\frac{q}{2} +1}h^{-d}
\label{20-4-27}
\end{multline}
as $T=\bar{T}$ because with our choice~(\ref{20-3-5}) for  $\varepsilon$, a mollification error does not exceed $Ch^{1-d}$ as well.

Thus we arrive to 

\begin{theorem}\label{thm-20-4-3} 
Let conditions \textup{(\ref{20-1-1})}, \textup{(\ref{20-1-5})}, \textup{(\ref{20-1-7})}, \textup{(\ref{20-1-8})}--\textup{(\ref{20-1-11})}, and  \textup{(\ref{20-2-31})} be fulfilled, $(\bar{l},\bar{\sigma})=(2,1)$, $(l,\sigma)=(1,1)$. 

Let $q\ge 2$ and condition \textup{(\ref{20-3-3})} be fulfilled. 

\begin{enumerate}[label=(\roman*), fullwidth]
\item Then there exist two framing approximations (see footnote \footref{book_new-foot-18-16} of Chapter~\ref{book_new-sect-18}) of \cite{futurebook} such that for both of them
\begin{equation}
\R^\MW  \le Ch^{1-d} + Ch^{-d}(\mu h)^{\frac{q}{2}+1};
\label{20-4-28}
\end{equation}
in particular $\R^\MW  \le Ch^{1-d}$ as
\begin{equation}
\mu \le \bar{\mu}^*_{(q)} \Def Ch^{-\frac{q}{q+2}}.
\label{20-4-29}
\end{equation}
\item Furthermore, the same estimates hold for $\R^\W_{(\infty)}$.
\end{enumerate}
\end{theorem}

Now we can assume  instead of (\ref{20-3-3}) that
\begin{equation}
\bar{\mu}^*_{(q)} \le \mu \le \epsilon (h|\log h|)^{-1}.
\label{20-4-30}
\end{equation}

In this case we will use the following corollary of (\ref{20-4-20}):

\begin{proposition}\label{prop-20-4-4} 
Let $q\ge2$ and \textup{(\ref{20-4-30})} hold and $Q=Q(hD''')$, $Q=Q(\xi ''')$ be operator with symbol supported in $\cZ_\rho$, 
$\bar{\rho}^*_1\le \rho \le \bar{\rho}_1$,
$\psi \in \sC_0^\infty (B(0,1))$. Then
\begin{equation} 
\R^\T_{Q}=|\Gamma (Q\tilde{e}) - h^{-1} \int_{-\infty}^0 
\Bigl(F_{t\to h^{-1}\tau}\bar{\chi}_T(t)\Gamma Q\psi U\Bigr)\,d\tau |\le
Ch^{1-d}\rho^{q-1} 
\label{20-4-31}
\end{equation}
as $Ch|\log h|/\rho^2\le T\le T^*=\epsilon\rho$ and
\begin{multline}
\R^\T_{\bar{T},Q}=|\Gamma (Q\tilde{e}) - h^{-1} \int_{-\infty}^0 
\Bigl(F_{t\to h^{-1}\tau}\bar{\chi}_{\bar{T}}(t)\Gamma Q\psi U\Bigr)\,d\tau | \le\\ Ch^{1-d}\rho^{q-1} + Ch^{-d}\bigl( \frac {\mu h} {\rho^2} \bigr)^s\rho^q.
\label{20-4-32}
\end{multline}
with $\bar{T}=\epsilon \mu^{-1}$. 
\end{proposition}

Note that the right-hand expression of (\ref{20-4-31}) integrated by ${d\rho}/\rho$ does not exceed $Ch^{1-d}$.

\subsection{Case $q\ge 2$ (some calculations)}
\label{sect-20-4-1-3}

Now we can apply the technique of the previous subsection to calculate the Tauberian
expression
\begin{equation}
h^{-1} \int_{-\infty}^0 \Bigl(F_{t\to h^{-1}\tau}\phi_T(t) \Gamma Q\psi U\Bigr)\,d\tau
\label{20-4-33}
\end{equation}
with $\phi =\bar{\chi}$. Let us consider the difference 
\begin{equation}
h^{-1}\int_{-\infty}^0 \Bigl(F_{t\to h^{-1}\tau}\phi_T(t) \Gamma Q\psi (U-U')\Bigr)\,d\tau
\label{20-4-34}
\end{equation}
between (\ref{20-4-33}) and the same expression for the perturbed operator $A'=A+O(\mu^{-1})$. 
Then expression (\ref{20-4-34}) is equal to $I_0+I_1+\ldots +I_{\bar n}$ where $I_0$ is defined by (\ref{20-4-34}) with $\phi=\bar{\chi}$ and $T$ replaced by 
$T_*=\epsilon \mu^{-1}$  and $I_n$ are defined by (\ref{20-4-34}) with 
$\phi =\chi$ and $T=2^n \bar{T}$, 
$n=1,\ldots, \bar{n}$, $\bar{n}= \lceil T/\bar{T}\rceil$.

Applying arguments of the previous subsection and noting that 
\begin{claim}\label{20-4-35}
An operator norm of $\Bigl(e^{ih^{-1}tA}-e^{ih^{-1}tA'}\Bigr)$ does not exceed $C\mu^{-1}h^{-1}T$ as $|t|\asymp T$, 
\end{claim}
we conclude that $I_n$ does not exceed expression (\ref{20-4-25}) multiplied by $C\mu^{-1}h^{-1}T$ as $n\ge 1$, i.e.
\begin{multline}
|h^{-1}\int_{-\infty}^0 \Bigl(F_{t\to h^{-1}\tau}\chi_T(t) \Gamma Q\psi (U-U')\Bigr)\,d\tau|\le \\
Ch^{-d}\rho^{q+2}\times  
\bigl(\frac{h} {\rho^2 T} \bigr)^s \times \mu^{-1}h^{-1}T.
\label{20-4-36}
\end{multline}
Recall that with cut-off $\chi_T$ we can replace $h^{-1}$ by $T^{-1}$ and $\chi$ by $\check{\chi}$ like in (\ref{20-4-24}), also $\rho^2 \ge \mu h \ge \mu^{-1}$ due to $\mu \ge h^{-\frac{1}{2}}$ as $q\ge 2$. 

After summation with respect to $n$ (i.e. $T$) we get the same expression $I_0$ with  $T=\bar{T}$:
\begin{multline}
|h^{-1}\int_{-\infty}^0 \Bigl(F_{t\to h^{-1}\tau}\bigl(\bar{\chi}_T(t)- \bar{\chi}_{\bar{T}}(t)\bigr) \Gamma Q\psi (U-U')\Bigr)\,d\tau|\le \\
C\mu^{-2} h^{-d-1}\rho^{q+2} 
\bigl(\frac {\mu h} {\rho^2} \bigr)^s.
\label{20-4-37}
\end{multline}
Note that the right-hand expression in (\ref{20-4-37}), integrated by 
$d\rho/ \rho$  from $\bar{\rho}^*_1$ to $1$, does not exceed its value as $\rho=\bar{\rho}^*_1$ i.e. $C\mu^{-2}h^{-d-1}(\mu h)^{\frac{1}{2}(q+2)} = Ch^{1-d}(\mu h)^{\frac{1}{2}(q-2)}$.

Therefore for $q\ge 2$ and $Q$ supported in the intermediate zone expression (\ref{20-4-34}) is equal to $I_0$ modulo $O(h^{1-d})$. However, for $I_0$ we can apply the standard theory rescaled and replace it by the difference of the Weyl expressions for $A$ and $A'$ (we will do it in Section~\ref{sect-20-6}). So far we proved 

\begin{proposition}\label{prop-20-4-5} 
In the framework of proposition \ref{prop-20-4-4}
\begin{multline}
|\Gamma (Q\psi \tilde{e}) - h^{-1} \int_{-\infty}^0 
\Bigl(F_{t\to h^{-1}\tau}\bar{\chi}_T(t)\Gamma Q\psi U^0 - \\
F_{t\to h^{-1}\tau}\bar{\chi}_{\bar{T}}(t)\Gamma Q\psi (U-U^0)
\Bigr)\,d\tau |\le Ch^{1-d}
\label{20-4-38}
\end{multline}
where $U^0$ is a Schwartz kernel of $e^{ih^{-1}tA^0}$, $A^0=\cT^*\cA^0\cT$, \begin{equation}
\cA^0\Def \mu^2\sum _{q+1\le j,k\le q+r}\zeta_j^{\dag\,\w}a^{jk\,\w}\zeta_k^\w + 
h^2\sum _{j,k\le q} D_j g^{jk\,\w}D_k + a_0^\w
\label{20-4-39}
\end{equation}
is the main part of $\cA$.
\end{proposition}

\subsection{Case $q=1$}
\label{sect-20-4-1-4}

Let us repeat arguments of Subsubsection~\ref{sect-20-4-1-1}.1 first as $q=1$, 
\begin{gather}
(h|\log h|)^{-\frac{1}{3}}\le \mu \le \epsilon (h|\log h|)^{-1},
\label{20-4-40}\\
|\xi_1| \asymp \rho \ge \bar{\rho}^*_1 =C(\mu h)^{\frac{1}{2}}+C\mu^{-1}.
\label{20-4-41}
\end{gather}

The first problem is that the factor $(\rho^2 +\mu^{-1})$ in the right-hand expression of (\ref{20-4-17}) is too large as $\mu \le Ch^{-\frac{1}{2}}$; further, even if there are no unremovable cubic terms, one should take factor $(\rho^2 + \rho^{-1}h|\log h|)$ instead because we need to take $C\rho^{-1}h|\log h|$ scale with respect to $x_1$. We handled this problem in Chapter~\ref{book_new-sect-18} as $d=3$ and the idea from there (may be combined with the idea used in Subsection~\ref{sect-20-2-3}) is applicable now.

\medskip\noindent
(a) The factor problem is rather simple.  Let us assume first that 
$g^{jk}$, $F_{jk}$ are constant (and in particular there no unremovable cubic terms).

Then this problem arrises only as $\rho^3\le h|\log h|$ which implies 
$\mu \le h^{-\frac{1}{3}}|\log h|^{\frac{2} {3}}$. 

Let us consider a $\gamma$-covering with $\gamma=\rho^\delta$ and consider separately two types of elements:

\medskip\noindent
(i) Elements with $|\nabla V|\ge \zeta = |\log h|^{-2}$. If 
$|\partial_{x_1} V|\ge \zeta $ we have that  $\xi_1$-shift for time 
$T\in [T_*, \rho^{\delta_1}]$, $T_*=\epsilon \mu^{-1}$  is of magnitude no less than $\zeta T$ and then since logarithmic uncertainty principle is fulfilled
\begin{equation*}
\mu^{-1} \zeta \times \varepsilon \asymp \mu^{-2}\zeta \ge Ch|\log h|
\end{equation*}
we conclude that $\chi_T\Gamma Q\psi u$ is negligible;  otherwise 
$|\nabla_\perp V|\ge \zeta$ and similar arguments work 
for $(x'',\xi'')$-shifts. 

\medskip\noindent
(ii) Elements with $|\nabla V|\le \zeta$. We can take on them
$\gamma = C(\mu ^{-1}h)^{\frac{1}{2}}|\log h|$ which is larger than 
$C\rho^{-1}h|\log h|$ and replace the term $\gamma$ by 
$C(\zeta \gamma + \gamma ^l|\log h|^{-\sigma})$ which does not exceed $C\mu h$.
In this case one can take $T_1= \epsilon \rho \zeta^{-1}$. 

\medskip\noindent
(b) Similar analysis combined with ideas of subsection~\ref{sect-20-2-3} works in the general case as long as there are no unremovable cubic terms. We leave all the details to the reader.

\medskip\noindent
(c) If there are unremovable terms of this type one can replace decomposition (\ref{20-4-8}) by a similar decomposition with $\Upsilon_\alpha$ orthonormal eigenfunctions of an auxiliary operator 
\begin{gather}
\mathbf{a}= \cA_0 + \mu^{-1}\sum_{j,k,m\le 2r} \beta_{jkm}L_jL_kL_m+ C_0\mu^{-2}\cA_0^2
\label{20-4-42}\\
\shortintertext{with} 
\cA_0=\sum _{1\le j\le r} \bar{f}_j \bigl(\mu^2 h^2D_{q+j}^2 + x_{q+j}^2\bigr),
\label{20-4-43}
\end{gather}
where $L_{2j-1}=\mu hD_{q+j} $, $L_{2j}=x_{q+j}$ and $\Upsilon_\beta$ are replaced by their complex conjugates). We added the last term in (\ref{20-4-42}) to make operator non-negative and self-adjoint; this term is $O(\mu^{-2})$ on $U$ and thus it is included in the approximation error estimate $O(\mu^{-l}|\log \mu|^{-\sigma})$ anyway. We got  operator (\ref{20-4-43}) by rescaling $x'\mapsto \mu^{-1}x'$, $D_j\mapsto \mu D_j$.

For an eigenvalue counting function $\mathbf{n}(\tau,\hbar)$ of this operator $\mathbf{a}$ in $\sL^2(\bR^r)$ the standard semiclassical eigenvalue asymptotics holds
\begin{gather}
\mathbf{n}(\tau, \hbar)=\mathsf{n}^\W (\tau, \hbar) + O(\hbar^{1-r}), 
\qquad \mathsf{n}^\W \asymp \hbar ^{-r}\qquad {\textrm  as\ } \hbar\to +0
\label{20-4-44}\\
\intertext{and in particular}
\mathbf{n}(\tau +\hbar , \hbar)- \mathbf{n}(\tau , \hbar) =O(\hbar^{1-r})
\label{20-4-45}
\end{gather}
where $\mathsf{n}^\W$ means Weyl expression for $\mathbf{a}$ and we plug 
$\hbar= \mu h$.

Then the same modifications hold for estimates (\ref{20-4-17}) and (\ref{20-4-19}) as well and then the standard Tauberian arguments improve estimate (\ref{20-4-20}) to 
\begin{equation}
\R^\T_{Q}\le 
C \Bigl(h^{-d} \mu^{-1} \rho ^2 \times 
\bigl( \frac {\mu h} {\rho^2 } \bigr)^s+Ch^{1-d}\Bigr) |\log \rho|^{-2}
\label{20-4-46}
\end{equation}
as $\mu\le h^{\delta-1}$ because $T^*\ge \epsilon_1\rho |\log \rho|^2$ now.

Integrating the right-hand expression over ${d\rho}/\rho $ from $\bar{\rho}^*_1$ to $\bar{\rho}_1$ we get $O(h^{1-d})$ then as long as $\mu \le h^{\delta-1}$.  However  for  $ h^{\delta -1}\le \mu \le \epsilon (h \log h|)^{-1}$ we can take only  $T^*=\epsilon_1 \rho$ and recover estimate $O(h^{1-d}|\log h|)$.

On the other hand, under microhyperbolicity assumption for 
$\mu \le \epsilon (h|\log h|)^{-1}$ we can take  $T^*=\epsilon_1$
and recover estimate $O(h^{1-d})$ again.

Thus we have proven

\begin{proposition}\label{prop-20-4-6}
Let $q=1$, $(l,\sigma)=(1,2)$. Then
\begin{enumerate}[label=(\roman*),fullwidth]
\item  Estimate \textup{(\ref{20-4-28})} holds. In particular for 
$\mu \le Ch^{-\frac{1}{3}}$
sharp remainder estimate $\R^\MW\le Ch^{1-d}$ holds.

\item  In the framework of proposition \ref{prop-20-4-2} $\R^\T_{1\,Q}$ defined by \textup{(\ref{20-4-22})} does not exceed $Ch^{1-d}|\log h|$ where again 
$T_n=\epsilon \rho_n$.

\item  If \underline{either} $\mu \le h^{\delta-1}$ \underline{or} $f_j$ have constant multiplicities \underline{or} microhyperbolicity condition is fulfilled 
at level $0$, then $\R^\T_{1\,Q}$ defined by \textup{(\ref{20-4-22})} does not exceed $Ch^{1-d}$ where $T_n= \rho_n |\log \rho_n|^2$.

\item Furthermore, if $Q$ is supported in $\{|\xi'''|\ge \rho\}$ and equal $1$
in $\{2|\xi'''|\ge \rho\}$ then $\R^\T_{\bar{T}\,Q}$ defined by \textup{(\ref{20-4-32})} does not exceed 
$Ch^{1-d} + Ch^{-d}\rho  \bigl( {\mu h}/{\rho^2}\bigr)^s$.
\end{enumerate}
\end{proposition}

\subsection{Case $q= 1$ (calculations)}
\label{sect-20-4-1-5}

Let us apply as $q=1$ the same arguments as in Subsubsection~\ref{sect-20-4-1-3}. However as perturbation is $O(\mu^{-1})$ it makes sense only for $\mu \ge h^{-\frac{1}{2}}$ resulting in 
$Ch^{1-d}(\mu h)^{-\frac{1}{2}}=C\mu^{-\frac{1}{2}}h^{\frac{1}{2}-d}$ which is less than the weak magnetic field estimate only in this case.

On the other hand, as perturbation is $O(\mu^{-1})$ these arguments make sense as $\mu \ge h^{-\frac{1}{3}}$ resulting in 
$C\mu^{-\frac{3}{2}}h^{\frac{1}{2}-d}=O(h^{1-d}$.

Leaving easy details to the reader we arrive to the following 

\begin{proposition}\footnote{\label{foot-20-10} cf. proposition~\ref{prop-20-4-5}} \label{prop-20-4-7}
In the framework of proposition~\ref{prop-20-4-6}

\begin{enumerate}[label=(\roman*), fullwidth]
\item Estimate \textup{(\ref{20-4-38})} holds as \underline{either} there are no $3$-rd order resonances \underline{or} we include in unperturbed operator $\cA^0$ non-removable $O(\mu^{-1})$ terms like in \textup{(\ref{20-4-42})};

\item If there are non-removable $O(\mu^{-1})$ terms and we do not include them into $\cA^0$ then as $\mu \ge h^{-\frac{1}{2}}$ the left-hand expression of \textup{(\ref{20-4-38})} does not exceed $C\mu^{-\frac{1}{2}}h^{\frac{1}{2}-d}$.
\end{enumerate}
\end{proposition}

We will return to the intermediate zone (for $q=1$) later under microhyperbolicity or non-degeneracy assumptions.

\section{Inner zone}
\label{sect-20-4-2}

Now let us consider  an \emph{inner zone\/}\index{zone!inner}
\begin{equation}
\Omega_\inn=\{|\xi'''|\le \bar{\rho}^*_1= (\mu h)^{\frac{1}{2}}\}.
\label{20-4-47}
\end{equation}
\subsection{Case $q\ge 2$}
\label{sect-20-4-2-1}

Assume that $q\ge 2$ and (\ref{20-4-30}) holds. Then one can see easily that 

\begin{claim}\label{20-4-48}
The perturbation of the magnitude $O(\mu^{-1})$ in 
$\Omega_\inn \cup \Omega_\interm$ leads to an error 
$O(\mu^{-1}\bar{\rho}^{*\, q}_1 h^{-d})=O(h^{1-d})$ in $h^{-d}\cN^\MW_Q$.
\end{claim}

Therefore as $q\ge 2$ instead of $\cA$ one can consider in 
$\Omega_\inn \cup \Omega_\interm$ a reduced operator $\cA^0$ with the coefficients
\begin{align*}
&a^{jk\,\w}=a^{jk}(x'',x''',\mu^{-1}hD''),
	&&g^{jk\,\w}=g^{jk}(x'',x''',\mu^{-1}hD''),\\
&a_0^\w =a_0(x'',x''',\mu^{-1}hD''), 
	&&\zeta_j^\w = \mu^{-1}hD_j + i x_j
	\end{align*}
where $a^{jk}=0$ unless $\exists \fm : j,k \in \fm $; this substitution leads to an error $O(h^{1-d})$ in $\R^\T _{Q}$. However we provided different arguments for inner zone $\Omega_\inn$ and intermediate zone   $\Omega_\interm$ to justify this replacement.

Moreover, if $f_j$ had constant multiplicities we would have even 
\begin{equation}
\cA^0=\mu^2 \sum _{j\le r} f_j \bigl(h^2 D_j^2 +\mu^2 x_j^2 - \mu h\bigr) +
h^2\sum _{j,k\ge 2r+1} D_j^\w g^{jk\,\w}D_k^\w + a_0^\w
\label{20-4-49}
\end{equation}
and decomposing $\cT U$ into (\ref{20-4-8}) series  we would arrive to the family of operators
\begin{equation}
\cA_\alpha\Def h^2\sum _{j,k\ge 2r+1} D_j g^{jk\,\w}D_k  + W_\alpha ^\w, \qquad 
W_\alpha = a_0+ \sum _{j\le r} (2\alpha_j+1) \mu h f_j
\label{20-4-50}
\end{equation}
with $\alpha\in \bZ^{+\,r}$ and we would be able to apply without any significant modifications analysis of Sections~\ref{book_new-sect-4-5} and~\ref{book_new-sect-5-3} of \cite{futurebook}.

However, we need to consider more general operator (\ref{20-4-39}). Consider operator $Q=Q^\w$ with symbol supported in $\cZ_\rho$ with $\rho \le \bar{\rho}^*_1$. 

Repeating arguments of the proof of proposition \ref{prop-20-4-1} we arrive to estimate (\ref{20-4-16}) where $Q'=\cT Q\psi '\cT^*$, $\psi'$ is an element of $\gamma$-admissible partition in $(x'',x''',\xi'')$, $\gamma \ge \mu h$. 

Restoring to $Q$ we arrive to the estimate
\begin{equation}
|F_{t\to h^{-1}\tau} \chi_T(t)\Gamma Q\psi U |\le 
C h^{-d} T \rho^q \gamma 
\Bigl(\bigl( \frac {h} {\rho^2 T} \bigr)^s + 
\bigl(\frac {h} {\rho \gamma} \bigr)^s + 
\bigl( \frac {\mu^{-1}h} { \gamma^2} \bigr)^s
\Bigr).
\label{20-4-51}
\end{equation}
As $q\ge 2$ and $\mu\ge \bar{\mu}^*_{(q)}$ we can take $\gamma \ge  \mu h$
and $\rho \ge \bar{\rho}_0 \ge \mu^{-\frac{1}{2}}$; remaining \emph{inner core\/} 
\begin{equation}
\Omega^0_\inn= \{|\xi'''|\le \bar{\rho}_0\}
\label{20-4-52}
\end{equation}
will be considered separately. Then both $h/(\rho \gamma) = 1/(\mu \rho)$ and 
$\mu^{-1}h/\gamma^2 = 1/(\mu^3 h)$ do not exceed $h^\delta$ and we arrive to the estimate
\begin{equation}
|F_{t\to h^{-1}\tau} \chi_T(t)\Gamma Q\psi U |\le 
C \mu h^{1-d} T \rho^q 
\bigl(  \frac {h} {\rho^2 T} \bigr)^s 
\label{20-4-53}
\end{equation}
which in turn implies
\begin{multline}
|F_{t\to h^{-1}\tau} \bar{\chi}_T(t)\Gamma Q\psi U |\le 
C \mu h^{1-d} \rho^q \times \frac{h}{\rho^2} = C\mu h^{2-d} \rho^{q-2}, 
\\ 
\text{as\ \ }  T_* \Def \epsilon \frac {h}{\rho^2} \le T \le T^*=\epsilon \rho.
\label{20-4-54}
\end{multline}
Therefore, the contribution of $\cZ_\rho$ to the remainder  is 
$O(\mu h^{2-d}\rho^{q-3})$; integrated over ${d\rho}/\rho $ it results in
$O(\mu h^{2-d})=O(h^{1-d})$ as $q\ge 4$,
$O(\mu h^{2-d}|\log h|)=O(h^{1-d})$ as $q= 3$, 
$O(\mu h^{2-d}\bar{\rho}_0 ^{-1})$ as $q=2$.

On the other hand, contribution of the inner core $\Omega^0_\inn$ 
to the remainder  is $O(\mu h^{1-d} \bar{\rho}_0^q)$.

Therefore we arrive to the remainder estimate 
$O\bigl(h^{1-d}+ \mu h^{1-d} \bar{\rho}_0^q\bigr)$ as $q\ge 3$ and
$O\bigl(\mu h^{1-d} (\bar{\rho}_0^2 + h\bar{\rho}_0^{-1})\bigr)$ as $q=2$
and we need to pick-up 
\begin{equation}
\bar{\rho}_0=\max\bigl(h^{\frac{1}{3}}, \mu ^{- \frac {1}{2}}\bigr),
\label{20-4-55}
\end{equation}
then finally arriving to the remainder estimate $O(h^{1-d})$ as $q\ge 3$ and $O\bigl( h^{1-d}+\mu h^{ \frac {5}{3} -d}\bigr)$ as $q=2$\,\footnote{\label{foot-20-11} Similarly we would get 
$O\bigl( h^{1-d}+\mu h^{ \frac {4}{3} -d}\bigr)$ as $q=1$ but there is a lot of other things to consider.}.

Thus we have proven

\begin{proposition}\label{prop-20-4-8}
Let $q\ge 2$, condition \textup{(\ref{20-4-30})} be fulfilled and $Q=Q(hD''')$ 
be operator with the symbol $Q(\xi ''')$ supported in 
$\{|\xi'''| \le C\bar{\rho}_1\}$, $\psi \in \sC_0^\infty (B(0,1))$.  Then
\begin{enumerate}[label=(\roman*), fullwidth]
\item Estimate 
\begin{equation}
\R^\T_{Q}\le  C h^{1-d}
\label{20-4-56}
\end{equation}
holds provided \underline{either} $q\ge 3$ \underline{or} $q=2$ and $\mu \le h^{-\frac{2}{3}}$;

\item Estimate 
\begin{equation}
\R^\T_{Q}\le  C \mu h^{\frac{2}{3}-d}
\label{20-4-57}
\end{equation}
holds as $q=2$ and $\mu \ge h^{-\frac{2}{3}}$;
\end{enumerate}
where in both cases $\bar{\chi} \in \sC_0^\infty ([-1, 1])$, $\bar{\chi}=1$ on 
$[-\frac{1}{2},\frac{1}{2}]$, $Q =\sum_{0\le n\le \bar {n}} Q_n$,
$Q_n$ are operators with symbols supported in 
$\{\frac{1}{2}\rho_n\le |\xi'''|\le 2\rho_n \}$ as $n=1,\ldots,\bar{n}$, 
and in $\{|\xi'''|\le 2\rho_0\}$ as $j=0$, $\rho_n= 2^n\rho^*_0$, 
${\bar n}=\lfloor \log \bar{\rho}_1/\rho^*_0\rfloor +1$,
$T_n = \epsilon \rho_n$ ($n\ge 1$), $T_0=\epsilon \mu^{-1}$.

Further, both statements remain true with $U$ replaced by $U^0$ which is the Schwartz kernel of $e^{ih^{-1}tA_0}$.
\end{proposition}

\subsection{Case $q=1$}
\label{sect-20-4-2-2}

In contrast to the previous case $q\ge 2$ skipping $O(\mu^{-1})$ terms is not now generally justified. Because of this we need to modify now arguments of Subsubsection~\ref{sect-20-4-2-1}.1, in the same manner as we modified arguments of  Subsubsections~\ref{sect-20-4-1-2}.2--\ref{sect-20-4-1-3}.3 in Subsubsections~\ref{sect-20-4-1-4}.4--\ref{sect-20-4-1-4}.4. 

So let us proceed without skipping such terms first, assuming that
\begin{gather}
Ch^{-\frac{1}{3}}\le \mu \le \epsilon (h|\log h|)^{-1},\label{20-4-58}\\
C\mu^{-1}\le |\xi_1| \asymp \rho \le \bar{\rho}^*_1 =(\mu h)^{\frac{1}{2}}. 
\label{20-4-59}
\end{gather}

Starting from (\ref{20-4-51}) with $\gamma = \mu h$ due to the above arguments we arrive to
\begin{multline}
|F_{t\to h^{-1}\tau} \bigl(\bar{\chi}_T(t)-\bar{\chi}_{T_*}(t)\bigr)
\Gamma Q\psi U |\le \\
C \mu h^{1-d} \rho \Bigl(T_* \bigl( \frac{h} {\rho^2 T_*} \bigr)^s + 
T\bigl(\frac {h} {\rho \gamma}\bigr)^s + 
T\bigl( \frac {\mu^{-1}h} { \gamma^2} \bigr)^s \Bigr)
\label{20-4-60}
\end{multline}
as $T_*= {h}/ {\rho^2 }  \le T\le T^*=\epsilon \rho $ and thus 
\begin{equation}
|F_{t\to h^{-1}\tau} \bar{\chi}_T(t) \Gamma Q\psi U |\le \\
C \mu h^{1-d} \rho 
\Bigl( \frac {h} {\rho^2}  + T\bigl( \frac{h} {\rho \gamma} \bigr)^s + 
T\bigl( \frac {\mu^{-1}h} { \gamma^2} \bigr)^s \Bigr).
\label{20-4-61}
\end{equation}
Therefore due to the standard Tauberian arguments 
\begin{multline}
|\Gamma (Q\psi \tilde{e}) - h^{-1} \int_{-\infty}^0 
\Bigl(F_{t\to h^{-1}\tau} \bar{\chi}_T(t)\Gamma Q\psi U\Bigr)\,d\tau |\le\\
\shoveright{C \mu T^{-1} h^{1-d} \rho 
\Bigl( \frac {h} {\rho^2}  + T\bigl( \frac {h} {\rho \gamma} \bigr)^s + 
T\bigl( \frac {\mu^{-1}h} {\gamma^2} \bigr)^s \Bigr)\asymp }\\
C \mu h^{1-d} \rho 
\Bigl( h\rho^{-3} + 
\bigl( \frac {h} {\rho \gamma} \bigr)^s + 
\bigl( \frac {\mu^{-1}h} {\gamma^2} \bigr)^s \Bigr).
\label{20-4-62}
\end{multline}
After integration the right-hand expression over ${d\rho}/{\rho}$ with $\rho$ ranging from $\bar{\rho}_0$ to $\bar{\rho}^*_1$,
\begin{gather}
\bar{\rho}_0= Ch^{\frac{1}{3}}\ge C\mu^{-1}\qquad \text{as\ \ } 
\mu\ge h^{-\frac{1}{3}},
\label{20-4-63}\\
\shortintertext{we get}
C \mu h^{1-d} \rho 
\Bigl( h\rho^{-3} + 
\bigl( \frac {h} {\rho \gamma} \bigr)^s \Bigr)\Bigr|_{\rho=\bar{\rho}_0}+
C \mu h^{1-d} \rho 
\bigl( \frac {\mu^{-1}h} {\gamma^2} \bigr)^s\Bigr|_{\rho=\bar{\rho}^*_1}
\label{20-4-64}
\end{gather}
where all terms except 
$C\mu h^{2-d} \bar{\rho}_0^{-2}= C\mu h^{\frac{4}{3}-d}$ do not exceed $Ch^{1-d}$.

Also, contribution of the inner core $\Omega^0_\inn=\{|\xi_1|\le \bar{\rho}_0\}$ does not exceed $C\mu h^{1-d} \bar{\rho}_0= C\mu h^{ \frac{4}{3}-d}$ as well.

Therefore, the final estimate is $Ch^{1-d}+C\mu h^{ \frac {4}{3}-d}$.
Thus we arrive to

\begin{proposition}\label{prop-20-4-9} 
Let $q=1$, condition \textup{(\ref{20-4-30})} be fulfilled and $Q=Q(hD''')$ 
be operator with the symbol $Q(\xi ''')$ supported in 
$\{|\xi'''| \le C\bar{\rho}_1\}$, $\psi \in \sC_0^\infty (B(0,1))$.  Then
\begin{equation}
\R^\T_{1\,Q} \le Ch^{1-d} + C\mu h^{ \frac{4}{3}-d }
\label{20-4-65}
\end{equation}
where $\R^\T_{1\,Q}$ is defined by \textup{(\ref{20-4-22})} with 
$Q =\sum_{0\le n\le \bar{n}} Q_n(hD''')$, $Q_n$ are operators with symbols supported in $\{ \frac{1}{2}\rho_n\le |\xi'''|\le 2\rho_n \}$ as 
$n=1,\ldots, \bar{n}$, 
and in $\{|\xi'''|\le 2\rho_0\}$ as $n=0$, $\rho_n= 2^n\bar{\rho}_0$, 
$\bar{n}=\lfloor \log \bar{\rho}_1/\bar{\rho}_0\rfloor +1$,
$T_n = \epsilon \rho_n$ ($n\ge 1$), $T_0=\epsilon \mu^{-1}$.
\end{proposition}

\section{Cases $q=1,2$ revised}
\label{sect-20-4-3}

We repeat arguments of Subsection~\ref{prop-20-4-2} as $q=1,2$, but we will use certain additional assumptions to improve remainder estimates. 

\subsection{Microhyperbolicity assumption}
\label{sect-20-4-3-1}
Let us assume first that the microhyperbolicity assumption (see definition~\ref{def-20-1-2}) is fulfilled. Note that (\ref{20-4-11}) implies that
\begin{multline}
|F_{t\to h^{-1}\tau} \chi_T(t)\Gamma Q'\psi U|\le 
C \mu ^r h^{r-d} T \rho^q \times \\
\Bigl(\bigl( \frac {h} {\rho^2 T} \bigr)^s + 
\bigl( \frac {h} {\rho \gamma} \bigr)^s + 
\bigl( \frac {\mu^{-1}h} { \gamma^2} \bigr)^s
\Bigr) \times \\
\gamma^s \int_{B(z,2\gamma)} 
\1\tr \bigl(\cA(x',x''',\xi'')+i\gamma \bigr)^{-s}\1 \,dx''\,dx'''\,d\xi''
\label{20-4-66}
\end{multline}
where $B(z,\gamma)\supset \supp Q'$, $\cA $ is an operator reduced and considered as an operator in $\sL^2(\bR^r, \bC)$ and $\1.\1$ is an operator norm there. Then the same estimate (with an integral over $B(0,1)$) holds as $Q'$ is replaced by $Q$.

Further, one can prove easily that
\begin{claim}\label{20-4-67}
Due to the microhyperbolicity assumption the last integral in (\ref{20-4-66}) is of magnitude  $\gamma (\mu h)^{-r}$ as $\gamma \ge \rho^2$; even condition 
$\gamma \ge \mu^{-1}$ is not needed anymore if we consider operator (\ref{20-4-41}) from the very beginning\footnote{\label{foot-20-12} As $q\ge 2$ we can skip all the $O(\mu^{-1})$ terms. As $q=1$ we preserve all the terms which do not contain factor $hD_1$, we skip all such terms containing $(hD_1)^2$ because this leads to the approximation error $\rho ^2\mu^{-1}$ and then we remove as in Chapter~\ref{book_new-sect-18} all such terms containing $hD_1$ exactly in power $1$.}.
\end{claim}

Then (\ref{20-4-66}) becomes
\begin{multline}
|F_{t\to h^{-1}\tau} \chi_T(t)\Gamma Q\psi U |\le \\
C h^{-d} \gamma T \rho^q \times 
\Bigl(\bigl( \frac {h} {\rho^2 T} \bigr)^s + 
\bigl( \frac {h} {\rho \gamma} \bigr)^s + 
\bigl( \frac {\mu^{-1}h} {\gamma^2} \bigr)^s\Bigr) 
\label{20-4-68}
\end{multline}
as $T_*=h/\rho^2 \le T\le T^*=\epsilon \rho$ and picking up $\gamma =\rho^2$, 
$\rho \ge \bar{\rho}_0= h^{\frac{1}{3}-\delta}\ge \mu^{-1}$  we rewrite the right-hand expression here as 
$\displaystyle{C h^{-d} T \rho^{q+2} \bigl( \frac {h} {\rho^2 T} \bigr)^s }$.

\begin{remark}\label{rem-20-4-10}
Note that we changed here (under microhyperbolicity assumption) definition of $\bar{\rho}_0$ and the corresponding notion of the inner core $\Omega^0_\inn$.
\end{remark}

Then (\ref{20-4-68}) implies that
\begin{equation}
|F_{t\to h^{-1}\tau} \bar{\chi}_T(t)\Gamma Q\psi U |\le 
C h^{1-d} \rho^q 
\label{20-4-69}
\end{equation}
as long as $T\le \epsilon \rho$, $\rho\ge \bar{\rho}_0$. 

However, due to the microhyperbolicity assumption 
$F_{t\to h^{-1}\tau} \chi_T(t)\Gamma Q\psi U$ is negligible as long 
$Ch|\log h|\varepsilon ^{-1}\le T\le \epsilon$ and $T=\epsilon \rho $ satisfies this inequality as long as 
\begin{equation}
\rho \ge C\bar{\rho}_0,\quad \varepsilon \ge Ch\rho^{-1}|\log h|.
\label{20-4-70}
\end{equation}
Therefore, (\ref{20-4-70}) yields that (\ref{20-4-69}) holds for 
$T\le \epsilon$.

Then we can apply the standard Tauberian arguments resulting in the estimate with the right-hand expression $Ch^{1-d}\rho^q$. This last estimate holds for 
$\rho \le \bar{\rho}_1$ as well. Integration of this expression over ${d\rho}/\rho$ implies the remainder estimate $O(h^{1-d})$.

Consider now the contribution of the inner core 
$\Omega^0_\inn=\{|\xi'''|\le \bar{\rho}_0\}$. In this zone we get estimate
\begin{equation}
|F_{t\to h^{-d}\tau} \bar{\chi}_T(t)\Gamma Q\psi U|\le 
C h^{-d} \bar{\rho}_0^{q +3},
\label{20-4-71}
\end{equation}
as $T=\rho$ and $\rho\Def \bar{\rho}_0$ there; then the standard microhyperbolicity arguments push it to $T=\epsilon$ under condition (\ref{20-4-70}) and then the standard Tauberian arguments imply that the contribution of this zone to the remainder does not exceed
$C h^{-d} \bar{\rho}_0^{q +3}$ which is $O(h^{1-d})$ as $q\ge 1$ and  $\delta$ small enough.

Thus under microhyperbolicity assumption the Tauberian remainder estimate is $O(h^{1-d})$ even as $q=1,2$ but we need to look at an approximation error; one can prove easily that for smallest $\varepsilon \Def C\rho^{-1}h|\log h|$ satisfying (\ref{20-4-70}) it is also $O(h^{1-d})$. 

Thus we arrive to

\begin{proposition}\label{prop-20-4-11}
Under the microhyperbolicity assumption $\R^\T_{Q}$ does not exceed $Ch^{1-d}$ where $\R^\T_{Q}$ is given by \textup{(\ref{20-4-20})} with  $T=\epsilon $.
\end{proposition}

\subsection{Case of constant $g^{jk}$, $F_{jk}$}
\label{sect-20-4-3-2}

Recall that if $F_{jk}$, $g^{jk}$ are constant then there are no non-removable $O(\mu^{-1})$ terms and the microhyperbolicity assumption is (\ref{20-2-41}).

Using arguments similar to those used in Subsection~\ref{book_new-19-4-7} of \cite{futurebook} one can prove easily that one can replace it by non-degeneracy condition (\ref{20-2-42}):

\begin{Problem}\label{Problem-20-4-12}
Let conditions \textup{(\ref{20-1-1})}, \textup{(\ref{20-1-5})}, \textup{(\ref{20-1-7})}, \textup{(\ref{20-1-8})}, $\textup{(\ref{20-1-9})}_3$, \textup{(\ref{20-1-11})}, and  \textup{(\ref{20-2-31})} be fulfilled, and $\mu \le \epsilon( h|\log h|)^{-1}$.

Further, let $g^{jk}$, $F_{jk}$ be constant, $(l,\sigma)=(2,0)$ and assumption (\ref{20-2-42}) be fulfilled.

Prove that  $\R^\T =O(h^{1-d})$.
\end{Problem}

\subsection{Case of $f_j$ having constant multiplicities}
\label{sect-20-4-3-3}

Assume that $f_j$ have constant multiplicities. Then in both  inner and intermediate zones after reduction and skipping $O(\mu^{-1})$ terms which are due to the $3$-rd order resonances and all smaller terms\footnote{\label{foot-20-13} The case of $q=1$ where there are such terms will be considered in details in the next Subsection~\ref{sect-20-4-3-4}.4 } instead of ``matrix'' operator we get a family of ``scalar'' operators 
\begin{multline}
A_\alpha (x'',x''', \mu^{-1}hD'',hD''')=\\
\sum_{q+1\le j,k\le q+r}(hD_j) g^{jk} (x'',x''',\mu^{-1}hD''')(hD_k) + V_\alpha (x'',x''',\mu^{-1}hD''')
\label{20-4-72}
\end{multline}
and we can study them separately. Let us apply the standard rescaling procedure: for each index $\alpha$ we introduce functions $\gamma=\gamma_\alpha$ and $\varrho=\varrho_\alpha$ by
\begin{multline}
\gamma =\\ 
\inf \bigr\{t: \vartheta (t)\ge \rho^2 +|V_\alpha (x'',x''',\xi'')|,
\vartheta (t) t^{-1}\ge |\nabla V_\alpha (x'',x''',\xi'')|\bigr\}+\bar{\gamma},
\label{20-4-73}
\end{multline}
\begin{equation}
\varrho = \vartheta (\gamma)^{\frac{1}{2}},\label{20-4-74}\\
\end{equation}
with
\begin{multline}
\bar{\varrho}=\vartheta(\bar{\gamma})^{\frac{1}{2}},\qquad
\bar{\varrho}\bar{\gamma}= Ch \implies \\
\bar{\varrho}=h^{\frac {l} {l+2}}|\log \mu|^{- \frac {l\sigma} {2(l+2)} },\quad
\bar{\gamma}=h^{\frac {2} {l+2}}|\log h|^{\frac {\sigma} {l+2} }
\label{20-4-75}
\end{multline}
where we assume that 
\begin{equation}
(1,1)\preceq (l,\sigma)\preceq (2,0).
\label{20-4-76}
\end{equation}
Note that  $\gamma^2 \ge C\mu^{-1}h|\log h|$.  

\begin{enumerate}[label=(\roman*),fullwidth]

\item Consider first elements with $\gamma \ge C{\bar\gamma}$. Calculating contribution of each  element of $(\gamma; \rho)$-partition with respect to $(x'',x''',\xi''; \xi''')$ we can apply elliptic arguments unless 
$\epsilon \varrho \le \rho \le C\varrho$; in the latter case we can apply microhyperbolic arguments; then the contribution of each element to the remainder estimate does not exceed 
$C\mu ^r h^{r-d+1}\varrho ^{q-1} \gamma^{d-1}$ in the latter case and
\begin{equation*}
C\mu ^r h^{r-d+1}\rho ^{q-1} \gamma^{d-1}\times 
\bigl(\frac {h} {(\rho+\varrho)\gamma} \bigr)^s
\end{equation*}
in the former one.

Then the contribution of each such element $\cW$ of $\gamma$-partition with respect to $(x'',x''',\xi')$\,\footnote{\label{foot-20-14} More precisely,  of 
$\cW\times \{|\xi'''|\le C\rho^*_1\}$.} 
to the remainder also does not exceed $C\mu ^r h^{r-d+1}\varrho ^{q-1} \gamma^{d-1}$. 

Then the total contribution of such elements to the Tauberian remainder  for given $\alpha $ does not exceed
\begin{equation*}
C\int \mu ^r h^{r-d+1}\varrho_\alpha  ^{q-1} \gamma_\alpha ^{-1} \,dz
\end{equation*}
with $z=(x'',x''',\xi'')$ and after summation over all indices we get
\begin{equation}
C\sum_{\alpha \in \bZ^{+\,r}}
\mu ^r h^{r-d+1}\int \varrho_\alpha  ^{q-1} \gamma_\alpha ^{-1} \,dz.
\label{20-4-77}
\end{equation}
\item  On the other hand, contribution of each element of $\gamma$-partition with $\gamma \le C{\bar\gamma}$ does not exceed 
$C\mu^r h^{r-d}\gamma^d \bar{\varrho}^q$ and their total contribution does not exceed (\ref{20-4-77}) as well.

Since for given $(z,\rho)$
\begin{equation}
\#\{\alpha \in\bZ^{+\,r}: \epsilon \rho^2\le \varrho_\alpha^2\le c\rho^2\}\le
C(\mu h)^{-r}(\rho^2+\mu h),
\label{20-4-78}
\end{equation}
expression (\ref{20-4-77}) does not exceed
\begin{equation}
C\mu  h^{2-d}\int \varrho_1 ^{q-1} \gamma_1  ^{-1} \,dz + 
C h^{1-d} \int \rho^q\gamma (\rho)^{-1}\,d\rho
\label{20-4-79}
\end{equation}
where $\gamma_1=\min_\alpha \gamma_\alpha$ and $\varrho_1$ corresponds to it,
$\gamma(\rho)$ is defined from equation $\vartheta(\gamma)=\rho$. 

One can see easily that the second term in (\ref{20-4-79}) does not exceed $Ch^{1-d}$, while the first term does not exceed 
\begin{equation*}
C\mu h^{2-d}\int \gamma^{ \frac{1}{2} (q-1)l -2}
|\log \gamma|^{- \frac{1}{2} (q-1)\sigma}\,d\gamma
\end{equation*}
which in turn does not exceed $C\mu h^{1-d}\bar{\varrho}^q$ unless $q=2$, $l=2$.
In the latter case $\sigma\le 0$ and the first term does not exceed 
$C\mu h^{2-d}|\log h|^{-\sigma}$. In both cases we arrive to the remainder estimate
\begin{equation}
\R^\T_{1\,Q}\le  Ch^{1-d} +
C \mu h^{1-d+ {\frac {ql} {l+2}}}|\log \mu|^{- \frac {q\sigma} {(l+2)} }
\label{20-4-80}
\end{equation}
which is exactly remainder estimate  (\ref{20-6-88}) below. Recall that it is 
$C\mu h^{2-d}$ as $(l,\sigma)=(2,0)$. Here $T=\gamma /\varrho$.
\end{enumerate}

Thus we have proven

\begin{proposition}\label{prop-20-4-13} 
Let $q= 1,2$, $\mu \ge \bar{\mu}^*_{(q)}$, $f_j$ have constant multiplicities and $\cA$ contain no unremovable cubic terms\footnote{\label{foot-20-15} May be we just dropped them in the interior zone causing an error $O(h^{1-d})$ if \underline{either} $q=2$ \underline{or} $q=1$ and each of these terms contained factor $\mu^{-1}hD_1$ and  an error  $O(\mu^{\frac{1}{2}} h^{1-d})$ in the general case as $q=1$.}. Let $Q=Q(hD''')$, $Q=Q(\xi ''')$ be operator with the symbol supported in 
$\{|\xi'''| \le C\bar{\rho}_1\}$, $\psi \in \sC_0^\infty (B(0,1))$. 

Then estimate \textup{(\ref{20-4-80})} holds with $\R^\T_{1\,Q}$ defined by \textup{(\ref{20-4-22})} with\linebreak 
$Q =\sum_{0\le n\le \bar {n}} Q_n(hD''')$, where
$Q_n$ are operators with symbols supported in 
$\{\frac{1}{2}\gamma_n\le \gamma\le 2\gamma_n \}$ as $n=1,\ldots,\bar{n}$, 
and in $\{\gamma \le 2\bar{\gamma}\}$ as $n=0$,
$\bar{n}=\lfloor \log \bar{\gamma}\rfloor +1$, $\gamma_n= 2^{-n}$,
$T_n = C \gamma_n ^{1-\frac{l}{2}} |\log \gamma_n|^{-\frac{\sigma}{2}}$
($n\ge 1$), $T_0=\epsilon \mu^{-1}$.
\end{proposition}

On the other hand, if $q=1$ skipping unremovable cubic terms in $\cA$ comes with a price of $O(\mu^{\frac{1}{2}}h^{1-d}+\mu^{-\frac{1}{2}}h^{\frac{1}{2}-d})$ error. If we would prefer not to remove them, we cannot treat separate equations.

\subsection{Case of constant $f_j$}
\label{sect-20-4-3-4}

However we can do as well as before assuming that 
\begin{equation}
f_j=\const, \quad j=1,\ldots,r
\label{20-4-81}
\end{equation}
which by no means excludes variable $g^{jk}$, $F_{jk}$ or unremovable cubic terms.

Really, under assumption (\ref{20-4-81}) we can introduce  
\begin{equation}
\gamma = \inf \bigr\{t: \vartheta (t)\ge \rho^2,
\vartheta (t)t^{-1}\ge |\nabla V (x'',x''',\xi'')|\bigr\}+\bar{\gamma},
\label{20-4-82}
\end{equation}
and $\varrho$ by (\ref{20-4-74}) and $\bar{\varrho}$, $\bar{\gamma}$ by (\ref{20-4-75}). Then we can apply our standard arguments as long as 
\begin{equation}
\rho^2+ |\nabla V| \ge C\mu^{-1}
\label{20-4-83}
\end{equation}
which will be fulfilled automatically for $\gamma \ge C\bar{\gamma}$ as long as
\begin{equation}
\bar{\rho}\ge C\mu^{-\frac{1}{2}}\iff 
\mu \ge C h^{-\frac {2l}{l+2} }|\log h|^{\frac {2\sigma}{l+2}}
\label{20-4-84}
\end{equation}
with the right-hand expression exceeding $h^{-\frac{1}{2} }$ for sure. 

Therefore, under condition (\ref{20-4-84}) estimate (\ref{20-4-80}) holds. On the other hand, if condition (\ref{20-4-84}) is violated, we need also consider elements with violated condition (\ref{20-4-83}), in particular, with $|\nabla V|\le C\mu^{-1}$ and
\begin{equation}
\gamma \le \hat{\gamma}=C\mu^{-\frac {1}{l} }|\log h|^{\frac {\sigma}{l} },
\label{20-4-85}
\end{equation}
Instead of $\gamma$ partition we consider $\max(\gamma, \hat{\gamma})$-partition 
and on  elements with $\gamma=\hat{\gamma}$ let us  introduce 
\begin{align}
&\gamma'=\hat{\gamma}\mu \rho^2+\bar{\gamma}', \qquad &&\varrho'= \rho+\bar{\rho}', \label{20-4-86}\\[2pt]
&\bar{\gamma}'= \mu^{\frac{1}{3}} \hat{\gamma}^{\frac{1}{3}} h^{\frac{2}{3}}=
h^{\frac{2}{3}} \mu^{\frac{l-1}{3l}}|\log h|^{\frac {\sigma} {3l}},\qquad
&&\bar{\rho}'=\mu^{-\frac{1}{3}} \hat{\gamma}^{-\frac{1}{3}} h^{\frac{1}{3}}
\label{20-4-87}
\end{align}
Then we arrive to the remainder estimate 
$Ch^{1-d}+C \mu h^{2-d}{\bar\gamma}^{\prime\,-1}$:

\begin{proposition}\label{prop-20-4-14} 
Let $q=1$, $\mu \ge \bar{\mu}^*_{(1)}$ and condition \textup{(\ref{20-4-81})} be fulfilled. Then 

\begin{enumerate}[label=(\roman*), fullwidth]
\item Under condition \textup{(\ref{20-4-84})} estimate \textup{(\ref{20-4-80})} holds.

\item If condition \textup{(\ref{20-4-84})} is violated then estimate
\begin{equation}
\R^\T_{1\,Q}\le  Ch^{1-d} +
C \mu h^{1-d+ \frac{l} {l+2} }|\log \mu|^{- \frac {\sigma} {(l+2)} }+
C\mu ^{\frac {2l+1} {3l} } h^{\frac{4}{3}-d}|\log \mu|^{-\frac{\sigma} {3l} }
\label{20-4-88}
\end{equation}
holds.
\end{enumerate}
\end{proposition}

\begin{remark}\label{rem-20-4-15}
\begin{enumerate}[label=(\roman*), fullwidth]
\item Condition (\ref{20-4-84}) means exactly that the third term in the right-hand expression of (\ref{20-4-88}) does not exceed the second one.

\item If condition (\ref{20-4-81}) is violated one needs to take in account that 
as the main part is of the form 
$\xi_1^2+ \sum_{1\le j\le r} f_j^\w (\mu^2 x_{1+j}^2+h^2D_{1+j}^2) +V^\w$
with $f_j=f_j(x'',x_1,\xi'')$, $V=V(x'',x_1,\xi'')$, its gradient with 
respect to $x'',x_d,\xi''$ depends on the localization with respect to
$\bigl\{(\mu^2 x_{1+j}^2+h^2\xi_{1+j}^2),j=1,\ldots,r\bigr\}$ and there is at least one obstacle to it: due to the logarithmic uncertainty principle these quantities are defined with $C\mu h|\log h|$ precision. 

Therefore one must take 
$\hat{\gamma}^{l-1}|\log \hat{\gamma}|^{-\sigma}|\ge C\mu h|\log h|$ in the arguments leading to proposition \ref{prop-20-4-14} which deteriorates estimate (\ref{20-4-80}). It does not look worth needed efforts.

\item Another approach would be to drop unremovable cubic terms as 
$\rho \ge \hat{\varrho}$ paying the price of $O(h^{2-d}\hat{\varrho}^{-1})$ for this and apply the same arguments as above. Then the remainder estimate is
$C\mu h^{\frac{4}{3}-d} \hat{\varrho}^{\frac{2}{3}} \hat{\gamma}^{\frac{1}{3}}$
where $\hat{\gamma}= \hat{\varrho}^{\frac {2} {l}} |\log h|^{\frac{\sigma}{l} }$
and one should optimize the answer with respect to $\hat{\varrho}$ which is not
$\mu^{- \frac{1}{2}}$ anymore.
\end{enumerate}
\end{remark}

We leave to the reader the following generalization of Problem~\ref{Problem-20-4-12}:

\begin{Problem}\label{Problem-20-4-16}
Under assumption (\ref{20-4-81}) (and, in particular, in the case of constant $g^{jk}$, $F_{jk}$) prove remainder estimate $\R^\T= O(h^{1-d})$ under non-degeneracy assumption (\ref{20-2-42}).
\end{Problem}

\section{Improved remainder estimates} 
\label{sect-20-4-4-4}

There are peculiar cases when without microhyperbolicity or  non-degeneracy assumptions one can improve general remainder estimate $\R^\T_{1\,Q}\le Ch^{1-d}+C\mu h^{1+\frac{q}{3}-d}$ or estimates (\ref{20-4-80}), (\ref{20-4-88}). It happens only for $r\ge 2$.

Note that the factor estimating the number of indices one needs to take in account is actually 
$\bigl(\mathbf{n}(\tau , \mu h) - \mathbf{n}(\tau -C\rho^2,\mu h)\bigr)$
where $\mathbf{n}(\tau , \mu h)$ is an eigenvalue counting function for operator
${\mathbf a}$ in $\sL^2(\bR^r,\bC)$ calculated at any point $(x'',x''',\xi'')$ at $\supp Q$ while so far we used estimate 
$C\bigl(\rho^2 + \mu h\bigr)(\mu h)^{-r}$ for this difference.

Therefore as $(l,\sigma)\preceq (2,0)$ one can replace term 
$C\mu h^{1-d}\bar{\varrho}^q$ which we have in the framework of proposition~\ref{prop-20-4-13} by
\begin{multline}
C\mu ^r h^{r-d+1} \times \\
\shoveright{\int \int_{\bar{\varrho}} ^{\rho^*_1}
\bigl(\mathbf{n}(x,\tau , \mu h) - \mathbf{n}(x,\tau -C\rho^2,\mu h)\bigr) 
\rho ^{q-1}\gamma^{-1} \,d\rho \,dx +}\\
C\mu ^{r+1} h^{r-d} \int  \bigl(\mathbf{n}(x,\tau , \mu h) - 
\mathbf{n}(x,\tau -C\bar{\varrho}^2,\mu h)\bigr) \bar{\varrho} ^q \,dx\qquad
\label{20-4-89}%
\end{multline}
with $\gamma=\varrho ^{2/l} |\log \rho |^{\sigma/l}$, 
$\bar{\varrho}=h^{l /(l+2) }|\log h|^{-\sigma/(l+2) }$ where we already returned to coordinates $x$.

Also, in the general settings one can replace term $C\mu h^{1+\frac{q}{3}-d}$ by expression (\ref{20-4-89}) with $\gamma = \rho^2$ and 
$\bar{\varrho}=h^{\frac{1}{3}}$. 

Consider two examples

\begin{example}\label{ex-20-4-17} 
Consider eigenvalue counting function  $\mathbf{n}_0=\mathbf{n}_0(\tau, \hbar)$ of operator $\mathbf{a}_0=\cA_0$ in $\sL^2(\bR^r,\bC)$; since $\cA_0$ depends on $(x'',x'',\xi'')$, so does  $\mathbf{n}_0$ and we can always return to variable $x$ by map $\Psi_0$. Since $\cA_0$ contains only quadratic terms one can calculate $\mathbf{n}_0$ explicitly
\begin{equation}
\mathbf{n}_0(\tau , \hbar ) = \#\bigl\{\alpha \in \bZ^{+\,r}: 
\sum_j  (2\alpha_j+1)f_j\hbar + V< \tau\bigr\}.
\label{20-4-90}
\end{equation}
If all $f_1,\dots,f_r$ are commensurable then 
\begin{equation*}
\mathbf{n}_0(\tau,\hbar)-\mathbf{n}_0(\tau-\lambda,\hbar)\le 
C\hbar ^{1-r}(\lambda+\hbar)
\end{equation*}
is the best possible estimate ($\lambda \ge 0$); otherwise from
\begin{multline}
\mathbf{n}_0(\tau,\hbar)-\mathbf{n}_0(\tau-\lambda,\hbar)=\sum_{n\ge 0}
\#\bigl\{\alpha'\in \bZ^{+\,r-1}: 
\sum_{j\ge 2} {\frac {f_j}{f_1}} (2\alpha_j+1)\hbar \in \\
[f_1^{-1}(\tau -V) -2(n+1)-f_1^{-1}\lambda, f_1^{-1}(\tau -V) -2(n+1))\bigr\}
\label{20-4-91}
\end{multline}
one can derive a better estimate $C\hbar^{-r}(\lambda+\nu(\hbar))$ with $\nu(\hbar)=C\hbar^\kappa$, $\kappa>1$
depending on Diophantine properties of $\{f_2/f_1,\dots,f_r/f_1\}$ (see Harman~\cite{harman}).
\end{example}

\begin{example}\label{ex-20-4-18} 
Let us assume that\footnote {\label{foot-20-16} Actually we need a quantitative version of this condition and only at points where $\nabla V$ is a linear combination of $\nabla f_1,\ldots, \nabla f_r$.} 
\begin{equation}
\rank\{\nabla {{\frac {f_2}{f_1}}, \ldots, \nabla {\frac {f_r}{f_1}}}\}\ge \kappa-1\qquad\forall x.
\label{20-4-92}
\end{equation}

One can easily recover an integrated version of (\ref{20-4-91}), namely \begin{multline}
|\int_{B(0,1)} \bigl(\mathbf{n}_0(x, \hbar, \tau') -\mathbf{n}_0(x, \hbar, \tau')\bigr)  \, dx| \le C\hbar^{-r} \bigl(|\tau-\tau'|+\nu (\hbar)\bigr)\\
\forall \hbar \in (0,1)\ 
\forall \tau,\tau':|\tau |\le \epsilon,|\tau'|\le\epsilon
\label{20-4-93} 
\end{multline}
with $\nu(\hbar)=C\hbar^\kappa$.  Really, one needs just to consider $\zeta$-admissible partition in $(\alpha \hbar, x)$-space with the scaling function
$\zeta = \epsilon |\nabla (\tau-V-\sum_j(2\alpha_j+1)f_j\hbar)|$.

Furthermore, these arguments show that for $\mathbf{n}$ the left-hand expression of (\ref{20-4-93} ) does not exceed 
$C\hbar^{1-r} \bigl(\lambda +\hbar^\kappa  +\mu^{-\kappa} \bigr)$.
\end{example}

Then we arrive to

\begin{proposition}\label{prop-20-4-19} 
Let condition \textup{(\ref{20-4-93})} hold for $\mathbf{n}_0$ which is an eigenvalue counting function for $\cA_0$. Let 
$\bar{\mu}^*_{(q)}\le \mu \le \epsilon (h\log h|)^{-1}$. 

\begin{enumerate}[label=(\roman*), fullwidth]
\item Let \underline{either} $q=2$ \underline{or} $q=1$ and $\cA$ contain no unremovable cubic terms\footnote{\label{foot-20-17}
Unremovable cubic terms containing $hD_1$ are allowed.}. Then 
\begin{equation}
\R^\T_{1\,Q} \le Ch^{1-d}+ C\nu (\mu h) h^{ \frac{q}{3}-d};
\label{20-4-94}
\end{equation}

\item Assuming instead that $\mathbf{n}$ which is an eigenvalue counting function for $\mathbf{a}$, satisfies \textup{(\ref{20-4-93})}, we  need  in (i)  neither assumption ``$\cA$ contains no unremovable cubic terms'' nor \textup{(\ref{20-4-93})} for $\mathbf{n}_0$;

\item Let $q=1$ and $\cA$ contain unremovable cubic terms and 
$\mu \ge (h |\log h|)^{-\frac{1}{2} }$. Then 
\begin{equation}
\R^\T_{1\,Q} \le Ch^{1-d}+ C\nu (\mu h) h^{\frac{1}{3}-d} + 
C\mu^{-1}h^{\frac{1}{3}-d}.
\label{20-4-95}
\end{equation}
\end{enumerate}
\end{proposition}

\begin{proof} Statements (i) and (ii) follow directly from the above arguments. In  (iii) one needs to notice that (\ref{20-4-93}) for $\mathbf{n}_0$ implies
\begin{equation}
\int \bigl(\mathbf{n}(x,\tau,\hbar )- \mathbf{n}(x,\tau-\lambda,\hbar)\bigr)\,dx \le 
C\bigl(\lambda +\mu^{-1}+\nu(\hbar)\bigr)\hbar^{-r}.
\label{20-4-96}
\end{equation}
Easy details we leave to the reader.
\end{proof}

\begin{remark}\label{rem-20-4-20}
\begin{enumerate}[label=(\roman*), fullwidth]
\item If $q=1$ and (\ref{20-4-93}) holds then skipping unremovable cubic terms leads to an error not exceeding 
\begin{equation}
C\nu (\mu h)\mu^{-\frac{1}{2}}h^{-d}+ C\mu^{-\frac{3}{2}}h^{-d}.
\label{20-4-97}
\end{equation}
\item  As $\mu \ge h^{-\frac{2}{3}}$ last terms in the right-hand expression of (\ref{20-4-95}) and in (\ref{20-4-97}) do not exceed $Ch^{1-d}$.
\end{enumerate}
\end{remark}

Similarly one can prove 

\begin{proposition}\label{prop-20-4-21} 
Let $q=1,2$ and $f_j$ have constant multiplicities. Let condition \textup{(\ref{20-4-93})} hold for $\mathbf{n}_0$. Let
$\bar{\mu}^*_{(q)}\le \mu \le \epsilon (h |\log h|)^{-1}$. 

\begin{enumerate}[label=(\roman*), fullwidth]
\item  Let \underline{either} $q=2$ \underline{or} $q=1$ and $\cA$ contain no unremovable cubic terms. Then 
\begin{equation}
\R^\T_{1\,Q} \le Ch^{1-d} + C \nu (\mu h) \cdot h^{-d + \frac {lq}{l+2} }
|\log \mu |^{- \frac {lq\sigma} {2(l+2)} };
\label{20-4-98}
\end{equation}
\item  Let $q=1$, condition \textup{(\ref{20-4-81})} be fulfilled, $\cA$ contain unremovable cubic terms, estimate  \textup{(\ref{20-4-93})} hold for $\mathbf{n}$ as well and
$\mu \ge (h|\log h|)^{- \frac{1}{2}}$. 

Then estimate \textup{(\ref{20-4-95})} holds as \textup{(\ref{20-4-84})} is fulfilled  and estimate
\begin{multline}
\R^\T_{1\,Q} \le \\ Ch^{1-d} + C \nu (\mu h) \cdot
\Bigl( h^{-d + \frac{l}{l+2} } |\log \mu |^{- \frac {l\sigma} {2(l+2)} }+
C\mu ^{-\frac {l-1}{3l} }h^{-\frac{2}{3}-d}|\log h|^{-\frac {\sigma}{3l} }\Bigr)
\label{20-4-99}
\end{multline}
holds as \textup{(\ref{20-4-84})} fails;

\item   Let $q=1$, condition \textup{(\ref{20-4-81})} be fulfilled, $\cA$ contain unremovable cubic terms and 
$\mu \ge \epsilon (h|\log h|)^{-\frac{1}{2}}$.  Then estimate
\begin{equation}
\R^\T_{1\,Q} \le Ch^{1-d} + C (\nu (\mu h)+\mu^{-1}\bigr) 
\cdot h^{-d + \frac {lq}{l+2} }
|\log \mu |^{- \frac {lq\sigma} {2(l+2)} };
\label{20-4-100}
\end{equation}
holds as \textup{(\ref{20-4-84})} is fulfilled and estimate
\begin{multline}
\R^\T_{1\,Q} \le Ch^{1-d} +\\  C \bigl(\nu (\mu h) +\mu^{-1}\bigr) \cdot
\Bigl( h^{-d + \frac {l} {l+2} } |\log \mu |^{- \frac {l\sigma} {2(l+2)} }+
C\mu ^{-\frac {l-1}{3l} }h^{-\frac{2}{3}-d}|\log h|^{-\frac{\sigma} {3l} }\Bigr)
\label{20-4-101}
\end{multline}
holds as \textup{(\ref{20-4-84})} fails.
\end{enumerate}
\end{proposition}

\chapter{Stronger magnetic field: estimates}
\label{sect-20-5}

In this section we derive Tauberian remainder estimates the remaining cases: the \emph{strong magnetic field\/}  $\epsilon (h|\log h|)^{-1}\le \mu \le \epsilon h^{-1}$,
\emph{very strong magnetic field\/}  $\epsilon h^{-1}\le \mu \le C_0 h^{-1}$ and 
\emph{super strong magnetic field\/}  $ \mu \ge C_0 h^{-1}$.

\section{Strong magnetic field}
\label{sect-20-5-1}

\subsection{General settings}
\label{sect-20-5-1-1}

Consider now the strong magnetic field  case 
\begin{equation}
\epsilon (h|\log h|^{-1})\le \mu \le \epsilon h^{-1}.
\label{20-5-1}
\end{equation}
Recall that the reduction to the canonical form was done with a greater 
\begin{equation}
\varepsilon = C\bigl(\mu ^{-1}h|\log h|\bigr)^{\frac{1}{2}}
\label{20-5-2}
\end{equation}
rather than $C\mu^{-1}$; furthermore, an outer zone is empty now (but intermediate zone still present) and $\bar{\rho}_1=C$ rather than $\bar{\rho}_1 = C\bigl(\mu h|\log h|\bigr)^{\frac{1}{2}}$.

As $q\ge 2$ throwing away $O(\mu ^{-1})$ terms brings an $O(h^{1-d})$ error of $h^{-d}\cN^\MW$. Then all the arguments of Subsection \ref{sect-20-4-1} and~\ref{sect-20-4-2} remain true and the Tauberian remainder estimate will be $O(h^{1-d})$ as $q\ge 3$, 
$O(\mu h^{\frac{5}{3}-d})$ as $q=2$ and $O(\mu h^{\frac{4}{3}-d})$ as $q=1$;
in the latter case contribution of the intermediate zone is estimated by $Ch^{1-d}|\log h|$ rather than by $Ch^{1-d}$. Therefore  in the general settings previous results hold without any modifications:

\begin{proposition}\label{prop-20-5-1} 
Propositions \ref{prop-20-4-1},  \ref{prop-20-4-2},  \ref{prop-20-4-4}--\ref{prop-20-4-9} and theorem  \ref{thm-20-4-3} remain valid for the strong magnetic field case \textup{(\ref{20-5-1})} as well.   
\end{proposition}

Therefore in what follows we need to consider cases $q=1,2$ only.

\subsection{Microhyperbolic case}
\label{sect-20-5-1-2}
We must distinguish two cases:

\begin{enumerate}[label=(\roman*), fullwidth]
\item The  microhyperbolicity direction $\ell$ does not depend on $\boldtau$. In this case all the arguments of Subsubsection~\ref{sect-20-4-3-1}.1 remain true and the remainder estimate is $O(h^{1-d})$. 

In particular, the case of $f_j$ having constant multiplicities is covered. Really, in this case microhyperbolicity condition means exactly that at each point $x$
\begin{equation}
\sum_{1\le j\le r_1} \tau_j \nabla \bigl( \frac {f_j} {V} \bigr) =0,
\quad \forall_j\ \tau_j\ge 0 \implies \forall_j\ \tau_j =0;
\label{20-5-3}
\end{equation}
then there exists a vector $\ell=\ell (x)$ such that
\begin{equation}
\langle \ell, \nabla \bigl({\frac {f_j} V}\bigr)\rangle > 0\qquad \forall j
\label{20-5-4}
\end{equation}
which exactly means microhyperbolicity in direction $\ell$.

\item Microhyperbolicity condition is fulfilled with $\ell$ depending on 
$\boldtau$. Then as in Subsection~\ref{book_new-sect-19-5-2} of \cite{futurebook} one needs to study partition of the phase space and  select $\epsilon$ small enough in the upper bound for $\mu$ in (\ref{20-5-1}); then all the previous results will still hold.
\end{enumerate}

So, we have proven

\begin{proposition}\label{prop-20-5-2}
In the strong magnetic field case proposition~\ref{prop-20-4-11} holds provided \underline{either} microhyperbolicity direction does not depend on $\boldtau$ \underline{or}  $\epsilon$ small enough in the upper bound for $\mu$ in \textup{(\ref{20-5-1})}.
\end{proposition}

\subsection{Case of constant $g^{jk}$, $F_{jk}$} 
\label{sect-20-5-1-3}

We leave to the reader the following 

\begin{Problem}%
\footnote{cf. Problem~\ref{Problem-20-4-12}.}\label{Problem-20-5-3}
Prove that in the strong magnetic field case (\ref{20-5-1}) $\R^\T = O(h^{-1})$ provided $g^{jk}$ and $F_{jk}$ are constant and non-degeneracy assumption~(\ref{20-2-42}) is fulfilled.
\end{Problem}

\subsection{Case of $f_j$ having constant multiplicities} 
\label{sect-20-5-1-4}

Assuming that there are no cubic terms we get a family of separate operators $\cA_\alpha$ and all the previous arguments hold bringing us remainder estimate
(\ref{20-4-80}). On the other hand, as $q=1$ dropping unremovable cubic terms
produces an error $C\mu^{\frac{1}{2}}h^{1-d}$ which is less than the right-hand expression of (\ref{20-4-80}) unless $l=2$, $\sigma \in (-4,0]$. In this case assuming (\ref{20-4-81}) one can prove easily (\ref{20-4-80}) as well. So we arrive to 

\begin{proposition}\label{prop-20-5-4}
\begin{enumerate}[label=(\roman*), fullwidth]
\item  In the strong magnetic field  case as $f_j$ have constant multiplicities estimate \textup{(\ref{20-4-80})} holds unless $q=1$, $(l,\sigma)\succeq (2,-4)$ and there are unremovable cubic terms; 

\item  In this exceptional case should include an extra term 
$C\mu^{\frac{1}{2}}h^{1-d}$ in the right-hand expression.
\end{enumerate}
\end{proposition}

\subsection{Case of constant $f_j$} 
\label{sect-20-5-1-5}

We leave to the reader the following 

\begin{proposition}\label{prop-20-5-5}
In the strong magnetic field case proposition~\ref{prop-20-4-14} holds.
\end{proposition}
and
\begin{Problem}%
\footnote{\label{foot-20-19} cf. Problem~\ref{Problem-20-4-16}.}\label{Problem-20-5-6}
Prove that in the strong magnetic field case (\ref{20-5-1}) $\R^\T = O(h^{-1})$ provided $f_j$ are constant and non-degeneracy assumption~(\ref{20-2-42}) is fulfilled.
\end{Problem}

\subsection{Improvement without microhyperbolicity} 
\label{sect-20-5-1-6}

Furthermore, under condition (\ref{20-4-93} ) for $\mathbf{n}_0$ we get remainder estimate (\ref{20-4-94})  in the framework of proposition~\ref{prop-20-4-19}(i) which is an extremely small improvement now.

Moreover, if $f_j$ have constant multiplicities we can recover  (\ref{20-4-98}) in frames of proposition \ref{prop-20-4-21} (i).

\begin{proposition}\label{prop-20-5-7}
In the strong magnetic field  propositions \ref{prop-20-4-19} and \ref{prop-20-4-21} hold.
\end{proposition}

\section{Very strong magnetic field}
\label{sect-20-5-2}

Consider now the very strong magnetic field case
\begin{equation}
\epsilon h^{-1}\le \mu \le C_0h^{-1}.
\label{20-5-5}
\end{equation}
Then the same arguments as before but in simplified form hold and one can prove easily 

\begin{proposition}\label{prop-20-5-8} 
In the very strong magnetic field case

\begin{enumerate}[label=(\roman*), fullwidth]
\item In the general settings remainder estimate 
\begin{equation}
\R^\T  \le C h^{1-d}\left\{\begin{aligned}
&1\qquad &&\text{as\ \ }q\ge 3,\\
&h^{-\frac{1}{3} } \qquad&&\text{as\ \ }q=2,\\
&h^{-\frac{2}{3} } \qquad&&\text{as\ \ }q=1
\end{aligned}\right.
\label{20-5-6}
\end{equation}
\item If $q=1,2$ and $f_j$ have constant multiplicities then the remainder estimate 
\begin{equation}
\R^\T_1 \le Ch^{1-d}+ C h^{-d+ \frac {lq}{l+2} }|\log h|^{-\frac {lq\sigma} {2(l+2)} }
\label{20-5-7}
\end{equation}
holds.
\end{enumerate}
\end{proposition}

Note that we can ignore cubic terms (which are $O(h)$ now). In the very strong magnetic field case we need to modify microhyperbolicity assumption 

\begin{definition}\label{def-20-5-9} 
In the case of the very strong magnetic field we call operator \emph{microhyperbolic\/} (on energy level $\tau$) if there exists vector 
$\ell= \ell (\bar{z})\in \bR^{2r+q}$ such that 
\begin{gather}
\mu^2 \sum_{j,k\ge q+1} \bigl(\ell (a_{jk}a_0^{-1}) \bigr)
\zeta _j^\dag \zeta _k\ge \epsilon_1\qquad \forall \zeta\in \bC^r
\tag{\ref{20-1-19}}\\
\shortintertext{as long as}
a= \mu^2 \tau ;
\label{20-5-8}
\end{gather}
\end{definition}
we also need to modify non-degeneracy assumption (\ref{20-2-42}) (see (\ref{20-5-10}) below).

\begin{proposition}\label{prop-20-5-10} In the very strong magnetic field case
\begin{enumerate}[label=(\roman*), fullwidth]
\item If $q=1,2$ and microhyperbolicity condition (see definition~\ref{def-20-5-9})  is fulfilled then the remainder estimate 
$\R^\T \le C\mu h^{2-d}$ holds.

\item If $f_j$ are constant and either microhyperbolicity condition
\begin{gather}
|V-\tau-\sum_j (2\alpha_j+1)\mu h f_j |+|\nabla V|\ge \epsilon_0
\label{20-5-9}\\
\intertext{or non-degeneracy condition }
|V-\tau-\sum_j (2\alpha_j+1)\mu h f_j |+ |\nabla V|\le \epsilon_0\implies |\det \Hess V|\ge \epsilon_0
\label{20-5-10}
\end{gather}
is fulfilled then the remainder estimate $\R^\T \le C\mu h^{2-d}$ holds.
\end{enumerate}
(in the latter case we assume that $(l,\sigma)\succeq (2,0)$).
\end{proposition}

\section{Superstrong Magnetic Field}
\label{sect-20-5-3}

The last and the easiest case to consider is the superstrong magnetic field case $\mu \ge C_0h^{-1}$. In this case we need to consider Schr\"odinger-Pauli operator. Then we arrive to a single operator
\begin{equation}
\cA_0= h^2\sum_{2r+1\le j,k\le d}D_jg ^{jk}(x'',x''',\mu^{-1}hD'')D_k + a_0(x'',x''',\mu^{-1}hD'')
\label{20-5-11}
\end{equation}
where we recall that 
\begin{equation}
a_0=W\circ \Psi_0,\qquad W= V-\sum_j f_j\mu h
\label{20-5-12}
\end{equation}
belongs to $\sF^{l,\sigma}$ uniformly.

Surely we get a system, but all other components $\sU_{\alpha,\beta}$ (with $(\alpha,\beta)\ne 0)$ could be expressed via $\sU_{0,0}$.

Then the principal part of asymptotics (as $|W|\le c$) is of magnitude 
$\mu^rh^{r-d}$. For the remainder estimate we have

\begin{proposition}\label{prop-20-5-11} 
In the superstrong magnetic field

\begin{enumerate}[label=(\roman*), fullwidth]
\item For $q\ge 3$ remainder estimate $\R^\T \le C\mu^rh^{r+1-d}$    holds;

\item  For $q=1,2$ and  microhyperbolicity condition \textup{(\ref{20-5-9})}
fulfilled,  the remainder estimate $\R^\T \le C\mu^r h^{r+1-d}$ holds;

\item For $q=2$ and  $(l,\sigma)\succeq (2,1)$  the remainder estimate $\R^\T \le C\mu^r h^{r+1-d}$ holds;

\item For $q=1$, $(l,\sigma)\succeq (2,1)$,  and non-degeneracy condition \textup{(\ref{20-5-10})} fulfilled the remainder estimate $\R^\T \le C\mu^r h^{r+1-d}$ holds;

\item For $q=2$, $(l,\sigma)\prec (2,0)$ and for $q=1$  the remainder estimate
\begin{equation}
\R^\T \le C\mu^r h^{r-d}\left\{\begin{aligned}
& h^{  \frac {2l} {l+2}  } 
|\log h|^{- \frac {2\sigma} {l+2} }\qquad &&\text{as\ \ }q=2 ,\\
&h^{ \frac {l} {l+2} } |\log h|^{- \frac {\sigma}{l+2}}\qquad &&\text{as\ \ }q=1 
\end{aligned}\right.
\label{20-5-13}
\end{equation}
holds.
\end{enumerate}
\end{proposition}

On the other hand, both   mollification and approximation errors in the operator are  $O(h)$; then if either $q\ge 2$ or microhyperbolicity or non-degeneracy assumptions are fulfilled then both mollification and approximation errors in the magnetic Weyl expression are $O(h)$ multiplied by $\mu^r h^{r-d}$ i.e. $\mu^r h^{r+1-d}$. 

Furthermore, if $q=1$ and we consider the general case, then both mollification and approximation errors are $O(h^{\frac{1}{2} })$ multiplied by  
$\mu^r h^{r-d}$ i.e. $\mu^r h^{r+ \frac{1}{2}-d}$. 

In all these cases both mollification and approximation errors do not exceed the remainder estimate. 

\chapter{Intermediate magnetic field: Calculations}
\label{sect-20-6}

Now we need to derive more explicit expressions rather than the Tauberian expression 
\begin{equation}
h^{-1}\int _{-\infty}^0 F_{t\to h^{-1}\tau}{\bar\chi}_T(t)\Gamma 
( U\psi_y \,^t\!Q_y)\,d\tau
\label{20-6-1}
\end{equation}
or the sum of those with different $T$ where $\psi_y=\psi(y)$, $Q_y=Q(y,hD_y)$.

\section{Weak magnetic field redone}
\label{sect-20-6-1}

Weak magnetic field case is when $\mu \le h^{\delta-1}$ and we derived asymptotics with the principal part given by (\ref{20-6-1}) with $Q=I$, $T= \epsilon \mu^{-1}$ and with some remainder estimate derived without canonical form reduction. 

Note that the mollified operator $\tilde{A}$ is microhyperbolic in the  direction 
$\langle \xi ,\partial_{\xi}\rangle$ due to condition (\ref{20-1-11}). Further, 
in the zone $\{|\xi'''|\ge \epsilon' \}$ this operator is a differential operator since we take there $\varepsilon = Ch|\log h|$ which does not depend on $|\xi'''|$.

On the other hand, in the zone $\{|\xi'''|\le {\frac{1}{2}}\epsilon\}$ operator $\tilde{A}$ is microhyperbolic in direction
$\langle \xi' ,\partial_{\xi'}\rangle + \langle \xi'' ,\partial_{\xi''}\rangle$
and it is a differential operator with respect to $x',x''$ since mollification parameter $\varepsilon = C\rho^{-1}h|\log h|$ does not depend on $\xi',\xi''$.

Therefore due to the standard results rescaled 
$F_{t\to h^{-1}\tau}\chi_T(t)\Gamma ( U\psi_y)$ is negligible as 
$|\tau|\le \epsilon'$ and $T\in [T_*,\bar{T}]$, $T_*=Ch|\log h|$, 
$\bar{T}_0=\epsilon \mu^{-1}$. Therefore, due to (\ref{20-2-29}) we can take $T=T_*$ in (\ref{20-6-1}) with only $O(h^s)$ difference. It happens in the proof of theorems~\ref{thm-20-2-10}, \ref{thm-20-2-15}, and \ref{thm-20-4-3}.

Then we can launch the successive approximation method with the unperturbed operator
\begin{multline}
\bar{A}=\sum _{j,k}\bar{P}_j\bar{g}^{jk}\bar{P}_k+\bar{V},\\
\bar{g}^{jk}=g^{jk}(y),\ \bar{V}=V(y),\ 
\bar{P}_k= hD_k-\mu V_k(y)-\mu\sum_j(\partial_kV)(y)(x_k-y_k)
\label{20-6-2}
\end{multline}
and plug it into (\ref{20-6-1}). In what we get then the first term does not exceed $CTh^{-1-d}=Ch^{-d}|\log h|$. Perturbation can be written as
\begin{equation*}
\sum B_j(x_j-y_j) + \mu \sum_{j,k}B_{jk} (x_j-y_j)(x_k-y_k)+\ldots
\end{equation*}
with $h$-differential operators $B_{\dots}$ and since each factor $(x_j-y_j)$ according to our standard approach  leads to an extra factor  $T$ in the estimate, one can see easily that  each next term in the successive approximations gains factor $Ch^{-1}T\bigl(T^2+\mu T^3 \bigr)
\le Ch|\log h|$\,\footnote{\label{foot-20-20} Recall  that factor $Th^{-1}$ comes from Duhamel principle.}.

Therefore as we are looking for $O(h^{1-d})$ error, only first two terms should be considered. The first term results exactly (\ref{20-6-1}) in with $U(x,y,t)$ replaced by $\bar{U}(x,y,t)$ and it is exactly 
\begin{equation}
h^{-d}\int \tilde{\cN}^\MW (y,\tau) \psi(y)\,dy
\label{20-6-3}
\end{equation}
with $O(h^s)$ error\footnote{\label{foot-20-21} After we replaced $U$ by $\bar{U}$ we can replace $T$ by $\infty$ with $O(h^s)$ error since then (\ref{20-2-29}) holds for any $T\ge \epsilon\rho$ with $O(h^sT^{-s})$ error and it will be exactly (\ref{20-6-3}).}.

The second term consists of two parts; one of them is generated by the perturbation
\begin{equation}
R_1=\sum _{j,k,m}\bar{P}_j (\partial_m g^{jk})(y)(x_m-y_m) \bar{P}_k+
\sum_m (\partial_mV)(y)(x_m-y_m)
\label{20-6-4}
\end{equation}
and it is obviously $0$, while the second part is generated by 
\begin{equation}
R_2=\sum _{j,k}\Bigl(
(P_j-\bar{P}_j) \bar{g}^{jk} \bar{P}_k + \bar{P}_j \bar{g}^{jk} (P_k-\bar{P}_k) + (P_j-\bar{P}_j) \bar{g}^{jk} (P_k-\bar{P}_k)\Bigr)
\label{20-6-5}
\end{equation}
and is obviously $O(h^{-d-2}\mu T^3)=O(h^{1-d})$. 

Therefore 
\begin{claim}\label{20-6-6}
In the case of the weak magnetic field approach we can replace the Tauberian expression (\ref{20-6-1}) by (\ref{20-6-3})
without deterioration of the remainder estimate\footnote{\label{foot-20-22} We need also to replace $\tilde{\cN}^\MW$ by $\cN^\MW$ but estimate of the approximation error is already done; in what follows we will not distinguish between $\cN^\MW$ and $\tilde{\cN}^\MW$ keeping in mind that we always must include an approximation error estimate in our final statement.}.
\end{claim}

\begin{remark}\label{rem-20-6-1}
This takes care of Theorems \ref{thm-20-2-10} and \ref{thm-20-4-3}--with $\R^\MW$ rather than $R^\W_{(\infty)}$ remainder. Also under microhyperbolicity or non-degeneracy conditions it takes care of Theorem \ref{thm-20-2-15} and Problem~\ref{Problem-20-2-17} respectively.
\end{remark}

\section{Decomposition}
\label{sect-20-6-2}

\subsection{Decomposition. Part I}
\label{sect-20-6-2-1}

Now let us consider  intermediate magnetic field 
$\bar{\mu}^*_{(q)}\le \mu \le \epsilon (h|\log h|)^{-1}$.

Then there is a \emph{standard zone\/}\index{zone!standard} where we use the already derived estimate for $\R^\T_Q$ with $T=\bar{T}$; this standard zone contains an outer zone  $\Omega_\out=\{|\xi '''|\ge \bar{\rho}_1\}$
but could be wider due to the microhyperbolicity or non-degeneracy assumptions. Recall that due to assumption (\ref{20-1-11}) in this standard zone we can replace $\bar{T}$ by $T_*=Ch|\log h|$.

The unexpectedly difficult (as $q=1$) problem is to join asymptotics derived in the standard zone with a cut-off operator $I-\cT\cQ\cT^*$ and in the remaining non-standard zone where we employ cut-off operator $\cT\cQ\cT^*$. To  tackle it better let us rewrite the formula we derived for the answer:
\begin{equation}
h^{-1}\sum_n \int _{-\infty}^0 F_{t\to h^{-1}\tau} \bar{\chi}_{T_n}(t)
\Gamma ( U\psi_y \,^t\!Q_{n\,y})\,d\tau 
\label{20-6-7}
\end{equation}
where $Q_n$ are appropriate elements of the partition of $I$ and $T_n$ are already chosen. The asymptotics with this principal part and different remainder estimates were derived in Section~\ref{sect-20-4}. Actually one can replace here $T_n$ by larger values without the affecting remainder estimate, but let us do a bit differently. Considering some term here with $Q=Q_n$ and $T=T_n$ one can rewrite it as 
\begin{equation}
h^{-1} \int_{-\infty}^0\sum_{0\le m <\infty} \phi_{m, L_m} ( \tau) \Bigl(F_{t\to h^{-1}\tau} \bar{\chi}_{T_n} (t)\Gamma 
\bigl(U\psi_y \,^t\!Q_y \bigr) \Bigr)\,d\tau
\label{20-6-8}
\end{equation}
where $\supp \phi_0\subset [-2,2]$, $L_0=C\bar{\rho}_1$, 
$\supp \phi_m\subset [-2,-1]$,  $L_m= 2^mL_0$ as $m\ge 1$  and $\phi_{m,L_m}$ form partition of 1:  $\sum_{0\le m <\infty} \phi_{m,L_m}(\tau)=1$
as $\tau \le 0$.

Note that one can rewrite each term with $m>0$ in (\ref{20-6-8})   as
\begin{equation}
\phi_{m, L_m} ( hD_t )\Bigl(\varphi_T (t)\Gamma 
\bigl(U \psi_y \,^t\!Q_y  \bigr)\Bigr)\Bigr|_{t =0}, \quad \text{with\ \ } \varphi=\bar{\chi},\; T=T_n.
\label{20-6-9}
\end{equation}
Let us consider (\ref{20-6-9}) with $\varphi =\chi$ and  $T\ge hL_m^{-1}$. Using arguments of Subsection~\ref{sect-20-4-1} one can estimate such term by 
\begin{equation}
Ch^{-d-1}\rho ^q L^2_m T \times \bigl( \frac {h} {T L_m} \bigr)^s 
\label{20-6-10}
\end{equation}
where factors $L_m$, $T$ are due to the integrations with respect to  $\tau$, $t$ respectively and another factor $L_m$ is due to the calculation of the number of contributing indices $\alpha$; we assume that 
\begin{equation}
L_m\ge C(\mu^{-1}+\rho^2).
\label{20-6-11}
\end{equation}
Further, the same $Ch^{-d}\rho ^q L_m $ estimate holds for (\ref{20-6-9}) terms with $\varphi={\bar\chi}$ and $TL_m=1$ as well. Then the sum with respect to $t$-partition results in $Ch^{-d}\rho ^q L_m $, and then the sum with respect to $\tau$ partition results in $C\rho^qh^{-d}$ because  $\max L_m\asymp 1$ here; terms with $L_m\ge C$ are less than $Ch^sL_m^{-s}$ due to the standard ellipticity arguments.

Let us employ a method of the successive approximations described in Subsection~\ref{sect-20-6-1} and plug it into (\ref{20-6-9}). Obviously,  estimate $C\rho^qh^{-d}$ for the first term is fulfilled.

A perturbation with factor $(x_j-y_j)$ will result  \underline{either} in the factor $CT^2h^{-1}$ if commuting with $A$ or $\bar{A}$, \underline{or} in the factor $CT\rho^{-1}$ if commuting with $Q$, \underline{or}, what is equivalent, in factors $hL_m^{-2}$ and $hL_m^{-1}\rho^{-1}$ respectively thus resulting  in the estimates
\begin{equation}
Ch^{1-d}\rho ^q L_m^{-1} \times \bigl(\frac{h} {T L_m} \bigr)^s \quad \text{and}\quad 
Ch^{1-d}\rho ^{q-1} \times \bigl(\frac{h} {T L_m} \bigr)^s
\label{20-6-12}
\end{equation}
respectively.

Note that the commutator with $Q_n$ (and thus the second of these expressions) appears in the sum  only as  $T\ge \epsilon T_n$; as for 
$T\le \epsilon T_n$  commutators from adjacent elements compensate one another.

Therefore taking the sum with respect to $t$-partition and then with respect to $m$ we get
\begin{equation}
Ch^{1-d}\rho ^q L_0^{-1},\qquad 
Ch^{1-d}\rho ^{q-1} \times \bigl( \frac{h} {T_n L_0} \bigr)^s
\label{20-6-13}
\end{equation}
respectively. We can always take $L_0=\bar{\rho}_1$ since $T_n\ge \mu^{-1}$ for sure and get $O(h^{1-d})$. As $q\ge 2$ we can do even better but it does not matter. Then both terms in (\ref{20-6-13}) do not exceed $Ch^{1-d}$ and remain this way after we sum over partition in zone $\{|\xi'''| \le \bar{\rho}_1\}$.

Therefore, in the terms with $m>0$ we can completely ignore partition and consider only the first term of the successive  approximations, Further,  after this we can take $T_n$ arbitrarily large (thus we can take $T_n=+\infty$); then after easy calculations we get
\begin{equation}
h^{-d} \int_{-\infty}^0 \int \bigl(1-\phi_{0, L_0} (\tau)\bigr) d_\tau 
\cN^\MW (x,\tau)\,\psi(x)\,dx.
\label{20-6-14}
\end{equation}
So we proved estimate
\begin{multline}
| h^{-1} \int_{-\infty}^0 \bigl(1-\phi_{0, L_0} ( \tau)\bigr) \sum_n\Bigl(F_{t\to h^{-1}\tau} \bar{\chi}_{T_n}(t)
(\Gamma U\psi_y \,^t\! Q_{n\,y} )\Bigr)\,d\tau - \\
h^{-d}\int_{-\infty}^0 \int \bigl(1-\phi_{0, L_0} (\tau)\bigr)\, d_\tau 
\cN^\MW (x,\tau)\,\psi(x)\,dx|\le Ch^{1-d}.
\label{20-6-15}
\end{multline}

Now we need to consider terms with $m=0$ i.e. 
\begin{equation}
h^{-1} \int_{-\infty}^0 \phi_{0, L_0} ( \tau) 
\Bigl(F_{t\to h^{-1}\tau} \bar{\chi}_T(t)\Gamma (U\psi_y\,^t\!Q_{ny}) \Bigr)\,d\tau, \qquad T=T_n.
\label{20-6-16}
\end{equation}
Instead we consider the same expression (\ref{20-6-16}) but with 
$T=\bar{T} =\epsilon \mu^{-1}$ and consider a correction later. We remember that $T_n\ge \bar{T}$  in the intermediate and interior zones. Thus,  let us  consider expression  (\ref{20-6-16}) with $T=\bar{T}$.

Now we want to replace $\bar{T}$ by a lesser value. We could replace  it by $T_*=Ch|\log h|$ with a negligible error before but it is not the case anymore since as a result of mollification and transformation the symbol of operator satisfies only an estimate
\begin{equation}
|\partial_\xi^\alpha a|\le C\bar{\rho}_1^{l-|\alpha|}
|\log \bar{\rho}_1|^{-\sigma}+C
\label{20-6-17}
\end{equation}
in the intermediate and interior zones. However then we can take there
\begin{equation}
T_*= C\bar{\rho}_1^{-1}h|\log h|\le \bar{T}
\label{20-6-18}
\end{equation}
Now let us apply the same successive approximations to calculate expression (\ref{20-6-16}) with $T=T_*$.

Then the perturbation with the factor $(x_j-y_j)$ results in the factor
$T^2h^{-1}=\bar{\rho}_1 ^{-2}h|\log h|^2$  and extra factor  $\bar{\rho}_1^{-2}|\log h|^2$ (in comparison with $h$) is well absorbed by $\rho^{q+1}$ at least as  $q\ge2$, $\mu \le h^{-1}|\log h|^{-2}$.

One can consider both cases $q=1$ and the exceptional case $\mu \ge h^{-1}|\log h|^{-2}$ either using rescaling arguments thus punishing  $T\ge h/\bar{\rho}_1$ by the factor $\bigl(h/(\bar{\rho}_1T)\bigr)^s$ or just by taking the two-terms approximation instead of the one-term and proving that the second term is identically $0$. 

Furthermore, a perturbation containing factor $\mu(x_j-y_j)(x_k-y_k)$ will get  an  extra factor $\mu T^3 h^{-1}=\mu h^2 |\log h|^3\bar{\rho}_1^{-3}$ which should be treated in the same way. We leave the easy standard details to the reader.

Therefore we can again consider only the first term of the successive approximations but in contrast to the previous part we cannot tend $T=T_n$ to infinity in this approximation term unless in the outer zone. Then the answer will be similar to one given in (\ref{20-6-14}), namely
\begin{multline}
h^{-1} \int_{-\infty}^0 \phi_{0, L_0} ( \tau) 
\Bigl(F_{t\to h^{-1}\tau} \bar{\chi}_{T_*}(t)\Gamma \bigl(\bar{U}\psi_y \bigr)\Bigr)\,d\tau =\\
h^{-d}\int_{-\infty}^0 \int \bigl(\phi (\tau)+ \omega (\tau)\bigr) 
d_\tau \cN^\MW (x,\tau)\,\psi(x)\,dx
\label{20-6-19}
\end{multline}
with
\begin{equation}
\omega(\tau) = T_* h^{-1}
\int_{-\infty}^0 \phi(\tau') ( \tau')
\hat{\bar{\chi}} \bigl( \frac {(\tau'-\tau)T_*} {h} \bigr)\,d\tau'-\phi(\tau), 
\label{20-6-20}
\end{equation}
where $\phi = \phi_{0, L_0}$. 

Here $\omega(\tau)\ne 0$ because we cannot replace $T_*$ by $+\infty$.
Instead we replaced back $T_*$ by $\bar{T}$ since operator $\bar{A}$ has the same propagation with respect to $x$ properties as $A$ and therefore
even with $\Gamma_y$ instead of $\Gamma$ the difference is negligible.

Let us define
\begin{multline}
\cN^\MW_\cQ (x,0) \Def\\
(2\pi)^{-d+2r} h^{r}\mu^r \sum_{\alpha\in \bZ^{+\, r}} 
\int_{\Omega_\alpha (x,\tau)} 
f_1\cdots f_r \, \frac {\sqrt g}{\sqrt g'''}  \, \cQ(\xi''')\,d\xi''',
\label{20-6-21}
\end{multline}
with
\begin{multline}
\Omega_\alpha (x,\tau)=\\
\bigl\{\xi: V+ \sum_j (2\alpha_j+1)\mu h + 
\sum_{1\le j,k\le q} g^{jk}\xi_j\xi_k+V\le \tau\bigl\},
\label{20-6-22}
\end{multline}
where $g^{jk}$ with $j,k= 1,\ldots, q$ denote coefficients after we reduce operator $\bar{A}$ to its canonical form 
\begin{equation}
\sum_{1\le j\le r} f_j(h^2D_{q+j}^2+ \mu^2x_{q+j}^2) + 
\sum_{1\le j,k\le q} g^{jk}h^2D_jD_k+V,
\label{20-6-23}
\end{equation}
and $g'''= \det (g^{jk})_{j,k=1,\ldots,q}^{-1}$.

Then $\cN^\MW = \cN^\MW_\cQ+ \cN^\MW_{1-\cQ}$. Plugging this into (\ref{20-6-14}), (\ref{20-6-19}) and adding we get in our calculation a candidate to the final answer
\begin{multline}
h^{-d} \int \cN^\MW(x,0)\psi(x)\,dx + \\
\shoveright{h^{-d} \int_{-\infty}^0 \int \omega (\tau) d_\tau \cN^\MW (x,\tau)\,\psi(x)\,dx \equiv}\\
h^{-d} \int \cN^\MW (x,0) \psi (x)\,dx + 
h^{-d} \int_{-\infty}^0 \int \omega (\tau) d_\tau \cN_\cQ^\MW (x,\tau)\, \psi(x) \,dx
\label{20-6-24}
\end{multline}
modulo negligible error  where $\cQ$ is supported in  
$\{|\xi'''|\le 2\bar{\rho}_1\}$  and equal $1$ in 
$\{|\xi'''|\le 2\bar{\rho}_1\}$.

Now it is a time to recall a correction
\begin{equation}
h^{-1} \sum_n \int_{-\infty}^0 \phi_{0, L_0} ( \tau) 
\Bigl(F_{t\to h^{-1}\tau}
\bigl( \bar{\chi}_{T_n }(t)- \bar{\chi}_{\bar{T}}(t)\bigr)
\Gamma  \bigl(U \psi_y \,^t\!Q_{n\,y}\bigr)\Bigr)\,d\tau 
\label{20-6-25}
\end{equation}
because we replaced $T_n$ by $\bar{T}$ in the intermediate and inner zones. Thus the following statement is proven:

\begin{proposition}\label{prop-20-6-2} 
In the framework of proposition \ref{prop-20-4-4}  expression \textup{(\ref{20-6-7})} with summation over all zones modulo $O(h^{1-d})$ is equal to
\begin{multline}
h^{-d} \int \cN^\MW (x,0)\psi (x)\,dx + \\
h^{-d} \int_{-\infty}^0 \int \omega (\tau) d_\tau \cN_\cQ^\MW (x,\tau)\,\psi(x)\,dx+\\
h^{-1} \sum_n \int_{-\infty}^0 \phi_{0, L_0} ( \tau) 
\Bigl(F_{t\to h^{-1}\tau}
\bigl( \bar{\chi}_{T_n}(t)- \bar{\chi}_{\bar{T} }(t)\bigr)\Gamma 
\bigl(\sU \tilde{\psi}_y \,^t\! \cQ_{n\,y} t\bigr)\Bigr)\,d\tau 
\label{20-6-26}
\end{multline}
where $\sU=\cT^* U\,^t\!\cT^{*}_y$, $\tilde{\psi}= \cT^*\psi\cT$, $\cQ_n=\cT^*Q_n\cT$
and here we can take $\cQ_n=\cQ_n(hD''')$ elements of partition in zone
$\{|\xi '''|\le 2\bar{\rho}_1\}$, $\cQ=\sum_n\cQ_n$.
\end{proposition}

Now what we need is to estimate the sum of the second and the third terms in the right-hand expression of (\ref{20-6-26}) (but not separately). 

\subsection{Decomposition. Part II}
\label{sect-20-6-2-2}

Let us consider the remaining inexplicit terms in the the the sum in expression (\ref{20-6-26}) 
\begin{equation}
h^{-1} \int_{-\infty}^0 \phi ( \tau) 
\Bigl(F_{t\to h^{-1}\tau}
\bigl(\bar{\chi}_{T_n} (t)-\bar{\chi}_{\bar{T}}(t)\bigr)
\Gamma \bigl(\sU \tilde{\psi}_y \,^t\!\cQ_y\bigr)\Bigr)\,d\tau 
\label{20-6-27}
\end{equation}
where we renamed $\phi_{0,L_0}$ into $\phi (\tau)$ (so, in (\ref{20-6-20}) $\phi$ is no more $\phi_{0,L_0}$) and  for the sake of simplicity we  skip index $n$. 

Let us apply the same decomposition technique to $\phi(\tau)$ as we did to 1:
$\phi (\tau)=\sum_{m\ge 0}\phi_{m,L_m}(\tau)$ with lesser $L_m$ than before and rewrite (\ref{20-6-27}) as the sum of (\ref{20-6-9})-type terms
\begin{equation}
h^{-1} \int_{-\infty}^0 \phi_{m, L_m} ( \tau) 
\Bigl(F_{t\to h^{-1}\tau}
\bigl(\bar{\chi}_{T_n}(t)-\bar{\chi}_{\bar{T}}(t)\bigr)
\Gamma \bigl(\sU \tilde{\psi}_y \,^t\!\cQ_y\bigr)\Bigr)\,d\tau.
\label{20-6-28}
\end{equation}
Replacing $\bigl(\bar{\chi}_{T_n}(t)-\bar{\chi}_{\bar{T}}(t)\bigr)$ by $\chi_T(t)$ we find that each obtained term with $m>0$ does not exceed (\ref{20-6-10}) again; however, there is a difference: after summation with respect to $T$ ranging from  $\bar{T}$ now to $T_n$ we get that  expression (\ref{20-6-28}) does not exceed 
\begin{equation}
C\rho ^q h^{-d}L_m \times \bigl( \frac {\mu h}{L_m} \bigr)^s
\label{20-6-29}
\end{equation}
as long $L_0\ge \mu h$ and $m>0$. After summation with respect to $L_m$ it gives us $C\mu \rho^q h^{1-d}$.

These arguments would enable us to set $L_0=\mu h$ albeit condition  
$L\ge \rho^2 +\mu^{-1}+\mu h$ prevents us. 

Let us assume now that $\mu \ge h^{-\frac{1}{2}}$; then we can take 
\begin{equation}
L_0=C\rho^2+C\mu h;
\label{20-6-30}
\end{equation}
analysis in the case $q=1$, $h^{-\frac{1}{3}}\le \mu \le h^{- \frac{1}{2}}$
will be done later with or without microhyperbolicity or non-degeneracy assumptions.

Let us employ the same method of the successive approximations with an unperturbed operator $\bar{\cA}=\cA_0(y'',y''',\mu^{-1}hD'',hD''')$ considered as an operator with operator valued symbol in $\sL^2(\bR^d)$. Then the perturbation is
\begin{multline}
\sum_{ 1\le j\le q+r} (x_j-y_j) a_j(x'',y'',x''',y''',\mu^{-1}hD'',hD''') + \\
b(x'',y'',x''',y''',\mu^{-1}hD'',hD''')
\label{20-6-31}
\end{multline}
with $b =O(\mu^{-1})$.

Consider the second term of the successive approximations plugged into (\ref{20-6-28}) with $m>0$ and transform it in our usual way. Then the factor $(x_j-y_j)$ in the first part of perturbation (\ref{20-6-28}) should commute with \underline{either} $\cG^\pm $ bringing thus factor $\rho T^2 h^{-1}$ as 
$1\le j\le q$ or the smaller factor $\mu^{-1} T^2 h^{-1}$ as $q+1\le j\le q+r$ \underline{or} with $\cQ$ bringing factor $\rho^{-1}T$. 

However again commuting with $\cQ$ in the final run (after summation with respect to $n$) should be considered only as $T\ge T_n$ and since 
$T_n\ge \epsilon h\rho^{-2}$ due to the analysis of Section~\ref{sect-20-4} the first factor is larger anyway.

Thus  from the first part of perturbation we get terms estimated by (\ref{20-6-10}) multiplied by one these factors: so we get
\begin{phantomequation}\label{20-6-32}\end{phantomequation}
\begin{align}
&Ch^{-d}\rho ^q L\times \bigl( \frac {h} {T L} \bigr)^s \times \rho T^2 h^{-1},
\label{20-6-32-1}\tag*{$\textup{(\ref*{20-6-32})}_1$}\\
&Ch^{-d}\rho ^q L\times \bigl( \frac {h} {T L} \bigr)^s \times \mu^{-1} T^2h^{-1},
\label{20-6-32-2}\tag*{$\textup{(\ref*{20-6-32})}_2$}\\
&Ch^{-d}\rho ^q L\times \bigl( \frac {h} {T L} \bigr)^s \times T\rho^{-1}
\label{20-6-32-3}\tag*{$\textup{(\ref*{20-6-32})}_3$}
\end{align}
respectively because in our analysis $T\ge \epsilon\mu^{-1}$, $L\ge \mu h$.

On the other hand, the second part of the perturbation (\ref{20-6-31}) brings factor  $C\mu^{-1}Th^{-1}$ or equivalently $C\mu^{-1}L^{-1}$ thus giving us 
\begin{equation}
C\mu^{-1}\rho ^q h^{-d}\times \bigl( \frac h {T L} \bigr)^s.
\label{20-6-33}
\end{equation}

Summation of \ref{20-6-32-1}--\ref{20-6-32-3}, (\ref{20-6-33}) with respect to $L$, $T$ results in the same expressions but with the minimal possible values, i.e. $T=\mu^{-1}$ and $L=\rho^2+\mu h$; so the result of the summation for each of the expressions $\textup{(\ref{20-6-32})}_{1-3}$ does not exceed
\begin{equation}
C\mu^{-1}\rho ^{q+1} h^{-d} \times\bigl( \frac {\mu h} {\rho^2+\mu h} \bigr)^s,
\label{20-6-34}
\end{equation}
while the result for (\ref{20-6-33}) does not exceed 
\begin{equation}
C\mu^{-1}\rho ^q h^{-d}\times\bigl( \frac {\mu h} {\rho^2+\mu h} \bigr)^s.
\label{20-6-35}
\end{equation}
Integration over $d\rho/\rho$ brings instead of (\ref{20-6-34}), (\ref{20-6-35}) their values as $\rho =\bar{\rho}^*_1$; the first one is always $O(h^{1-d})$ while the second one is $O(h^{1-d})$ as $q\ge 2$ only and is
$O(\mu^{-\frac{1}{2} }h^{ \frac{1}{2}-d})$ as $q=1$. We can do better than this but marginally. We will do it later.

Therefore, replacing $\sU$ by its one-term approximation 
\begin{equation}
\bar{\sU}=-i h^{-1} \bar{\cG}^+\updelta (x-y)\updelta(t)+
i h^{-1} \bar{\cG}^-\updelta (x-y)\updelta(t)
\label{20-6-36}
\end{equation}
in each term of (\ref{20-6-28}) with $m>0$ brings the total error $O(h^{1-d})$ unless  $q=1$ and $\cA$ contains (unremovable) cubic terms; in the latter case the total error is $O(\mu^{- \frac{1}{2}} h^{ \frac{1}{2}-d})$.

Thus, let us plug $ \bar{\sU}$ into expression (\ref{20-6-28}). Recall that the symbol of $\tilde{\psi}$ is given by (\ref{20-3-34}). But then the remainder estimate $O(\mu^{-1}|\log h|^{-1})$ in this expression should be multiplied by $C\bar{\rho}_1\rho ^q h^{-d}$  (recall that an extra factor $\bar{\rho}_1$ comes from the support of $\phi$) and it results in $O(h^{1-d})$. 

Furthermore, terms containing factors $\eta_j^\w$ with $j=q+1,\ldots,d$ result in $0$  because these terms applied to $v(x'',x''')\Upsilon_\alpha(x')$ result in 
\begin{equation*}
\sum_{\beta:|\alpha-\beta|=1}v_\beta (x'',x''')\Upsilon_\alpha(x')
\end{equation*}
while $\cA_0$ ``honors'' $|\alpha|$.

Finally, terms containing $\eta_j^\w$ with $j=1,\ldots,q$ also result in $0$ since we can take $\cQ_n$ even with respect to $\xi'''$ and $\cA_0$ is even as well. 

Therefore it suffices to plug only $\psi_0$ instead of $\tilde{\psi}$. We can also replace $T_n$ by any larger value, say infinity.  But then using the change of the coordinates $\Psi_0$  we get exactly the second term of (\ref{20-6-26}) with the opposite sign and with $\omega$ defined by (\ref{20-6-20}) for $\phi=\phi_{m,L_m}$. Then we arrive to (\ref{20-6-26}) with $L_0=\mu h+\rho^2$ (which now depends on $n$).

We also need to look at the case $m=0$. So far we considered both intermediate and inner zones. Consider now intermediate zone only. Then all the above arguments remain true with the only difference that factor 
$\bigl( h/(\rho^2 T)\bigr)^s$ is now due to estimate (\ref{20-4-4}) rather than $\tau$-mollification. 

Therefore we conclude that 

\begin{enumerate}[label=(\roman*), fullwidth]
\item If \underline{either} $q\ge 2$ \underline{or} $q=1$ and $\cA$ contains no (unremovable) cubic terms then the expression
\begin{multline}
|\sum_n h^{-d}\int_{-\infty}^0 \int \omega (\tau) d_\tau \cN_{\cQ_n}^\MW (x,\tau)\,\psi(x)\,dx+\\
h^{-1} \sum_n \int_{-\infty}^0 \phi_{0, L_0} ( \tau) 
\Bigl(F_{t\to h^{-1}\tau}
\bigl( \bar{\chi}_{T_n}(t)- \bar{\chi}_{ \bar{T} }(t)\bigr)\Gamma \sU \tilde{\psi}_y \,^t\!\cQ_{n\,y} \Bigr)\,d\tau | 
\label{20-6-37}
\end{multline}
with summation over $n$ such that $\cQ_n$ belongs to the intermediate zone does not exceed $Ch^{1-d}$; 

\item If $q=1$ and $\cA$ contains (unremovable) cubic terms  then this expression (\ref{20-6-37}) does not exceed 
$C\mu^{- \frac{1}{2}}h^{ \frac{1}{2}-d}$ 
\end{enumerate}

Combining  proposition~\ref{prop-20-6-2} we arrive to 

\begin{proposition}\label{prop-20-6-3} 
In the framework of proposition \ref{prop-20-4-4} let us assume that
$h^{- \frac{1}{2}}\le \mu \le \epsilon (h |\log h|)^{-1}$. 

Then one can rewrite expression \textup{(\ref{20-6-7})} as \textup{(\ref{20-6-26})} with  $L_0= C\mu h$ and $\phi$, $\omega$ defined accordingly and one should take $\cQ_n=\cQ_n(hD''')$ elements of partition in the inner zone  $\{|\xi '''|\le \rho^*_1\}$, $\cQ=\sum_n\cQ_n$
and the error is 
\begin{enumerate}[label=(\roman*), fullwidth]
\item $O(h^{1-d})$ if \underline{either} $q\ge 2$ \underline{or} $q=1$ and $\cA$ contains no (unremovable) cubic terms;
\item $O(\mu^{- \frac{1}{2}}h^{\frac{1}{2}-d})$ if $q=1$ and $\cA$ contains (unremovable) cubic terms.
\end{enumerate}
\end{proposition}
Recall that the only cubic terms which matter are those which do not contain $hD_j$ with $1\le j\le q$.

\section{Inner zone: general settings}
\label{sect-20-6-3}

Now in virtue of proposition \ref{prop-20-6-3} we need to consider only an inner zone $\Omega_\inn=\{|\xi'''|\le \bar{\rho}^*_1\}$. In the general setting there is no point to take $L_0\le \mu h$ because we cannot improve estimate 
$O\bigl((\mu h)^{-r}\bigr)$ for the number of indices making contributions as $\tau$ runs an interval of length $L_0$. Because of this (as $m=0$) one should replace $\textup{(\ref{20-6-32})}_{1-3}$ by 
\begin{phantomequation}\label{20-6-38}\end{phantomequation}
\begin{align}
&C\mu \rho^q h^{1-d}\times 
\min\Bigl(\bigl( \frac {h} {T\rho^2} \bigr)^s,1\Bigr)\times \rho T^2h^{-1},
\tag*{$\textup{(\ref*{20-6-38})}_1$}\label{20-6-38-1}\\
&C\mu \rho^q h^{1-d}\times 
\min\Bigl(\bigl( \frac{h} {T\rho^2} \bigr)^s, 1\Bigr) \times \mu^{-1} T^2h^{-1}, \tag*{$\textup{(\ref*{20-6-38})}_2$}\label{20-6-38-2}\\
&C\mu \rho^q h^{1-d}\times 
\min\Bigl(\bigl( \frac {h} {T\rho^2} \bigr)^s,1\Bigr) \times T\rho^{-1}
\tag*{$\textup{(\ref*{20-6-38})}_3$}\label{20-6-38-3}
\end{align}
respectively where factor $(h/{T\rho^2})^s$ as  $T\ge h/\rho^2$ is due to the estimate (\ref{20-4-4}). Furthermore, \ref{20-6-38-3} one should count with $T\ge \epsilon T_n \asymp \rho$ only.

After summation with respect to $T$ we get instead of expressions $\textup{(\ref*{20-6-38})}_{1-2}$ their values as $T=h/\rho^2$ and instead of expression $\textup{(\ref*{20-6-38})}_3$ its value as $T=\rho$ i.e. 
\begin{phantomequation}\label{20-6-39}\end{phantomequation}
\begin{equation}
C\mu \rho^{q-3} h^{2-d},\quad C \rho^{q-4}h^{2-d},\quad 
C\mu \rho^q h^{1-d}\bigl({\frac h {\rho^3}}\bigr)^s
\tag*{$\textup{(\ref*{20-6-39})}_{1-3}$}\label{20-6-39-*}
\end{equation}
respectively. 

Finally, after summation with respect to $\rho \ge h^{\frac{1}{3}}$ we get 
$O(h^{1-d})$ as $q\ge 3$, $O(\mu h^{\frac{5}{3}-d})$ as $q=2$, and
$O(\mu h^{ \frac{4}{3}-d})$ as $q=1$.

On the other hand, instead of (\ref{20-6-33}) we get
\begin{equation}
C\mu \rho ^q h^{1-d}\times
\min\Bigl(\bigl(\frac {h} {T\rho^2} \bigr)^s,1\Bigr) \times \mu^{-1}Th^{-1}
\label{20-6-40}
\end{equation}
and summation with respect to $T$ results in $C\rho^{q-2}h^{1-d}$.

Summation of this expression with respect to $\rho$ should be taken over 
$\rho\ge \max \bigl(\mu^{-\frac{1}{2}}, h^{\frac{1}{3}}\bigr)$ only; it results in  $O(h^{1-d})$ as $q\ge 3$, $O(h^{1-d}|\log \mu^2 h|)$ as $q=2$, and $O(\mu^{\frac{1}{2}} h^{1-d})$ as $q=1$. 

To get rid of the factor $|\log \mu^2 h|$ as $q=2$ one should consider two-term approximation; then (\ref{20-6-40}) will be replaced  by 
\begin{equation}
C\mu \rho ^q h^{1-d}\times
\min\Bigl(\bigl({\frac h {T\rho^2}}\bigr)^s,1\Bigr) \times \mu^{-2}T^2h^{-2}
\tag*{$\textup{(\ref*{20-6-40})}^*$}\label{20-6-40-*}
\end{equation}
which results in $C\mu^{-1}\rho^{-2}h^{1-d}$ after summation with respect to
$T$ and in $O(h^{1-d})$ after summation with respect to $\rho$. 

On the other hand, consider the second part of perturbation (\ref{20-6-31}); one can rewrite it as $b(y'',y''',\mu^{-1}hD'',hD''')$. Really, since this part contains unremovable cubic terms, symbol $b$ belongs to $\sF^{1,1}$ and thus freezing $x''=y''$, $x'''=y'''$ one makes an error which could be accommodated in the first part of (\ref{20-6-31}) and treated correspondingly. 

However, the second term in the approximation of $\sU$ is
\begin{equation}
\bar{\sU}' = - i h \bar{\cG}^+  \bar{B} \bar{\cG}^+ \updelta (x-y) \updelta(t) 
+ i h \bar{\cG}^- \bar{B}\bar{\cG}^- \updelta (x-y) \updelta(t).
\label{20-6-41}
\end{equation}
and plugging $\psi_0$ instead of $\tilde{\psi}$ results in the trace equal to $0$.

Really, operators $\cA_0$ and thus $\bar{\cG}^\pm $ transform 
$v(x'',x''')\Upsilon_\alpha (x')$ into the sum of terms $v_\beta(x'',x''')\Upsilon_\beta (x')$ with $|\beta-\alpha|=0,2$ while operator $\bar{B}'$ transforms $v(x'',x''')\Upsilon_\alpha (x')$ into the sum of terms $v_\beta(x'',x''')\Upsilon_\beta (x')$ with $|\beta-\alpha|=1,3$ where 
\begin{equation*}
\bar{b}'=\sum_{(j,k,m)\in \fJ_3\cup\fJ_1} b_{jkm}\eta_j\eta_k\eta_m ,\qquad
\bar{b}''=\sum_{(j,k,m)\in \fJ_2\cup\fJ_0}  b_{jkm}\eta_j\eta_k\eta_m
\end{equation*}
and $\fJ_s= (j,k,m)$ with exactly $s$ of $j,k,m$ in $\{q+1,\ldots,d\}$.

Meanwhile symbol $\bar{b''}$ is odd with respect to $\xi'''$ while everything else is even and thus the corresponding part of the second term vanishes after integration with respect to $\xi'''$. 

Furthermore, plugging $(\tilde{\psi}-\psi_0)$ instead of
$\tilde{\psi}$ produces an extra factor $\mu^{-1}$ in the estimate which is more than enough to compensate the logarithmic factor.

Now we need to estimate the contribution of the inner core
$\Omega^0_\inn = \{|\xi'''|\le \bar{\rho}_0\}$ with 
\begin{enumerate}[label=(\alph*)]
\item $\bar{\rho}_0= h^{\frac{1}{3}}$ if \underline{either} there are no unremovable cubic terms \underline{or} $\mu \ge h^{-\frac{2}{3}}$,

\item $\bar{\rho}_0= \mu^{-\frac{1}{2}}$ if there are unremovable cubic terms 
\underline{and} $\mu \le h^{- \frac{2}{3}}$. 
\end{enumerate}

However it follows from Section~\ref{sect-20-4} that we can consider this zone as a single partition element and set  $T_n=\bar{T}=\epsilon \mu^{-1}$ here
thus bringing an extra term 
$C\bar{\rho}_0^q \bar{T}^{-1} h^{1-d}=C\bar{\rho}_0^q \mu h^{1-d}$ into the remainder estimate; the latter expression is  
$O(h^{1-d})$ as $q\ge 3$, and $C\mu h^{ \frac{5}{3}-d} $ as $q=2$. 

As  $q=1$ this expression is equal to $C\mu h^{\frac{4}{3}-d}$ if either 
$\mu \ge h^{-\frac{2}{3}}$ or there are no unremovable cubic terms in $\cA$; otherwise  (if  $\mu \le h^{-\frac{2}{3}}$ and there are unremovable cubic terms in $\cA$)  this expression is equal to $C\mu^{\frac{1}{2}}h^{1-d}$ . 

Therefore with the error described one can replace $\sU$ by $\bar{\sU}$ in (\ref{20-6-27}). Furthermore, we can replace $\tilde{\psi}$ by $\psi_0$ there with $O(h^{1-d})$ error. 

Now we can preserve the remainder while increasing $T_n$ to $\infty$. Really, consider 
\begin{equation}
h^{-1} \int_{-\infty}^0 \phi_{0, L_0} ( \tau) 
\Bigl(F_{t\to h^{-1}\tau}\chi_T(t)\Gamma \bar{\sU}\psi_{0\,y} \,^t\!\cQ_y\Bigr)\,d\tau
\label{20-6-42}
\end{equation}
with $T\ge Ch\rho^{-2}$; it does not exceed 
\begin{equation}
Ch^{-d}\times \mu h\times \rho ^q \times \bigl(\frac{h} {\rho^2 T} \bigr)^s
\label{20-6-43}
\end{equation}
exactly as would be for $\sU$ but now we can take any $T\ge Ch\rho^{-2}$ because
$\bar{\cA}$ is a constant coefficient operator-valued operator. Summation with respect to $T\ge T_n$ results in the value of (\ref{20-6-43}) with $T=T_n$ which is  $\rho$ in these settings i.e. we get
\begin{equation}
C h^{1-d}\bigl(\rho^2+\mu h\bigr) 
\rho ^q \times \bigl( \frac {h} {\rho^3}  \bigr)^s
\label{20-6-44}
\end{equation}
and after summation with respect to $\rho$ we get its value as 
$\rho = h^{\frac{1}{3}}$ which is $C\mu h^{\frac{q}{3}+1-d}$.

But then after obvious calculations the second and the new third term in formula (\ref{20-6-26}) just cancel one another and we are left with just the first term. 

So we had proven 

\begin{proposition}\label{prop-20-6-4}  
In the framework of proposition~\ref{prop-20-4-4} let us consider
$h^{- \frac{1}{2}} \le \mu \le \epsilon (h|\log h|)^{-1}$. 

Then one can rewrite expression \textup{(\ref{20-6-7})} as
$h^{-d}\int \cN^\MW (x,0)\,\psi(x)\,dx$ with an error 
\begin{enumerate}[label=(\roman*),fullwidth]
\item $O(h^{1-d})$ as $q\ge 3$;
\item $O(h^{1-d}+ \mu h^{ \frac{5}{3}-d})$ as $q=2$;
\item $O(\mu h^{\frac{4}{3}-d})$ as $q=1$ and either $\mu\ge h^{-\frac{2}{3}}$
or $\cA$ contains no (unremovable) cubic terms;
\item $O(\mu^{\frac{1}{2}} h^{1-d})$ as $q=1$, $\mu\ge h^{-\frac{2}{3}}$ $q=1$ and $\cA$ contains (unremovable) cubic terms.
\end{enumerate}
\end{proposition}

\begin{remark}\label{rem-20-6-5} 
\begin{enumerate}[label=(\roman*),fullwidth]
\item When we say ``unremovable cubic terms'' we mean \emph{only} terms
$b_{jkm} \eta_j\eta_k\eta_m$ with $j\ge q+1$, $k\ge q+1$, $m\ge q+1$. Really, any other term bears an extra factor $\rho$ and if $\cA$ contains only these terms one can take $\bar{\rho}_0=Ch^{\frac{1}{3}}$ and the above analysis of Subsection~\ref{sect-20-6-2} and this Subsection produces an error $O(h^{1-d})$ (we get rid off factor $|\log \mu^2h|$ in the same way as we did it for $q=2$ in the proof of proposition \ref{prop-20-6-4}).

\item In the next Subsection~\ref{sect-20-6-4} we show that under either microhyperbolicity or non-degeneracy assumptions these unremovable cubic terms do not cause any corrections even as $q=1$. Further, in Subsection~\ref{sect-20-6-4} we also derive asymptotics with a better remainder estimate. 

\item   In Subsection~\ref{sect-20-6-5} we derive asymptotics with a correction term due to the unremovable cubic terms and with a better remainder estimate.
\end{enumerate}
\end{remark}

Proposition \ref{prop-20-6-4}  together with the remainder estimate of proposition \ref{prop-20-4-4} and with a mollification error estimate imply theorem \ref{thm-20-6-17} below.

\section{Improved error estimates}
\label{sect-20-6-4}

\subsection{Microhyperbolicity assumption}
\label{sect-20-6-4-1}

If either $q=2$ or $q=1$ and there are no unremovable cubic terms, everything is easy: under microhyperbolicity condition one can replace in the estimates $\textup{(\ref*{20-6-38})}_{1-3}$ factor $\mu h$ by factor $\rho^2$ as $\rho \ge \bar{\rho}_0=C h^{\frac{1}{3}}$; then instead of expressions \ref{20-6-39-*} we will get expressions
\begin{phantomequation}\label{20-6-45}\end{phantomequation}
\begin{equation}
C\rho^{q-1}h^{1-d},\qquad 
C\mu^{-1}\rho^{q-2}h^{1-d},
\qquad C\rho^{q+2}h^{-d}\bigl( \frac{h} {\rho^3} \bigr)^s
\label{20-6-45-*}\tag*{$\textup{(\ref*{20-6-45})}_{1-3}$}
\end{equation}
respectively which will result after integration over $d\rho/\rho$ in $O(h^{1-d})$ in each of them\footnote{\label{foot-20-23} Sure in the first expression as $q=1$ we get $Ch^{1-d}|\log h|$  but we consider two-term approximation as above  and also we replace ${\tilde\psi}$ by $\psi_0$; then the error term results in $O(h^{1-d})$ while the second term (\ref{20-6-41}) resulting in $0$.}. 

Further, contribution of the inner core is estimated in Section~\ref{sect-20-4} by $C\bar{T}^{-1}\bar{\rho}_0^{q+2}h^{1-d}= O(\mu h^{2-d})$ as we take 
$T_n= \bar{T}$ there. So we can again replace $\sU$ by $\bar{\sU}$ and then $\tilde{\psi}$ by $\psi_0$ and after this we replace $T_n$ by $\infty$. Then after obvious calculations we estimate  expression (\ref{20-6-28}) by $Ch^{1-d}$ and thus  expression (\ref{20-6-26}) will be reduced  modulo $O(h^{1-d})$ to its first term.

On the other hand, factor $\mu h$ is replaced by $\rho^2$ in the expression(\ref{20-6-40}) as well thus resulting after summation with respect to $T$ in $C\mu^{-1}\rho^q h^{-d}$; after summation with respect to $\rho$ we get 
$C\mu^{-\frac{1}{2}}h^{\frac{1}{2}-d}$ as $q=1$. Even use of two-term approximation gives us \ref{20-6-40-*} multiplied by $\rho^2/(\mu h)$ and after summation with respect to $T$ and then by $\rho$ we get 
$C\mu^{-2}h^{- \frac{1}{3}-d}$ which is $O(h^{1-d})$ only for 
$\mu \ge h^{-\frac{2}{3}}$.

However, if $\cA$ contains no unremovable cubic terms in the sense of remark \ref{rem-20-6-5}(i), one should replace there factor $\mu^{-1}T$ by 
$\mu^{-1}\rho T$; thus we get
\begin{equation*}
C\mu^{-1}T\rho^4 h^{-d}\min\Bigl(\bigl( \frac{h} {T\rho^2} \bigr)^s,1\Bigr)
\end{equation*}
in both cases instead of (\ref{20-6-40}). Then after summation with respect to $T$  we get $C\mu^{-1} \rho^2 h^{-d}$ and after integration over $d\rho/\rho$ we get $C\mu^{-1}\bar{\rho}_1^{*\,2}h^{-d}=C h^{1-d}$.

Also,  we need to cover also the 
case $(h|\log h|)^{-\frac{1}{2} } \le \mu \le h^{-\frac {1}{2} }$ 
(lesser $\mu$ are covered by the  weak magnetic field case due to the microhyperbolicity) but one can easily weaken an assumption 
$\mu\ge h^{-\frac{1}{2}}$ if $\cA$ contains no unremovable cubic terms.

So, we can plug $\bar{\sU}$ instead of $\sU$ (and $\psi_0$ instead of $\tilde{\psi}$) and thus we have proven

\begin{proposition}\label{prop-20-6-6}
Let $ (h|\log h|)^{- \frac{1}{2}}\le \mu \le \epsilon (h|\log h|)^{-1}$. Then in the framework of proposition \ref{prop-20-4-4}  under microhyperbolicity assumption (see definition~\ref{def-20-1-2}) one can rewrite expression \textup{(\ref{20-6-7})} as $h^{-d}\int \cN^\MW (x,0)\psi (x)\,dx$ with the $O(h^{1-d})$ error provided \underline{either} $q=2$ \underline{or} $q=1$ and $\cA$ contains no unremovable cubic terms (in the sense of remark \ref{rem-20-6-5}(i)).
\end{proposition}

However, we want to get rid off this extra assumption and thus we want to consider case ``$q=1$, $\cA$ contains unremovable cubic terms''. To do this we must reexamine an intermediate zone as well.

Our main idea is to estimate decay of 
$F_{t\to h^{-1}\tau}\chi_T(t)\Gamma \psi_0 \cQ U$ more accurately, using microhyperbolicity conditions it was done in Subsections~\ref{book_new-sect-19-4-1}--\ref{book_new-sect-19-4-4} of \cite{futurebook}.

\begin{proposition}\label{prop-20-6-7} 
Let $ (h|\log h|)^{- \frac{1}{2}}\le \mu \le \epsilon (h|\log h|)^{-1}$. Then in the framework of proposition \ref{prop-20-4-4} under microhyperbolicity condition (see definition~\ref{def-20-1-2})
\begin{multline}
|F_{t\to h^{-1}\tau} \chi_T(t) \Gamma \tilde{\psi}\cQ U|\le \\
CT\rho^{q+2} h^{-d}
\bigl(\frac{h}{T}\bigr)^{l-1} |\log \bigl(\frac{h}{T}\bigr)|^{-\sigma}\times 
\min\Bigl( \bigl( \frac {h} {\rho^2 T} \bigr)^s, 1\Bigr)
\label{20-6-46}
\end{multline}
as $\bar{\rho}_0= h^{-\frac{1}{3}}\le \rho \le \bar{\rho}_1$, 
$|\tau|\le \epsilon$.
\end{proposition}

\begin{proof} 
The proof repeats those of propositions~\ref{book_new-prop-19-4-3} and \ref{book_new-prop-19-4-4} of \cite{futurebook} and is based on $\eta$-approximation of $\cA$ with 
$hT^{-1}\le \eta \le ChT^{-1}|\log h|$. Easy details are left to the reader.
\end{proof}

After we improved our estimate by an extra factor 
$(h/T )^{l-1} |\log (h/T)|^{-\sigma}$ we can do more precise estimates.

Recall that
\begin{multline}
h^{-1}\int_{-\infty}^0 \phi_{0,L_0}(\tau) 
\Bigl(F_{t\to h^{-1}\tau}\chi_T(t)\tilde{\psi}\cQ_n \sU\Bigr)\,d\tau =\\
iT^{-1} \Bigl(F_{t\to h^{-1}\tau}\check{\chi}_T(t)
\tilde{\psi}\cQ_n \sU\Bigr)\Bigr|_{\tau =0}-
iT^{-1}L_0^{-1} 
\int_{-\infty}^0 \phi'_{0,L_0}(\tau) 
\Bigl(F_{t\to h^{-1}\tau}\check{\chi}_T(t)\tilde{\psi}\cQ_n U\Bigr)\,d\tau
\label{20-6-47}
\end{multline}
with $\check{\chi}(t)=t^{-1}\chi(t)$ and $\phi'_0$ derivative of $\phi_0$; recall that $\phi_0(0)=1$. Here $L_0$ could be larger than $\rho^2$ to accommodate a weaker assumption $\rho\le \mu^{-\frac{1}{2}}$. Then proposition \ref{prop-20-6-7} implies that for $q=1$ the right-hand expression of (\ref{20-6-47}) does not exceed 
\begin{equation}
C\rho^3 h^{-d}
\bigl(\frac{h}{T} \bigr)^{l-1} |\log \bigl( \frac{h}{T} \bigr)|^{-\sigma}\times 
\min\Bigl( \bigl( \frac {h} {\rho^2 T} \bigr)^s, 1\Bigr).
\label{20-6-48}
\end{equation}
Integrating this expression (\ref{20-6-48}) over $dT/T$ from $T=\bar{T}$ to $T=\infty$ we get this expression as $T=\bar{T}$ and then  integrating over $d\rho/\rho$ from $\bar{\rho}_)=h^{\frac{1}{3}}$ to $\rho =\bar{\rho}_1$
as $(h|\log h|)^{- \frac{1}{2}}\le \mu\le h^{- \frac{2}{3}}$
we get this expression as $\rho =\bar{\rho}_0$ \ i.e.
\begin{equation}
C (\mu h)^{l+  \frac{1}{2}} |\log h|^{-\sigma}h^{-d}
\label{20-6-49}
\end{equation}
which does not exceed an approximation error 
\begin{equation}
C\mu^{-l}\bar{\rho}_1|\log h|^{-\sigma} h^{-d}= C\mu^{ \frac{1}{2}-l}
h^{ \frac{1}{2}-d}|\log h|^{\frac{1}{2}-\sigma}
\label{20-6-50}
\end{equation}
as $\mu \le h^{- \frac{1}{2}}$. Therefore, in this case we can plug $\bar{\sU}$ instead of $\sU$ (and $\psi_0$ instead of $\tilde{\psi}$).

Thus it remains to consider the case  
$h^{-\frac {1}{2}}\le \mu \le h^{- \frac{2}{3} }$.
Consider expression (\ref{20-6-48}) with $T=\mu^{-1}$ and the last factor equal to $1$. Integrating it with respect to $d\rho/\rho $ from $\bar{\rho}_0=h^{- \frac{1}{3} }$ to $\rho=\mu^{- \frac{1}{2} }$, we get its value as 
$\rho= \mu^{-\frac{1}{2}}$ which is 
\begin{equation}
C \mu^{- \frac{3}{2}} (\mu h)^{l-1} |\log h|^{-\sigma}h^{-d}
\label{20-6-51}
\end{equation}
and which does not exceed approximation error as 
$l\le \frac{3}{2}$ and $Ch^{1-d})$ as $l> \frac{3}{2}$. These arguments estimate contribution of the inner core as well.  Therefore we estimated properly contribution of the inner zone $\{|\xi'''|\le C\mu^{- \frac{1}{2}}\}$.

Consider now $\rho$ ranging from $\mu^{-\frac{1}{2}}$ to $\bar{\rho}_1$.
We also apply the above arguments as $T\ge \mu h \ge  h/\rho^2 $, resulting after integration over $dT/T$ in 
\begin{equation}
C\rho^3\mu^{1-l}|\log h|^{-\sigma}\bigl( \frac{1} {\mu \rho ^ 2} \bigr)^s h^{-d}.
\label{20-6-52}
\end{equation}

As   $\bar{T}\le T\le \mu h$ we apply the two-term approximation method with the second term resulting in $0$ in the final answer and with an error estimate
\begin{equation}
C\rho^3 h^{-d}
\bigl( \frac{h}{T} \bigr)^{l-1} |\log \bigl( \frac{h}{T} \bigr)|^{-\sigma}\times 
\bigl( \frac{T} {\mu h} \bigr)^2\times 
\min\Bigl( \bigl( \frac{h} {\rho^2 T} \bigr)^s, 1\Bigr)
\label{20-6-53}
\end{equation}
which can be easily proven by arguments of propositions~\ref{book_new-prop-19-4-3}, \ref{book_new-prop-19-4-4}  and \ref{book_new-prop-19-4-15} of \cite{futurebook}. After integration over $dT/T$ we get 
\begin{equation}
C\rho^{2l-3} \mu^{-2}|\log \rho|^{-\sigma}h^{-d}\times 
\min\Bigl(\bigl( \frac {\mu h}{\rho^2} \bigr)^s,1\Bigr).
\label{20-6-54}
\end{equation}
Finally, integration of (\ref{20-6-52}) over $d\rho/\rho$ in the indicated limits results in its value at the lower limit $\rho =\mu^{- \frac{1}{2}}$; this result does not exceed  $\mu^{- \frac{1}{2}-l}|\log h|^{-\sigma}h^{-d}$ which again does not exceed an approximation error (\ref{20-6-50}). This is also true for expression (\ref{20-6-54}) as $l< \frac{3}{2}$; as $l=\frac{3}{2}$ we get 
$C\mu^{-2}|\log h|^{-\sigma}h^{-d}\times |\log \mu^2h|$
which also does not exceed approximation error.

So again we can replace $\sU$ by $\bar{\sU}$ and $\tilde{\psi}$ by $\psi_0$.
Moreover, after this replacement we can replace $T_n=\epsilon$ by $T_n=\infty$ with a negligible error. Thus we have proven

\begin{proposition}\label{prop-20-6-8} 
Let $q=1$ and 
$ (h|\log h|)^{- \frac{1}{2}}\le \mu \le \epsilon (h|\log h|)^{-1}$. Then in the framework of proposition \ref{prop-20-4-4}  under microhyperbolicity assumption (see definition~\ref{def-20-1-2}) one can rewrite expression \textup{(\ref{20-6-7})} as $h^{-d}\int \cN^\MW (x,0)\psi (x)\,dx$ with 
$O\bigl(h^{1-d}+
\mu^{\frac{1}{2}-l} h^{ \frac{1}{2}-d}|\log h|^{ \frac{1}{2}-\sigma}\bigr)$ error.
\end{proposition}

Proposition \ref{prop-20-6-8} together with remainder estimate of proposition \ref{prop-20-4-11} and with mollification error estimate imply theorem \ref{thm-20-6-15} below.

\subsection{Constant multiplicities of $f_j$}
\label{sect-20-6-4-2}

If we drop all junior terms we get a family of scalar operators and everything is easy. However while covering case $q=2$ this would not cover case $q=1$ even if the unremovable cubic terms contain $hD_1$. So we need more delicate arguments. 

Note that instead of  $\cQ_n(hD''')$ we have now $\cQ_n(x'',x''',\mu^{-1}hD'')\cQ'_m (hD''')$ 
where $\cQ'_m$ are elements of $\rho$-partition,  $\cQ_n$ are elements of $\gamma$-partition  with $\gamma,\varrho$  introduced by (\ref{20-4-73})--(\ref{20-4-74}).

\begin{enumerate}[label=(\roman*), fullwidth]
\item Consider first elements with $\varrho \ge C\bar{\varrho}$ and 
$\rho\ge \epsilon \varrho$.

Then one should multiply (\ref{20-4-49}) by $\gamma^d\times \varrho^2/\gamma$ where factor $\gamma^d$ comes from $\cQ_n$ and factor $\varrho^2/\gamma$ (which is greater than $\mu^{-1}$) comes from the estimate of $|\nabla_{x'',x''',\xi'''}A|$. 

Then expressions \ref{20-6-39-*} become
\begin{phantomequation}\label{20-6-55}\end{phantomequation}
\begin{multline}
C\mu \varrho^{q-1}h^{2-d}\gamma^{d-1}, \qquad C\varrho^{q-2}h^{2-d}\gamma^{d-1},\\
C\mu\varrho^{q+2}h^{1-d}
\min\Bigl(\bigl(\frac {h} {\varrho^3} \bigr)^s,1\Bigr)\gamma^{d-1}
\tag*{$\textup{(\ref*{20-6-55})}_{1-3}$}\label{20-6-55-*}
\end{multline}
respectively and the second and the third expressions do not exceed the first one.

\item On the other hand, as $\rho\le \epsilon \varrho$, 
$\varrho \ge C\bar{\varrho}$  we can rescale problem to the microhyperbolic one with $x\mapsto x\gamma$, $h\mapsto h/(\varrho\gamma)$, 
$\mu^{-1}h\mapsto \mu^{-1}h\gamma^{-2}$. Then the total contribution of elements with fixed $n$ and different $m$ does not exceed $\textup{(\ref{20-6-55})}_1$ again. therefore, $\textup{(\ref{20-6-55})}_1$ estimates the contribution of element $\cQ_n$.

After summation over partition $\textup{(\ref{20-6-55})}_1$ results in 
\begin{equation}
\mu h^{2-d}\int \varrho^{q-1}\gamma^{-1}\,dx''dx'''d\xi''
\label{20-6-56}
\end{equation}
which does not exceed the right hand expression in (\ref{20-6-89}) (or (\ref{20-4-80})).

\item 
Finally, the same estimate remains true as $\varrho \asymp \bar{\varrho}$
because now we estimate the contribution to each of the second and the third term of (\ref{20-6-7}) rather than to their sum.

\item
Similarly, if unremovable cubic terms contain factor $hD_1$ then we get 
that the contribution of $\cQ_n$ instead of (\ref{20-6-50})  is estimated by 
\begin{equation}
C\mu h^{1-d} \varrho^q\gamma^d \min \bigl(\mu^{-1}\varrho^{-1}, 1\bigr) \le Ch^{1-d}\varrho^{q-1}\gamma^d
\label{20-6-57}
\end{equation}
and the total contribution of all elements does not exceed $Ch^{1-d}$ because
then $\mu$ scales as $\mu \varrho$\,\footnote{\label{foot-20-24} Parameter $\mu^{-1}h$ scales independently from $\mu $ and $h$ because we deal with the reduced form of operator.}.

\item
On the other hand, if there are unremovable cubic terms without factor $hD_1$
then (\ref{20-6-57}) is replaced by
\begin{equation}
C\mu h^{1-d}\varrho^q\gamma^d \min \bigl(\mu^{-1}\varrho^{-2}, 1\bigr) \le Ch^{1-d}\varrho^{q-1}\gamma^d
\label{20-6-58}
\end{equation}
and summation results in $C\mu^{\frac{1}{2}}h^{1-d}$  and $Ch^{1-d}$ as $q=1,2$
respectively.
\end{enumerate}

Thus we arrive to

\begin{proposition}\label{prop-20-6-9}
Let $q=1,2$ and $f_j$ have constant multiplicities. 

\begin{enumerate}[label=(\roman*), fullwidth]
\item Let \underline{either} $q=2$ \underline{or} $q=1$ and there are no unremovable cubic terms. Then with an error not exceeding  the right-hand expression of \textup{(\ref{20-4-79})} one can rewrite  \textup{(\ref{20-6-7})} as  $h^{-d}\int \cN^\MW (x,0)\psi (x)\,dx$;

\item Let  $q=1$ and there are  unremovable cubic terms. Then with an error not exceeding  the right-hand expression of  \textup{(\ref{20-4-80})} plus $C\mu^{\frac{1}{2}}h^{1-d}$ one can rewrite  \textup{(\ref{20-6-7})} as 
$h^{-d}\int \cN^\MW (x,0)\psi (x)\,dx$.
\end{enumerate}
\end{proposition}

Proposition \ref{prop-20-6-9} together with remainder estimate of proposition \ref{prop-20-4-13} and with mollification error estimate imply theorem \ref{thm-20-6-23} below.

We leave to the reader  the following easy:

\begin{Problem}\label{Problem-20-6-10}
Let $f_j$ be constant and non-degeneracy condition (\ref{20-2-42})  be fulfilled. Prove that with $O(h^{1-d})$ error  one can rewrite  \textup{(\ref{20-6-7})} as  $h^{-d}\int \cN^\MW (x,0)\psi (x)\,dx$.
\end{Problem}

\subsection{Number theoretical arguments}
\label{sect-20-6-4-3}

Now it is easy to prove

\begin{proposition}\label{prop-20-6-11} 
In the framework of proposition \ref{prop-20-4-19}(i)
one can rewrite  expression \textup{(\ref{20-6-7})} as 
$h^{-d}\int \cN^\MW (x,0)\psi (x)\,dx$ with an error not exceeding the right-hand expressions of \textup{(\ref{20-6-81})}, \textup{(\ref{20-6-83})} below for $q=2,1$ respectively. 
\end{proposition}

Proposition \ref{prop-20-6-11} together with remainder estimate of proposition \ref{prop-20-4-19} and with mollification error estimate imply theorem \ref{thm-20-6-21}(i),(iii) below.

A bit more sophistication requires

\begin{proposition}\label{prop-20-6-12} 
In the framework of proposition \ref{prop-20-4-21}(i)
one can rewrite expression \textup{(\ref{20-6-7})} as
$h^{-d}\int \cN^\MW (x,0)\psi (x)\,dx$ with an error not exceeding the right-hand expressions of \textup{(\ref{20-6-89})}, \textup{(\ref{20-6-90})} below  as $q=2,1$ respectively.
\end{proposition}

\begin{proof} 
One needs just to replace $\tilde{\psi}$ by $\psi_0$, then consider two term approximations (the second term will result in 0 in the end), and then replace $T_n$ by $\infty$. In the intermediate zone estimates repeat those of the generic case. 

In the inner zone we need to consider integral (\ref{20-6-56}) and notice that contribution in it of all elements with $c^{-1}t\le \varrho  \le ct$, 
$t\le (\mu h)^{\frac{1}{2}}$ does not exceed
\begin{equation*}
Ch^{1-d}\varrho^{q-1}\gamma^{-1}\bigl(\varrho^2+\nu (\mu h)\bigr)
\end{equation*}
with $\varrho$ defined by given $\gamma$ and thus this integral does not exceed 
\begin{equation*}
Ch^{1-d}\int \varrho^{q+1}\gamma^{-1}\,dx + 
C\nu(\mu h)h^{1-d} \int \varrho^{q-1}\gamma^{-1}\,dx
\end{equation*}
and  plugging $\varrho=\gamma^{\frac{l}{2}}|\log \gamma|^{-\frac{\sigma}{2}}$ and replacing $dx$ by $\gamma^{-1}d\gamma$
we find that the first term does not exceed $Ch^{1-d}$ provided \underline{either} $q=2$ \underline{or} $q=1$, $(l,\sigma)\succ (1,1)$ while the second integral does not exceed 
\begin{equation*}
C\nu(\mu h)h^{1-d}\bar{\varrho}^{q-1}{\bar\gamma}^{-1}\asymp 
C\nu(\mu h)\bar{\varrho}^q h^{-d}
\end{equation*}
unless $q=2$, $l=2$ in which case the second term does not exceed 
$C\nu(\mu h)h^{1-d}|\log h|^{-\sigma}$.
\end{proof}

Proposition \ref{prop-20-6-12} together with remainder estimate of proposition \ref{prop-20-4-21} and with mollification error estimate imply theorem \ref{thm-20-6-25} below.

\section{$q=1$: asymptotics with correction}
\label{sect-20-6-5}

We want to improve remainder estimates in asymptotics as $q=1$ and $\cA$ has non-removable cubic terms. To do this we need to include some correction in
$h^{-d}\cN^\MW$. 

First of all, note that if both $\bar{\chi}$ and $\phi$ are proper functions\footnote{\label{foot-20-25}  In Section~\ref{book_new-sect-2-3} of \cite{futurebook} sense.} and $L_0\ge C\mu h|\log h|$, we can rewrite Tauberian expression (\ref{20-6-7}) in (\ref{20-6-26}) form with $O(h^s)$ error because (\ref{20-6-9}) with $\varphi_T=\chi_T $ and with $m>0$ is negligible as  $TL_m\ge Ch|\log h|$. 

To calculate the last term in (\ref{20-6-26}) we apply the same successive approximation method as before but with an unperturbed operator $\bar{\cA}=\cA_1(y'',y''',\mu^{-1}hD'',hD''')$ where in $\cA_1$ we include unremovable cubic terms 
\begin{equation}
\cB = \mu^2 \sum_{j,k,m,\ge 2}\bigl(b_{jkm}\eta_j\eta_k\eta_m\bigr)^\w
\label{20-6-59}
\end{equation}
but do not include the similar terms with at least one of the indices $j,k,m$ equal $1$\,\footnote{\label{foot-20-26} We also include $C\mu^{-2}\cA_0^2$ to have operator semibounded from below}; obviously such unremovable terms contain either one or three such indices. 

Then the last term in the perturbation (\ref{20-6-31}) will contain at least one factor $\eta_1$ and thus we get an extra factor $\rho$ in the corresponding estimates.  Then instead of an unwanted term 
$C\rho \mu h^{1-d}\times (\mu \rho^2)^{-2}$ calculated as 
$\rho =\max \bigl(\mu^{- \frac{1}{2}},\bar{\rho}_0\bigr)$ in an estimate,  we get $C\rho  \mu h^{1-d}\times (\mu \rho )^{-2}$, also  calculated as 
$\rho =\max \bigl(\mu^{-\frac{1}{2}},\bar{\rho}_0\bigr)$, which results in the error estimate $O(h^{1-d})$; here we take the same two-term approximation because we still need to consider the second term to avoid extra $|\log h|$ factor. 

However, we need to consider it only as $b$ contains exactly one factor $hD_1$ because otherwise we gain one more extra factor $\rho$. We also can replace $\tilde{\psi}$ by $\psi_0$ with impunity because the we get  an error not exceeding $CL_0\rho \mu^{-1} h^{-d} =C\rho |\log h| h^{1-d}$ which does not exceed $Ch^{1-d}$ as $\mu \le h^{-1}|\log h|^{-3}$ and is much less than the remainder estimate otherwise (actually we can replace in the above estimate $L_0$ by $\mu h$ but we do not need it). But then the second term in this approximation will be odd with respect to $\xi_1$ and result in $0$ in the final answer. 

So, all estimates of propositions \ref{prop-20-6-4} , \ref{prop-20-6-9}, and \ref{prop-20-6-11} related to the case ``$q=1$ and $\cA$ contains no unremovable cubic terms'' hold for an error term arising if we replace in the last term of (\ref{20-6-26}) $\sU$ by  $\hat{\sU}$ and also replace $\tilde{\psi}$ by $\psi_0$. Here and below $\hat{\sU}$ is our ``new'' first term, i.e. the first term in our new approximations. 

Therefore in comparison with our previous calculations we get an extra  term
\begin{multline}
h^{-1} \sum_n \int_{-\infty}^0 \phi_{0, L_0} ( \tau) \times \\
\Bigl(F_{t\to h^{-1}\tau}
\bigl(\bar{\chi}_{T_n}(t)-\bar{\chi}_{\bar{T}}(t)\bigr)
\Gamma \bigl((\hat{\sU}-\bar{\sU}) \psi_{0\,y} \,^t\!\cQ_{n\,y}\bigr)\Bigr)\,d\tau 
\label{20-6-60}
\end{multline}
where $\bar{\sU}$ is the ``old'' first term in the successive approximations. 

This is a correction term in question but we want to rewrite it in a more  explicit form. First of all, we can replace here $T_n$ by $+\infty$ (using the same arguments as before) and also replace $\bar{T}$ by $\bar{T}'=Ch|\log h|$ 
because $\bigl(F_{t\to h^{-1}\tau}\chi_T(t)
\Gamma (\sU' \psi_{0\,y} \,^t\!\cQ_{n\,y})\bigr)$
is negligible as $T\in [\bar{T}',\bar{T}]$, $|\tau|\le CL_0$ for both  
$\sU'=\bar{\sU}$ and $\sU'=\hat{\sU}$.

However if we consider expression
\begin{equation}
h^{-1} \sum_n \int_{-\infty}^0 \phi_{0, L_0} ( \tau) 
\Bigl(F_{t\to h^{-1}\tau} \bar{\chi}_{\bar{T}'}(t)\Gamma 
\bigl((\hat{\sU}-\bar{\sU}) \psi_{0\,y} \,^t\!\cQ_{n\,y} \bigr)\Bigr)\,d\tau 
\label{20-6-61}
\end{equation}
and apply the method of successive approximation to calculate it (considering $\bar{\cA}_0$ as an unperturbed operator and $\bar{\cA}_1$ as a perturbed one, we can get easily an error estimate  $C\rho |\log h|^2 h^{1-d}$ for it which is again less than the remainder estimate we are looking for. 

Thus modulo term not exceeding the remainder estimate we can rewrite expression (\ref{20-6-60}) as the same expression but with $\bar{\chi}_T$ replaced by $1$ and with $\bar{\chi}_{\bar{T}}$ replaced by $0$ i.e. as
\begin{equation}
h^{-1} \sum_n \int_{-\infty}^0 \phi_{0, L_0} ( \tau) 
\Bigl(F_{t\to h^{-1}\tau}
\Gamma (\hat{\sU}-\bar{\sU}) \psi_{0\,y} \,^t\!\cQ_{n\,y}\Bigr)\,d\tau 
\label{20-6-62}
\end{equation}
which in turn is exactly equal to
$h^{-d}\int \cN^\MW_\corr (x,0)\psi (x)\,dx$ with
\begin{multline}
\cN^\MW_\corr \Def 
(2\pi)^{-d+r}(\mu h)^r  \int \bigl(1-\phi_{0,L_0}(\tau)\bigr)\times \\
d_\tau \Bigl( \int \bigl(\mathbf{n}(x,\tau -z^2) -\mathbf{n}_0(x,\tau -z^2)\bigr) \cQ(z)\,dz \times f_1\cdots f_r \sqrt{g}\Bigr).
\label{20-6-63}
\end{multline}

Thus we arrive to 

\begin{proposition}\label{prop-20-6-13} 
Let $q=1$, $(h|\log h|)^{-\frac{1}{3}}\le \mu \le (h|\log h|)^{-1}$ and $\cA$ contain unremovable cubic terms. Let $\R^\MW_\corr$ be the difference between expression \textup{(\ref{20-6-7})} and 
\begin{equation}
h^{-d} \int \bigl(\cN^\MW(x,0)+\cN^\MW_\corr(x,0)\bigr)\,\psi (x)dx
\label{20-6-64}
\end{equation}
where $\cN^\MW_\corr$ is defined by \textup{(\ref{20-6-63})} with 
$L_0=C_0\mu h|\log h|$ with large enough constant $C_0$ and a proper symbol $\cQ$ equal $1$ in  $\{|z|\le C\bar{\rho}_1\}$ and  $\mathbf{n}_0$ and $\mathbf{n}$ are eigenvalue counting functions for $1\D$-operators $\mathbf{a}_0$ and $\mathbf{a}$ respectively. Then

\begin{enumerate}[label=(\roman*), fullwidth]
\item
In the general settings $\R^\MW_\corr$ does not exceed 
$C\mu h^{\frac{4}{3}-d}$;

\item In the framework of proposition \ref{prop-20-4-14} $\R^\MW_\corr$ does not exceed the right-hand expression of \textup{(\ref{20-4-88})};

\item In the framework of proposition \ref{prop-20-4-19}(iii) $\R^\MW_\corr$ does not exceed the right-hand expression of \textup{(\ref{20-4-95})};

\item  In the framework of proposition \ref{prop-20-4-21}(ii) $\R^\MW_\corr$ does not exceed the right-hand expression of \textup{(\ref{20-4-99})} as condition \textup{(\ref{20-4-85})} is violated;

\item  In the framework of proposition \ref{prop-20-4-21}(iii) $\R^\MW_\corr$ does not exceed the right-hand expression of \textup{(\ref{20-4-101})} as condition \textup{(\ref{20-4-85})} is violated.
\end{enumerate}
\end{proposition}

This proposition together with corresponding results of Section~\ref{sect-20-4} implies all statement of theorems dealing with $\R^\MW_\corr$ rather than $\R^\MW$ of theorems \ref{thm-20-6-19}(i) and \ref{thm-20-6-23}(ii) below.

It follows from proposition \ref{20-6-12}(iii) that in the framework of (\ref{20-4-93}) with $\nu(\hbar)=\hbar^2$ we do not need any correction. On the other hand, the following example shows the need in the correction term. 

\begin{example}\label{ex-20-6-14} 
Consider the case $q=1$, $r=2,3$ and $f_j=\const$ and commensurable and there are unremovable cubic terms in $\cA$  so that the $\cA_1=\cA_0+ \mu^{-1}\cE$ with $\cE$ commuting with $\cA$ (at each point $x$). Furthermore, let us assume that one of the following conditions is fulfilled:
\begin{phantomequation}\label{20-6-65}\end{phantomequation}
\begin{align}
&r=2,\quad f_1=2f_2,
\tag*{$\textup{(\ref*{20-6-65})}_{2}$}\label{20-6-65-2}\\
&r=3,\quad f_1=f_2+f_3, \quad f_2\ne f_3.
\tag*{$\textup{(\ref*{20-6-65})}_{3}$}\label{20-6-65-3}
\end{align}
Then 
\begin{align*}
&\cA_0=2Z_1^*Z_1+Z_2^*Z_2+D_1^2,\qquad
&&\cE= (\omega  Z_1^*Z_2^2)+(\omega  Z_1^*Z_2^2)^* ,\\[3pt]
&\cA_0=f_1Z_1^*Z_1+f_2Z_2^*Z_2+f_3Z_3^*Z_3+D_1^2,\; 
&&\cE= (\omega  Z_1^*Z_2Z_3)+(\omega  Z_1^*Z_2Z_3)^* 
\end{align*}
respectively with 
$\omega = \omega (x'',x''',\mu^{-1}hD'')$\,\footnote{\label{foot-20-27} One should not forget to add $C\mu^{-2}\cA_0^2$ to avoid some non-semi-boundedness related problems.}. 

We can already apply proposition \ref{prop-20-6-13}, getting for $\R^\MW_\corr$ better estimate than estimate $C\mu^{\frac{1}{2}}h^{1-d}$ for $\R^\MW$, but we can improve further the former one. 

Let us assume that $|\omega| \ge \epsilon$. Then (assuming large enough  smoothness of symbols) one can apply a transformation by operator $e^{i\hbar^{-1}\cL}$ with $\cL=\sum_j \beta_j Z_j^*Z_j$ with Hermitian 
$\beta_j = \beta_j (x'',x''',\mu^{-1}hD'')$ and transform operator (modulo  terms containing factors $\mu^{-1}hD_1$ or $\mu^{-2}$) to the same form but with Hermitian operator $\omega$. So, $\cE= \omega\cE_0$ with 
\begin{equation}
\cE_0=\left\{\begin{aligned}
&Z_1^*Z_2^2+Z_1Z_2^{*\,2}\qquad &&r=2\\
&Z_1^*Z_2Z_3+Z_1Z_2^*Z_3^*\qquad &&r=3\end{aligned}\right.
\label{20-6-66}
\end{equation}
commuting with $\cE_0$. Then one can decompose with respect to common eigenfunctions of $\cA_0$ and $\cE_0$ and to break the systems into separate equations and one can easily prove that
\begin{equation}
\R^\MW_\corr \le C\mu h^{\frac{3}{2}-d}
\label{20-6-67}
\end{equation}
without any microhyperbolicity or non-degeneracy condition; note that the right-hand expression  is much less than  
$C\mu^{\frac{1}{2}}h^{1-d}$ as $\mu h\ll 1$).

Let us select $V= -\sum_j (2 \bar{\alpha}_j+1)f_j\mu h$ and  
$\mu \ge h^{-\frac{1}{2}}$.

Further, let us replace $\mathbf{n}(x, -z^2,\hbar)$ by 
\begin{equation}
\#\bigl\{\alpha: \sum_j (2\alpha_j+1)f_j\mu h +V<0\bigr\} +
\mathbf{m}(x, -\xi_d^2,\hbar)\uptheta (\rho^2-z^2 )
\label{20-6-68}
\end{equation}
where $\mathbf{m}(x, \lambda,\hbar)$ is an eigenvalue counting function for operator $\cE$ restricted to the linear span of $\Upsilon_\alpha$ such that 
$\sum_j (2\alpha_j+1)f_j\mu h +V=0$. Then the error in $h^{-d}\cN^\MW_\corr$ will not exceed $C \mu \rho h^{1-d} \times (\mu \rho^2)^{-2}$. Taking $\rho=\bar{\rho}^*_1$ we get this error not exceeding 
$C\mu ^{-\frac{5}{2}}h^{-\frac{1}{2}-d}\ll \mu ^{\frac{1}{2}}h^{1-d}$ as long as
$\mu \gg h^{-\frac{1}{2}}$. 

But then 
\begin{multline}
\int \cN^\MW_\corr(x,0)\psi(x)\,dx =\\ 
(\mu h)^{r}\iint 
\mathbf{m}(x, -\xi_d^2,\hbar)\uptheta (\rho^2-\xi_d^2 )\psi(x)\,dx\,d\xi_1 +
O(C\mu ^{-\frac{5}{2}}h^{-\frac{1}{2}})
\label{20-6-69}
\end{multline}
which is exactly of magnitude 
$\mu^{-\frac{1}{2}}\times \mu h = \mu^{\frac{1}{2}}h^{1-d}$ and one cannot skip it without increasing the remainder estimate.
\end{example}

\section{Main theorems}
\label{sect-20-6-6}

\subsection{Microhyperbolicity assumption}
\label{sect-20-6-6-1}

We assume first that microhyperbolicity condition (see definition~\ref{def-20-1-2}) is fulfilled. Assume that 
\begin{equation}
\bar{\mu}_{(2)}=C(h|\log h|)^{-\frac{1}{2}}\le \mu 
\le  \epsilon (h|\log h|)^{-1}.
\label{20-6-70}
\end{equation}
The low bound is justified as $q\ge 2$ since otherwise theorem~\ref{thm-20-2-10} provides $O(h^{1-d})$ remainder estimate. For $q=1$ justification is that we will need to take $\varepsilon = C\mu^{-1}$ and it is smaller than 
$\varepsilon = C\mu h |\log h|$ used in the previous results if and only if 
$\mu \ge \bar{\mu}_{(2)}$. The only exception (when we can take 
$\varepsilon = C(\mu^{-1}h|\log h|)^{\frac{1}{2}}$) is the case $r=1$ but then 
$q\ge 2$ since we assume that $d\ge 4$.

We also assume that 
\begin{claim}\label{20-6-71} 
\underline{Either} $q\ge 2$ and $(l, \sigma) \succeq (1, 1)$ \underline{or} $q = 1$ and $(l, \sigma) \succeq (1, 2)$.
\end{claim}

\begin{theorem}\label{thm-20-6-15} 
Let conditions \textup{(\ref{20-1-1})}--\textup{(\ref{20-1-11})}, \textup{(\ref{20-6-61})}  and \textup{(\ref{20-6-70})} be fulfilled. Further, let microhyperbolicity condition  (see definition~\ref{def-20-1-2}) be fulfilled. Then there exist two framing approximations (see footnote \footref{book_new-foot-18-16} of Chapter~\ref{book_new-sect-18} of \cite{futurebook}) such that
\begin{equation}
\R^\MW \le Ch^{1-d} + 
Ch^{-d}(\mu h|\log h|)^{\frac{q}{2}} \mu ^{-l}|\log h|^{-\sigma}.
\label{20-6-72}
\end{equation}
\end{theorem}

\begin{corollary}\label{cor-20-6-16} 
Under microhyperbolicity condition for $\mu \le \epsilon (h|\log h|)^{-1}$ the sharp remainder estimate $\R^\MW\le Ch^{1-d}$  holds as \underline{either} 
$q\ge 2$ and $(l,\sigma)\succeq (1,1)$  \underline{or} $q=1$ and
$(l,\sigma)\succeq (\frac{3}{2}, 1)$.
\end{corollary}

\subsection{General case}
\label{sect-20-6-6-2}

Now let us consider intermediate magnetic field without microhyperbolicity assumption. 

\begin{theorem}\label{thm-20-6-17} 
Let conditions \textup{(\ref{20-1-1})}--\textup{(\ref{20-1-11})}, \textup{(\ref{20-6-71})} and \textup{(\ref{20-6-70})} be fulfilled.
Then there exist two framing approximations  such that

\begin{enumerate}[label=(\roman*), fullwidth]
\item
For $q\ge 3$, $(l,\sigma)=(1,1)$, $(\bar{l},\bar{\sigma})=(2,1)$ estimate
$\R^\MW \le Ch^{1-d}$ holds;

\item
For $q=2$, $(l,\sigma)=(1,1)$, $(\bar{l},\bar{\sigma})=(2,1)$ estimate
\begin{equation}
\R^\MW \le Ch^{1-d} + C\mu h^{\frac{5}{3} -d} 
\label{20-6-73}
\end{equation}
holds; in particular, $\R^\MW \le Ch^{1-d}$ as $\mu \le h^{-\frac{2}{3}}$;

\item For $q=1$, $(l,\sigma)=(1,2)$ estimate
\begin{equation}
\R^\MW  \le Ch^{1-d} + C\mu h^{ \frac{4}{3} -d} +  C\mu^{\frac{1}{2}} h^{1-d}
\label{20-6-74}
\end{equation}
holds.
\end{enumerate}
\end{theorem}

\begin{remark}\label{rem-20-6-18}
\begin{enumerate}[label=(\roman*), fullwidth]
\item 
Note that  estimates (\ref{20-6-73}), (\ref{20-6-74}) are better than (\ref{20-2-32}) with $q=2,1$ respectively iff  
$\mu \ge (h |\log h|)^{- \frac{1} {2}}$.

\item
The second terms in (\ref{20-6-73}), (\ref{20-6-74}) are due to the remainder estimates  for $q$-dimensional Schr\"odinger operator (it is not important that operator in question is also $r$-dimensional $\mu^{-1}h$-pseudo-differential operator). These estimates correspond to the $\sC^{1,1}$-smoothness and cannot be improved in the general case even if we assume much the larger smoothness because we have a matrix, not a scalar $h$-pseudo-differential operator as for $d=3$. However, in certain cases we can do better than this; we need to consider $q=1,2$ only.

\item
In (\ref{20-6-74}) the third term estimates an approximation error arising when we skip irreducible $O(\mu^{-1})$ terms due; we can prevent these terms under certain conditions. 
\end{enumerate}
\end{remark}

\begin{theorem}\label{thm-20-6-19} 
Let $q=1$, $(1,2)\preceq (l,\sigma)\preceq (2,0)$ and conditions \textup{(\ref{20-1-1})}--\textup{(\ref{20-1-11})}, and \textup{(\ref{20-6-70})} be fulfilled. Then there exist two framing approximations  such that

\begin{enumerate}[label=(\roman*), fullwidth]
\item 
Estimate
\begin{multline}
\R^\MW_\corr \Def \\
|\int \Bigl(e(x,x,\tau)-h^{-d}\cN^\MW (x,\tau )- h^{-d}\cN^\MW_\corr (x,\tau ) \Bigr)\psi(x)\,dx|\\
\le Ch^{1-d} + C\mu h^{\frac{4}{3} -d} + 
C\mu^{1-\frac{l}{2}} h^{1-d}|\log h|^{-\frac{\sigma}{2}}
\label{20-6-75}
\end{multline}
holds where 
\begin{equation}
| \cN^\MW_\corr (x,\tau )|\le C\mu^{\frac{1}{2}}h 
\label{20-6-76}
\end{equation}
is defined by formula \textup{(\ref{20-6-69})} and in the general case one cannot improve estimates \textup{(\ref{20-6-74})}, \textup{(\ref{20-6-76})};

\item
If  one of the following conditions
\begin{claim}\label{20-6-77}
There are no third-order resonances;
\end{claim}
\vskip-20pt
\begin{claim}\label{20-6-78}
$g^{jk}=\const$, $F_{jk}=\const$
\end{claim}
is fulfilled then estimate \textup{(\ref{20-6-75})} holds for $\R^\MW$.
\end{enumerate}
\end{theorem}

As $r\ge 2$  and condition (\ref{20-6-70}) is fulfilled let us consider 
\begin{equation}
\mathbf{n}_0(x, \hbar,\tau)\Def
\#\bigl\{\alpha \in \bZ^{+\, r},  
\sum_j (2\alpha _j+1)f_j(x) \hbar +V(x) < \tau\bigr\}.
\label{20-6-79}
\end{equation}
Let us assume that estimate (\ref{20-4-93}) holds with  $\nu (\hbar)=o(\hbar)$. 

\begin{remark}\label{rem-20-6-20}
Estimate (\ref{20-4-93})  holds with $\nu(\hbar)=\hbar$ for sure and it holds  with  $\nu(\hbar)=\hbar^\kappa$, $\kappa >1$  if \underline{either} $f_j$ have constant multiplicities and 
\begin{equation}
\rank\{\nabla \frac {f_2}{f_1} , \ldots, \nabla \frac {f_r}{f_1} \}\ge \kappa-1   \qquad \forall x
\label{20-6-80}
\end{equation}
\underline{or}  $\frac {f_j}{f_1}=\const$ and the system 
$\bigl\{  \frac {f_2}{f_1},\ldots,  \frac {f_r}{f_1} \bigr\}$ has some Diophantine properties\footnote{\label{foot-20-28} See examples \ref{ex-20-4-17}, \ref{ex-20-4-18} for details.}.
\end{remark}

Under condition (\ref{20-4-93}) we can improve theorems \ref{thm-20-6-17},  \ref{thm-20-6-19} as $q=1,2$: 

\begin{theorem}\label{thm-20-6-21}
Let $q=1,2$ and conditions \textup{(\ref{20-1-1})}--\textup{(\ref{20-1-11})}, \textup{(\ref{20-6-70})}, \textup{(\ref{20-4-93})}   be fulfilled. Then there exist two framing approximations  such that
\begin{enumerate}[label=(\roman*), fullwidth]
\item 
For $q=2$, $(l,\sigma)=(1,1)$, $(\bar{l},\bar{\sigma})=(2,1)$ estimate
\begin{equation}
\R^\MW \le Ch^{1-d} + C\nu (\mu h) h^{\frac{2}{3} -d} 
\label{20-6-81}
\end{equation}
holds;
\item
For $q=1$, $(1,2)\preceq (l,\sigma)\preceq (2,0)$ estimates
\begin{align}
&\R^\MW \le Ch^{1-d}+ 
C\mu ^{- \frac{1}{2} }h^{ \frac{1}{2}-d}|\log h|^{\frac{1}{2}}+ 
C\nu(\mu h)  \bigl( h^{ \frac{1}{3}-d}+ \mu ^{-\frac{1}{2}}h^{-d}\bigr),
\label{20-6-82}\\[3pt]
&\R^\MW_\corr \le Ch^{1-d}+ 
C\mu ^{ \frac{1}{2}-l} h^{ \frac{1}{2}-d}|\log h|^{ \frac{1}{2}- \sigma }+ \label{20-6-83}\\[2pt]
&\quad C\bigl(\nu(\mu h) +\mu ^{-1}\bigr) \times 
\bigl( h^{\frac{1}{3}-d}+ 
\mu ^{-\frac{l}{2}} h^{-d}|\log h|^{-\frac{\sigma}{2}}\bigr)h^{-d}\notag
\end{align} 
hold;

\item
For $q=1$, $(1,2)\preceq (l,\sigma)\preceq (2,0)$ as one of conditions \textup{(\ref{20-6-77})}, \textup{(\ref{20-6-78})} is fulfilled the following estimate holds:
\begin{align}
\R^\MW \le &Ch^{1-d}+ 
C\mu ^{\frac{3}{2}-l}h^{\frac{3}{2}-d}|\log h|^{\frac{1}{2}- \sigma }+ \label{20-6-84}\\
&C\nu(\mu h) \bigl( h^{\frac{1}{3}-d}+
\mu ^{-\frac{l}{2}}h^{-d}|\log h|^{-\frac{\sigma}{2}}\bigr).
\notag
\end{align} 
\end{enumerate}
\end{theorem}

Recall that proposition \ref{prop-20-6-13} estimates $\R^\MW_\corr$ defined as a difference between expressions (\ref{20-6-7}) and (\ref{20-6-64}). One can prove easily that

\begin{theorem}\label{thm-20-6-22}
In the framework of statements (i)--(v) of proposition \ref{prop-20-6-13} the corresponding estimates hold also for $\R^\MW_\corr$ defined by \textup{(\ref{20-6-75})}, and \textup{(\ref{20-6-63})}.
\end{theorem}

\subsection{Constant multiplicities of $f_j$}
\label{sect-20-6-6-3}

Assume now that 
\begin{claim}\label{20-6-85}
Matrix $(F^j_k)$ has eigenvalues $\pm i f_m$ of constant multiplicities.
\end{claim}
Then the microhyperbolicity assumption is equivalent to
\begin{equation}
|V+\sum_{1\le j\le r}f_j\tau_j| + |\nabla (V+\sum_{1\le j\le r}f_j\tau_j )|\ge \epsilon\qquad \forall \tau_1,\dots,\tau_r\in\bR^+.
\label{20-6-86}
\end{equation}
Further, if 
\begin{equation}
f_2/f_1=\const , \ldots, f_r/f_1=\const
\label{20-6-87}
\end{equation}
then (\ref{20-6-86}) is equivalent to
\begin{equation}
|\nabla \frac {V} {f_1} |\ge \epsilon.
\label{20-6-88}
\end{equation}

\begin{theorem}\label{thm-20-6-23} 
Let conditions \textup{(\ref{20-1-1})}--\textup{(\ref{20-1-11})}, \textup{(\ref{20-6-71})}, and \textup{(\ref{20-6-85})} be fulfilled. Then there exist two framing approximations  such that

\begin{enumerate}[label=(\roman*), fullwidth]
\item For $q=2$, $(1,1)\preceq (l,\sigma)\preceq (2,0)$, $(\bar{l},\bar{\sigma})=(2,1)$ estimate
\begin{equation}
\R^\MW \le Ch^{1-d} +
C \mu h^{1-d+  \frac {ql} {l+2} }|\log \mu|^{- \frac {q\sigma} {(l+2)} } 
\label{20-6-89}
\end{equation}
holds. 

In particular, $R^\MW\le Ch^{1-d}$ as $(l,\sigma)=(2,0)$.

\item As  $q=1$, $(1,2)\preceq (l,\sigma)\preceq (2,0)$, $(\bar{l},\bar{\sigma})=(2,1)$
estimate
\begin{multline}
\R^\MW_\corr\le
Ch^{1-d} + C\mu h^{2-d-\frac{2} {l+2} }|\log h|^{- \frac {\sigma} {(l+2)} } + \\[3pt]
Ch^{1-d}(\mu h|\log h|)^{\frac{1}{2}} \mu ^{1-l}|\log h|^{-\sigma} +
C\mu^{1- \frac {l}{2}} h^{1-d} |\log h|^{- \frac {\sigma}{2} } +\\[3pt] 
C\mu ^{\frac {2l+1}{3l} }h^{ \frac {4}{3}-d} |\log h|^{- \frac {\sigma} {3l} }
\label{20-6-90}
\end{multline}
holds. 

\item
As $q=1$, $(1,2)\preceq (l,\sigma)\preceq (2,0)$, $(\bar{l},\bar{\sigma})=(2,1)$
and one of  conditions \textup{(\ref{20-6-77})}, \textup{(\ref{20-6-78})} is fulfilled estimate
\begin{multline}
\R^\MW\le
Ch^{1-d} + C\mu h^{2-d- \frac 2 {l+2} } |\log h|^{- \frac {\sigma} {(l+2)} } + \\
Ch^{1-d}(\mu h|\log h|)^{\frac{1}{2}} \mu ^{1-l} |\log h|^{-\sigma} +
C\mu^{1- \frac {l} {2} } h^{1-d}|\log h|^{- \frac {\sigma}{2}}
\label{20-6-91}
\end{multline}
holds.
\end{enumerate}
\end{theorem}

\begin{remark}\label{rem-20-6-24}

\begin{enumerate}[label=(\roman*), fullwidth]
\item
Note that the right-hand expression of (\ref{20-6-90}) contains an extra term (in comparison with (\ref{20-6-91}).

\item 
Also note  that for $q=1$ estimate (\ref{20-6-89}) is better than (\ref{20-2-32}) as
\begin{equation}
\mu^{l+1}h |\log h|^{\sigma+1}\ge 1;
\label{20-6-92}
\end{equation}
in particular, as $(l,\sigma)=(2,0)$ we get estimate 
$\R^\MW \le Ch^{1-d}+C\mu h^{ \frac{3}{2}-d}$ which always is better than (\ref{20-2-32}) as $\mu\ge \bar{\mu}_{1(1)}$.
\end{enumerate}
\end{remark}

Under condition (\ref{20-4-93}) one can improve theorem \ref{thm-20-6-23}:

\begin{theorem}\label{thm-20-6-25} 
Let conditions \textup{(\ref{20-1-1})}--\textup{(\ref{20-1-11})}, \textup{(\ref{20-6-71})} and  \textup{(\ref{20-6-85})}, and \textup{(\ref{20-4-93})} be fulfilled. Then there exist two framing approximations  such that

\begin{enumerate}[label=(\roman*), fullwidth]
\item As $q=2$, $(1,1)\preceq (l,\sigma)\preceq (2,0)$, $(\bar{l},\bar{\sigma})=(2,1)$
estimate
\begin{equation}
\R^\MW\le Ch^{1-d} + 
C\nu(\mu h) h^{1-d+ \frac {l-2} {l+2} } |\log h|^{- \frac {2\sigma} {l+2} } 
\label{20-6-93}
\end{equation}
holds. 

\item
As $q=1$, $(1,2)\preceq (l,\sigma)\preceq (2,0)$, $(\bar{l},\bar{\sigma})=(2,1)$
and  one of conditions \textup{(\ref{20-6-77})}, \textup{(\ref{20-6-78})} is fulfilled estimate
\begin{multline}
\R^\MW\le
Ch^{1-d} + Ch^{-d}(\mu h|\log h|)^{\frac{1}{2}} \mu ^{-l}|\log h|^{-\sigma}
\\[3pt]
C\nu(\mu h)\Bigl(h^{1-d- \frac{2} {l+2} }|\log h|^{-\frac {\sigma} {(l+2)} } + 
C\mu^{- \frac {l} {2} } h^{-d} |\log h|^{- \frac {\sigma} {2} }\Bigr)
\label{20-6-94}
\end{multline}
holds.

\item
If $f_1=\const$, \ldots, $f_r=\const$, $(l,\sigma)\succeq (2,0)$ and nondegeneracy assumption \textup{(\ref{20-2-42})} is fulfilled then $\R^\MW =O(h^{1-d})$.
\end{enumerate}
\end{theorem}

\chapter{Stronger magnetic field: calculations and main results}
\label{sect-20-7}

In this section we consider the cases of the strong magnetic field as
$\epsilon (h|\log h|)^{-1}\le \mu \le \epsilon h^{-1}$, of the very strong magnetic field as $\epsilon h^{-1}\le \mu \le C_0 h^{-1}$ and of superstrong magnetic field as $\mu \ge C_0 h^{-1}$.

\section{Strong magnetic field}
\label{sect-20-7-1}

\subsection{Calculations}
\label{sect-20-7-1-1}

Consider now case of a strong magnetic field $\epsilon (h|\log h|)^{-1}\le \mu \le \epsilon h^{-1}$.  The same arguments work again, with no outer zone and necessity to join outer and intermediate zones. Because of this one should deal directly with (\ref{20-6-7}) rather than (\ref{20-6-26}). 

Therefore  (\ref{20-6-7}) is preserved modulo $O(h^s)$ as $T_n$ are replaced by some lesser values which depend on the assumptions  described in the previous subsections; then we apply the method of the successive approximations in the same way as above replacing $\sU$ by $\bar{\sU}$ and $\tilde{\psi}$ by $\psi_0$ and estimating an error in (\ref{20-6-7}) cause by this; finally we replace 
$T'_n$ by $T_n$ again with $O(h^s)$ error and then $T_n$ by $+\infty$. We leave easy details to the reader.

Then one can easily prove several propositions below; we leave easy and standard details to the reader. The first proposition holds for both strong and very strong magnetic field:

\begin{proposition}\label{prop-20-7-1} 
Let  $\epsilon (h|\log h|)^{-1}\le \mu \le C_0 h^{-1}$, $(l,\sigma)=(1,1)$ as $q\ge 2$, $(l,\sigma)=(1,2)$ as $q=1$. Then one can rewrite \textup{(\ref{20-6-7})} as $h^{-d}\int \cN^\MW(x,0)\psi (x)\,dx$ with an error not exceeding $Ch^{1-d}$ as  $q\ge 2$, $C\mu^{\frac{1}{2}}h^{1-d}$ as $q=1$\,\footnote{\label{￼foot-20-29} One can improve it as there are no unremovable cubic terms but there is no point as the general remainder estimate 
$O(\mu h^{\frac{4}{3}-d})$ is larger.}.
\end{proposition}

\begin{proposition}\label{prop-20-7-2} 
Let  $\epsilon h^{-1}|\log h|^{-1}\le \mu \le \epsilon h^{-1}$, $(l,\sigma)\succeq(1,1)$ as $q\ge 2$, $(l,\sigma)\succeq(1,2)$ as $q=1$. Then under microhyperbolicity  condition (see definition~\ref{def-20-1-2}) one can rewrite \textup{(\ref{20-6-7})} as  $h^{-d}\int \cN^\MW(x,0)\psi (x)\,dx$ with an error not exceeding $Ch^{1-d}$.
\end{proposition}

\begin{proposition}\label{prop-20-7-3} 
Statements of propositions \ref{prop-20-6-9}--\ref{prop-20-6-13}, theorem \ref{thm-20-6-22}, problem~\ref{Problem-20-6-10} and example \ref{ex-20-6-14} remain true as $\epsilon (h|\log h|)^{-1}\le \mu \le \epsilon h^{-1}$. 
\end{proposition}

\subsection{Main results}
\label{sect-20-7-1-2}

Recall that in the case of the strong magnetic field we pick up
\begin{equation}
\varepsilon \asymp C(\mu^{-1}h|\log h|)^{\frac{1}{2}} + Ch|\log h|\asymp 
h|\log h|
\label{20-7-1}
\end{equation}

The first of our theorems hold for both strong and very strong magnetic field:

\begin{theorem}\label{thm-20-7-4} 
Let $\epsilon (h|\log h|)^{-1}\le \mu \le C_0 h^{-1}$. Let  conditions \textup{(\ref{20-1-1})}--\textup{(\ref{20-1-11})}, and \textup{(\ref{20-6-71})}  be fulfilled.  Then there exist two framing approximations (see footnote \footref{book_new-foot-18-16} of Chapter~\ref{book_new-sect-18} of \cite{futurebook}) such that

\begin{enumerate}[label=(\roman*), fullwidth]
\item For $q\ge 3$, $(l,\sigma)=(1,1)$, $(\bar{l},\bar{\sigma})=(2,1)$ sharp remainder estimate $\R^\MW\le Ch^{1-d}$ holds.

\item For $q=2$, $(l,\sigma)=(1,1)$, $(\bar{l},\bar{\sigma})=(2,1)$ estimate
\textup{(\ref{20-6-73})} holds: $\R^\MW\le C\mu h^{ \frac{5}{3}-d}$.

\item For $q=1$, $(1,2)\preceq (l,\sigma)\preceq (\bar{l},\bar{\sigma})=(2,1)$ estimate
\textup{(\ref{20-6-74})} holds: $\R^\MW\le C\mu h^{\frac{4}{3}-d}$.
\end{enumerate}
\end{theorem}

In all other cases we need to restrict ourselves here to the strong magnetic field case. Here conclusion is very simple: all results of Subsection~\ref{sect-20-6-6} remain valid:

\begin{theorem}\label{thm-20-7-5}
Let $\epsilon (h|\log h|)^{-1}\le \mu \le \epsilon h^{-1}$. Let  conditions \textup{(\ref{20-1-1})}--\textup{(\ref{20-1-11})}, and \textup{(\ref{20-6-71})}  be fulfilled. 

Further, let \underline{either} microhyperbolicity assumption (see definition~\ref{def-20-1-2}) \underline{or} $(l,\sigma)=(2,0)$,
$f_1=\const,\ldots, f_r=\const$ and nondegeneracy assumption \textup{(\ref{20-2-42})} be fulfilled.

Then there exist two framing approximations such that $\R^\MW \le Ch^{1-d}$.
\end{theorem}

\begin{theorem}\label{thm-20-7-6}
In the case of the strong magnetic field $\epsilon (h|\log h|)^{-1}\le \mu \le \epsilon h^{-1}$ all  statements of theorems     \ref{thm-20-6-23}  and \ref{thm-20-6-25} remain true.
\end{theorem}

We leave all  easy details to the reader.

\section{Very strong magnetic field}
\label{sect-20-7-2}

\subsection{Calculations}
\label{sect-20-7-2-1}

Consider now case of the very strong magnetic field 
$\epsilon h^{-1}\le \mu \le C_0h^{-1}$. The same arguments work again, with no intermediate zone anymore. We do not need consider the general case as it is covered by theorem~\ref{thm-20-7-4}. Thus we should consider only $q=1,2$.

On the other hand microhyperbolicity and non-degeneracy assumptions should be modified to definition~\ref{def-20-5-9} and (\ref{20-5-10}) (again non-degeneracy assumption is considered as $f_1,\ldots,f_r/f_1$ are constant). 

Finally, number theoretical arguments are no more applicable at all.

Again one can prove easily propositions below; we leave all easy details to the reader.

\begin{proposition}\label{prop-20-7-7} 
Let  $\epsilon h^{-1}\le \mu \le C_0 h^{-1}$, and $(l,\sigma)=(1,1)$ as $q\ge 2$, $(l,\sigma)=(1,2)$ as $q=1$. Then one can rewrite expression \textup{(\ref{20-6-7})} as 
$h^{-d}\int \cN^\MW (x,0)\psi (x)\,dx$ with an error not exceeding $Ch^{1-d}$ as 
$q\ge 2$, $Ch^{\frac{1}{2}-d}$ as $q=1$.
\end{proposition}

\begin{proposition}\label{prop-20-7-8} 
Let  $\epsilon h^{-1}\le \mu \le C_0 h^{-1}$,  $(l,\sigma)=(1,1)$ as $q\ge 2$, $(l,\sigma)=(1,2)$ as $q=1$. Then under microhyperbolicity condition (see definition~\ref{def-20-5-9}) with $\ell$ independent on $\mathbf{\tau}$ one can rewrite \textup{(\ref{20-6-7})} as $h^{-d}\int \cN^\MW(x,0)\psi (x)\,dx$ with an error not exceeding $Ch^{1-d}$.
\end{proposition}

\begin{proposition}\label{prop-20-7-9} 
Statements of propositions \ref{prop-20-6-9} and problem~\ref{Problem-20-6-10}
remain true as $\epsilon h^{-1}\le \mu \le C_0 h^{-1}$. 
\end{proposition}

\subsection{Main theorems}
\label{sect-20-7-2-2}

We arrive then to the following theorems:

\begin{theorem}\label{thm-20-7-10}
Let $\epsilon h^{-1}\le \mu \le C_0 h^{-1}$. Let  conditions \textup{(\ref{20-1-1})}--\textup{(\ref{20-1-11})}, and \textup{(\ref{20-6-71})}  be fulfilled. 

Further, let \underline{either} microhyperbolicity assumption (see definition~\ref{def-20-5-9}) \underline{or} $(l,\sigma)=(2,0)$,
$f_1=\const,\ldots, f_r=\const$ and nondegeneracy assumption \textup{(\ref{20-5-10})} be fulfilled.

Then there exist two framing approximations such that $\R^\MW \le Ch^{1-d}$.
\end{theorem}

\begin{theorem}\label{thm-20-7-11}
In the case of the strong magnetic field $\epsilon (h|\log h|)^{-1}\le \mu \le \epsilon h^{-1}$ all  statements of theorem    \ref{thm-20-6-23} remain true.
\end{theorem}

We leave all  easy details to the reader.

\section{Superstrong Magnetic Field}
\label{sect-20-7-3}

Consider now case of superstrong magnetic field $\mu \ge C_0h^{-1}$ and a magnetic-Schr\"odinger-Pauli operator. 

Then we need to modify  condition $\textup{(\ref{20-1-9})}_3$, replacing \ref{20-1-9-1} by
\begin{equation}
g^{jk}, F_{jk}\in \sC^{\bar{l},\bar{\sigma}}, \qquad 
W=V+\mu h \sum_jf_j \in \sC^{l,\sigma}
\label{20-7-2}
\end{equation}
and assuming that 
\begin{equation}
W+2f_j\mu h\ge \epsilon\qquad \forall j.
\label{20-7-3}
\end{equation}

Under these assumptions we have really a scalar operator, and should take

\begin{equation}
\varepsilon = Ch|\log h|
\label{20-7-4}
\end{equation}
and microhyperbolicity condition (see definition~\ref{def-20-5-9}) transforms into 
\begin{equation}
|W|+|\nabla W |\ge \epsilon
\label{20-7-5}
\end{equation}
while non-degeneracy assumption (\ref{20-5-10}) transforms into 
\begin{equation}
|W|+|\nabla W |\le \epsilon \implies |\det \Hess W|\ge \epsilon.
\label{20-7-6}
\end{equation}

Skipping calculations (which are not much different from what we did before, so we leave the easy details to the reader) we arrive immediately to

\begin{theorem}\label{thm-20-7-12} 
Let $\mu \ge C_0h^{-1}$ and let conditions \textup{(\ref{20-1-1})}--\textup{(\ref{20-1-8})},  \ref{20-1-9-2},  \textup{(\ref{20-6-71})},  \textup{(\ref{20-7-2})}, and \textup{(\ref{20-7-3})} be fulfilled. Then

\begin{enumerate}[label=(\roman*), fullwidth]
\item For $q\ge 3$ remainder estimate $\R^\MW \le C\mu^rh^{r+1-d}$    holds;

\item  For $q=1,2$ and  microhyperbolicity condition \textup{(\ref{20-5-9})}
fulfilled,  the remainder estimate $\R^\MW \le C\mu^r h^{r+1-d}$ holds;

\item For $q=2$ and  $(l,\sigma)\succeq (2,1)$  the remainder estimate 
$\R^\MW \le C\mu^r h^{r+1-d}$ holds;

\item For $q=1$, $(l,\sigma)\succeq (2,1)$,  and non-degeneracy condition \textup{(\ref{20-5-10})} fulfilled the remainder estimate 
$\R^\MW \le C\mu^r h^{r+1-d}$ holds;

\item As \underline{either} $q=2$, $(l,\sigma)\prec (2,0)$ \underline{or} $q=1$  then the remainder estimate
\begin{equation}
\R^\MW \le C\mu^r h^{r-d}\left\{\begin{aligned}
& h^{  \frac {2l} {l+2}  } 
|\log h|^{- \frac {2\sigma} {l+2} }\qquad &&\text{as\ \ }q=2 ,\\
&h^{ \frac {l} {l+2} } |\log h|^{- \frac {\sigma}{l+2}}\qquad &&\text{as\ \ }q=1 
\end{aligned}\right.
\label{20-7-7}
\end{equation}
holds.
\end{enumerate}
\end{theorem}

\begin{remark}\label{rem-20-7-13}
Under assumption (\ref{20-7-3}) he principal part of asymptotics is 
\begin{multline}
h^{-d}\int \cN_d^\MW (x,\tau)\psi(x)\, dx=\\ 
\omega_{d-2r}(2\pi )^{-d+r} 
\mu ^rh^{-d+r}
\int \Bigl(\tau - \sum_j f_j\mu h -V\Bigr)_+^{\frac{q}{2}}
f_1\cdots f_r\sqrt g\psi(x)\,dx.
\label{20-7-8}
\end{multline}
\end{remark}

\chapter{Degenerating magnetic field}
\label{sect-20-8}

\begin{remark}\label{rem-20-8-1}
One can get rid off assumption $|V|\ge \epsilon$ as $\mu \le \epsilon h^{-1}$ easily (see Subsubsection~\ref{book_new-sect-18-9-5-1}.1 of \cite{futurebook}).
\end{remark}

Condition that non-zero eigenvalues of $F^j_k$ are uniformly disjoint from $0$ is more subtle. Basically we want to apply arguments of Subsubsection~\ref{book_new-sect-18-9-5-2}.2 of \cite{futurebook} but we should take into account that different eigenvalues have different magnitudes.

Usually we assume that $q=1$ (so at generic point only $1$ eigenvalue is $0$), but at this moment we just assume
\begin{claim}\label{20-8-1}
At each point matrix $F^j_k(x)$ has at least $q$ eigenvalues equal to $0$ and let $\pm i f_j(x)$ denote other eigenvalues (which may be also vanish at some points\footnote{\label{foot-20-30} We will need to assume that these eigenvalues do not vanish in the generic points.}
\end{claim}
Let us introduce
\begin{equation}
\gamma (x) = \epsilon_0 \min _j f_j(s) 
\label{20-8-2}
\end{equation}
where in partitions but not in statements we replace $\gamma$ by $\max(\gamma,\bar{\gamma})$ and $\bar{\gamma}$ will be chosen later.

\section{Weak magnetic field}
\label{sect-20-8-1}

Consider some point $y$, its $\gamma(y)$ vicinity and classify $f_1,\ldots, f_r$ into two groups: those which are greater than $4\nu \gamma$ (we denote them by $f_{r'+1},\ldots, f_{r}$) and those which are less than $\nu \gamma$ (we denote them by $f_1,\ldots, f_{r'}$); an appropriate constant  $\nu=\nu(y)\in [c_0,c_1]$ exists with arbitrarily large $c_0\ge 1$ and $c_1=c_1(c_0)$. Then the similar inequalities 
\begin{equation}
f_j(x) \le 2\nu \gamma (y)\quad \forall j=1,\ldots,r',\qquad
f_j\ge 3\nu \gamma(y) \quad \forall j=r'+1,\ldots,r
\label{20-8-3}
\end{equation}
hold in $\gamma$-vicinity of $y$; here $r'=r'(y)$. 

We claim that
\begin{claim}\label{20-8-4}
As $\gamma (y)\ge C_0\mu^{-1}$ the contribution of $B(y,\gamma(y))$ to the remainder  does not exceed
\begin{equation}
C\bigl(\gamma^{-1}+ \mu (\mu h|\log h|)^{\frac{q}{2}}\bigr) h^{1-d}\gamma^d.
\label{20-8-5}
\end{equation}
\end{claim}
Really, one can use the weak magnetic field approach and take 
$T_*=\epsilon \mu^{-1}$. Further, one can take 
$T^* = T^*(y, |\xi'''|)=\epsilon \gamma(y)|\xi'''|$ and after summation with respect to $\xi'''$ we get (\ref{20-8-5}) as $q\ge 2$; as $q=1$ we get the same expression albeit with an extra logarithmic factor in the first term and we get rid off it exactly in the same manner as in Subsubsection~\ref{sect-20-2-2-2}.2. 

Here  $C\mu (\mu h|\log h|)^{\frac{q}{2}} h^{1-d}\gamma^d$ estimates a contribution of the complementary  zone  
$\Omega_\out^c =\{|\xi'''|\le \bar{\rho}_1\}$ with
$\bar{\rho}_1=C\max \bigl((\mu h|\log h|)^{\frac{1}{2}}, \mu^{-1}\gamma^{-1}\bigr)$ as 
$\bar{\rho}_1=C(\mu h|\log h|)^{\frac{1}{2}} \ge C\mu^{-1}\gamma^{-1}$; otherwise contribution of the complementary zone does not exceed 
$C\bar{\rho}_1\mu h^{1-d}\gamma^d\asymp C\gamma^{-1}h^{1-d}$ which  is  the first term in (\ref{20-8-5}).

This estimate (\ref{20-8-5}) could be improved but it does not make any difference. Summation  with respect to all balls with $\gamma \gtrsim \bar{\gamma}$ results in 
\begin{equation}
Ch^{1-d} \int_{\{\gamma (x)\ge \bar{\gamma}\}} \gamma^{-1}\,dx + 
C\mu (\mu h|\log h|)^{\frac{q}{2}} h^{1-d};
\label{20-8-6}
\end{equation}
while contribution of all balls with $\gamma \lesssim \bar{\gamma}$ does not exceed
\begin{equation}
C\mu h^{1-d} \int_{\{\gamma (x)\le C_0\bar{\gamma}\}}\, dx;
\label{20-8-7}
\end{equation}
Thus the total remainder $\R^\T$ does not exceed 
\begin{multline}
Ch^{1-d}+ C\mu (\mu h|\log h|)^{\frac{q}{2}}h^{-d}h^{1-d}+\\[3pt]
Ch^{1-d} \Bigl(\int_{\{\gamma (x)\ge \bar{\gamma}\}} \gamma^{-1}\,dx +
C\mu  \int_{\{\gamma (x)\le C_0\bar{\gamma}\}}\,dx\Bigr)
\label{20-8-8}
\end{multline}
where we included the first term for ``compatibility'' only (it does not exceed the second line anyway) and the optimal results are achieved as $\bar{\gamma}\asymp \mu^{-1}$.

We leave to the reader to prove by our standard arguments that the same estimate holds for $\R^\MW$ and $\R^\W_{(\infty)}$. Then we arrive to the following

\begin{theorem}\label{thm-20-8-12}
\begin{enumerate}[label=(\roman*), fullwidth]
\item 
Theorem \ref{thm-20-2-10}  remains true without assumption \textup{(\ref{20-1-4})}, replaced by \textup{(\ref{20-8-1})}, but with the right-hand expression \textup{(\ref{20-8-8})};

\item
Under assumption 
\begin{gather}
\int_{\mu^{-1}}^1  \gamma^{-1} d\upmu (\gamma) <\infty
\label{20-8-9}\\
\shortintertext{with}
\upmu (\gamma)\Def \mes \bigl( \{x:\, \gamma (x) <\gamma \}\bigr)
\label{20-8-10}
\end{gather}
an extra term (the second line in \textup{(\ref{20-8-8})}) is $O(h^{1-d})$ and one can skip it.
\end{enumerate}
\end{theorem}

\begin{Problem}\label{Problem-20-8-3}
In our usual manner get rid off logarithmic factor in the second term of (\ref{20-8-8}).
\end{Problem}

\begin{Problem}\label{Problem-20-8-4}
Prove that in the generic setting $\upmu (\gamma)= O(\gamma^3)$.
\end{Problem}

\section{Intermediate magnetic field}
\label{sect-20-8-2}

\subsection{General case}
\label{sect-20-8-2-1}

Situation is rather simple if $r=1$ (or all eigenvalues  $f_j$ have the same magnitude $\asymp \gamma$): then we can use the same simple rescaling technique as in Subsubsection~\ref{book_new-sect-18-9-5-2}.2 of \cite{futurebook}. However in the general case there are some eigenvalues $\asymp \gamma$ and there are some much larger eigenvalues, may be of magnitude $1$. 

Notice that after rescaling the eigenvalue of magnitude $\gamma$ would have magnitude $\mu \gamma^2$ and to deal with this properly we need $\mu \gamma^2 \gg 1$ (or at least $\mu \gamma^2 \ge C_0$) i.e. 
$\gamma \gg \mu^{-\frac{1}{2}}$ (or at least $\gamma \ge C_0\mu^{-\frac{1}{2}}$).

One possible approach is to apply in the \emph{singular zone\/}\index{zone!singular} $\{x:\ \gamma (x) \le \bar{\gamma}\}$ the rough estimate $C\mu h^{1-d} \upmu (\bar{\gamma})$ which works rather well as 
\begin{equation}
\upmu (\gamma)= O(\gamma^{\kappa}) \qquad \text{with\ \ } \kappa >2;
\label{20-8-11}
\end{equation}
then we can take $\bar{\gamma}=\mu^{-1/\kappa}$. 

In the \emph{regular zone\/}\index{zone!regular} 
$\{x:\, \gamma (x) \ge \bar{\gamma}\}$  one can go to the same canonical form as without degeneration but instead of error $O(\mu^{-m})$ error would become $O(\mu^{-m}\gamma^{-2m})$ in the principal part, and instead of $O(\mu^{-m}h^n)$ error would become $O(\mu^{-m}h^n \gamma^{-2m-n})$ in the lower terms.

Then we conclude that
\begin{claim}\label{20-8-12}
As $\gamma (y)\ge C_0\mu^{-\frac{1}{2}}$ the contribution of $B(y,\gamma(y))$ to the Tauberian remainder  does not exceed $O\bigl(h^{1-d}\gamma^{d-1}\bigr)$ as $q\ge 3$, 
$O\bigl(h^{1-d}\gamma^{d-1}+ \mu h^{\frac{5}{3}-d}\gamma^{d-\frac{2}{3}}\bigr)$ as $q=2$, and 
$O\bigl(h^{1-d}\gamma^{d-1}+ \mu h^{\frac{4}{3}-d}\gamma^{d-\frac{1}{3}}\bigr)$ as $q=1$.
\end{claim}

Then summation with respect to partition results in 
\begin{claim}\label{20-8-13}
Under assumption (\ref{20-8-11})  Tauberian remainder estimate is  $O\bigl(h^{1-d}\bigr)$ as $q\ge 3$,  
$O\bigl(h^{1-d}+\mu h^{\frac{5}{3}-d}\bigr)$ as $q=2$, and 
$O\bigl(h^{1-d}+\mu h^{\frac{4}{3}-d}\bigr)$ as $q=1$.
\end{claim}

Here we took into account that the contribution of the \emph{singular zone\/}\index{zone!singular} 
$\{x:\, \gamma (x) \ge \bar{\gamma}\}$ is
$O(\bar{\gamma}^\kappa \mu h^{1-d})=O(h^{1-d})$.

Let us calculate errors. First of all, in the regular zone removal of $O(\mu^{-1}\gamma^{-2})$ terms  brings an error 
$O( (\mu^{-1}\gamma^{-2})^{q/2} h^{-d}\gamma^{d})$ as $q\ge 2$ (we apply all the arguments used to and then summation with respect to partition results in $O(h^{1-d})$ under assumption (\ref{20-8-11}). The same is true if we remove $O(\mu^{-2}\gamma^{-4})$ terms as $q=1$.

As $q=1$ and there are unremovable $O(\mu^{-1}\gamma^{-2})$ terms we just go after $\R^\MW_\corr$; as our goal is to estimate a correction we use the weak magnetic field estimate as $\mu \le h^{-\frac{1}{2}}$ and also as 
$\mu \ge  h^{-\frac{1}{2}}$ but $\gamma \le (\mu^2 h)^{-\frac{1}{2}}$; in the latter case summation over such partition elements results in $O(\mu^{\frac{1}{2}}h^{1-d})$. 

Further, as $\gamma \ge (\mu^2 h)^{-\frac{1}{2}}$ and we estimate a correction we follow the arguments of Section~\ref{sect-20-6}. Finally we use a mollification parameter $\varepsilon =(\mu\gamma^2)^{-1}$ after rescaling or 
$\varepsilon =(\mu\gamma )^{-1}$ before. Leaving many rather not-very difficult but still delicate details to the reader we formulate the following

\begin{Problem}\label{Problem-20-8-5}
Prove that under assumption \textup{(\ref{20-8-11})} theorems~\ref{thm-20-6-17} and~\ref{thm-20-6-19}(i), (ii) remain true\footnote{\label{foot-20-31} Obviously, condition \textup{(\ref{20-6-78})} in theorem~\ref{thm-20-6-19}(ii) fails.}.
\end{Problem}

\subsection{Results under microhyperbolicity or non-degeneracy assumptions}
\label{sect-20-8-2-2}

It becomes a bit difficult to formulate a microhyperbolicity or non-degeneracy assumptions. First we assume that 
\begin{claim}\label{20-8-14}
$f_j$ ($j=2,\ldots,r$) are disjoint from $0$, $f_1\asymp \gamma(x)=\dist (x, Y)$ where $Y$ is $\sC^{2,1}$-manifold of codimension $3$ and $|\nabla f_1|\asymp 1$.
\end{claim}

Then without any loss of the generality one can assume that 
\begin{claim} \label{20-8-15}
At points of $Y$ matrices $(g^{jk})$ and $(F_{jk})$ satisfy $g^{jk}=0$  as $j=1,2,3$, $k=4,\ldots, d$ (and symmetrically) and $F_{jk}=0$ as $j=1,2,3$, $k=1,\ldots, d$ (and symmetrically).
\end{claim}

Let us consider submatrices $\mathsf{g}''$, $\mathsf{F}''$ of these matrices, consisting of elements $g^{jk}$, $F_{jk}$ with $j,k=4,\ldots,d$. Then canonical form is
\begin{multline}
\mathsf{a} (y,z;\xi';\zeta_1,\ldots,\zeta_r)=\\
\mathsf{g}(y,z; \xi') + f_1 (y,z)|\zeta_1|^2 + 
\mathsf{a}''(y,z; z_2,\ldots, \zeta_r) + V(y,z)+ \ldots
\label{20-8-16}
\end{multline}
where $y$ are coordinates along $Y$ and $z=(z_1,z_2,z_3)$ are additional coordinates (so $Y=\{z=0\}$) and  
\begin{equation}
|\nabla _z f_1|\asymp 1,\qquad |\nabla_y f_1|=O(|z|),
\label{20-8-17}
\end{equation}
$\mathsf{g}(y,z; \eta)$ is a non-degenerate quadratic form with respect to $\eta$,   $\mathsf{a}'(y,z; z_2,\ldots, \zeta_r)$ is a Hermitian form with respect to $\zeta_2,\ldots,\zeta_r$, .

Further, $f_1^2= b(y;z) +O(|z|^3)$ where $b(y;z)$ is a non-degenerate quadratic form with respect to $z$ and without any loss of the generality one can assume that $b(y;z)=|z|^2$. Then in the ball $B((y,z), \epsilon_0 |z|)$ we can use microhyperbolicity arguments as long as either $|\xi'|\ge \epsilon_0$ or $|\zeta_1|\ge C_0$.

Therefore if we are close to the energy level $0$, we need to consider only $|\xi'|\le \epsilon_0$ and $|\zeta_1|\le C_0$ but then the derivatives with respect to $z$ of the first two terms in the right-hand expression of (\ref{20-8-16}) are $O(|z|)$ and we need to consider only derivatives of  
\begin{equation}
\mathsf{a}''(y,z; \zeta_2,\ldots, \zeta_r)+V(y,z)
\label{20-8-18}
\end{equation}
as long as it close to our energy level $0$.

Therefore we can apply microhyperbolicity arguments there as well provided $\psi$ is supported in the small tubular vicinity of $Y$ and (\ref{20-8-18}) satisfies microhyperbolicity assumption with respect to $y$ (as $z=0$). Leaving many rather not-very difficult but still delicate details to the reader we formulate Problem~\ref{Problem-20-8-6}(i) below.

On the other hand, assume that (\ref{20-8-18}) is not microhyperbolic. Assume that
\begin{claim}\label{20-8-19}
Either $r=2$ or $f_2|_Y,\ldots, f_r|_Y$ are constant.
\end{claim}
Further, assume that $V/f_2|_Y$ has only non-degenerate critical points. Without any loss of the generality one can assume that $f_2=\const$ and there is just one critical point $0$. One can see easily that if $\mathsf{g}(y,z;\xi')=\bar{\mathsf{g}}(z;\xi')+O(|y|^2)$ and $b(y;z)=\bar{b}(z)+O(|y|^2)$ then 
$|\nabla_y \mathsf{a} (y,z;\xi';\zeta_1,\ldots,\zeta_r)|\asymp |y|$ at energy level $0$ and we can apply non-degeneracy arguments easily. In the general case however one can see easily that at the energy level $0$
\begin{equation*}
|\nabla_y \mathsf{a} (y,z;\xi';\zeta_1,\ldots,\zeta_r)| \ge \epsilon |y-\bar{y}(z,\xi',|\zeta_1|)|
\end{equation*}
and again we can apply non-degeneracy arguments.  Leaving many rather not-very difficult but still delicate details to the reader we formulate Problem~\ref{Problem-20-8-6}(ii) below. 

\begin{Problem}\label{Problem-20-8-6}
Let assumptions (\ref{20-8-14})--(\ref{20-8-15}) be fulfilled and  let  a cut-off function $\psi$ be supported in the sufficiently small vicinity of $Y$.

\begin{enumerate}[label=(\roman*), fullwidth]
\item Introduce a microhyperbolicity assumption as a microhyperbolicity assumption for these matrices $\mathsf{g}''$, $\mathsf{F}''$ and $V$  along $Y$ and prove under this assumption the remainder estimate $O(h^{1-d})$; 

\item Also, as (\ref{20-8-19}) is fulfilled, formulate a non-degeneracy assumption as non-degeneracy of critical points of $V/f_2|_Y$ and  prove under this assumption the remainder estimate $O(h^{1-d})$.
\end{enumerate}
\end{Problem}

\subsection{Other improved results}
\label{sect-20-8-2-3}

Situation becomes even more delicate as only assumptions (\ref{20-8-14})--(\ref{20-8-15}) are fulfilled. Let us assume that either $r=2$ or $f_2|_Y, \ldots,f_r|_Y$ are constant. Then we can use ``microhyperbolicity'' arguments and take $T_*=\epsilon \mu^{-1}$ unless gradient with respect to $z$ is small enough which is the case as in the decomposition $|\alpha_1-\bar{\alpha}_1|\ge C|\log h|$ which boils down to an extra factor $C\mu h|\log h|$ in the remainder estimate. 

Using rescaling technique we can get of the logarithm improving this factor to $C\mu h$.

\begin{Problem}\label{Problem-20-8-7}
Let assumptions (\ref{20-8-14})--(\ref{20-8-15}), and \textup{(\ref{20-8-19})} be fulfilled and let  a cut-off function $\psi$ be supported in the sufficiently small vicinity of $Y$.
\begin{enumerate}[label=(\roman*), fullwidth]
\item Prove remainder estimates 
\begin{equation}
\R^\MW \le Ch^{1-d} + C\mu^2
\left\{\begin{aligned}
&h^{\frac{8}{3}-d}\qquad &&\text{as\ \ } q=2,\\
&h^{\frac{7}{3}-d}\qquad &&\text{as\ \ }q=1;
\end{aligned}\right.\label{20-8-20}
\end{equation}
\item Provided there are no $3$-rd- order resonances prove remainder estimates $\R^\MW \le Ch^{1-d} $ as $q=2$ and
\begin{equation}
\R^\MW\le  Ch^{1-d} + C\mu^2  h^{\frac{5}{2}-d}\qquad \text{as\ \ }q=1;\label{20-8-21}\\
\end{equation}
\item Improve  estimate (\ref{20-8-21}) depending on $(l,\sigma)$;

\item Derive estimates for $\R^\W_{(\infty)}$ improving those of Subsection~
\ref{sect-20-8-1}; here and in (v) one should consider also $q\ge 3$;

\item Under microhyperbolicity or non-degeneracy assumptions of Subsubsection~\ref{sect-20-8-2-1}.1 further improve these estimates.
\end{enumerate}
\end{Problem}

\section{Strong and very strong magnetic field}
\label{sect-20-8-3}

Strong magnetic field basically repeats the intermediate magnetic field case leading us to part (i) below; in the case of the very strong magnetic field microhyperbolicity or non-degeneracy conditions should be modified and Problem~\ref{Problem-20-8-7} should be skipped, leading us to part (ii) below:

\begin{Problem}\label{Problem-20-8-8}
\begin{enumerate}[label=(\roman*), fullwidth]
\item Solve problems~\ref{Problem-20-8-5}, \ref{Problem-20-8-6}, and \ref{Problem-20-8-7} in the case of the strong magnetic field;

\item Problems~\ref{Problem-20-8-5} and \ref{Problem-20-8-6} in the case of the very strong magnetic field; in this case  microhyperbolicity and non-degeneracy conditions are restricted to $\alpha'=(\alpha_2,\ldots,\alpha_r)$ such that 
$|V+ \sum_j (2\alpha_j+1) f_j\mu h|\le \epsilon$.
\end{enumerate}
\end{Problem}

\section{Superstrong strong magnetic field}
\label{sect-20-8-4}

Theory of the superstrong magnetic field comes in two flavors: for the Schr\"odinger operator and for the Schr\"odinger-Pauli operator. 

\subsection{Schr\"odinger operator}
\label{sect-20-8-4-1}
In the  case of the Schr\"odinger operator the classically allowed zone is empty unless all $f_j$ vanish at $Y$; we assume  for simplicity hat 
\begin{equation}
f_j(x)\asymp \gamma(x) =\dist (x,Y) \qquad \forall j=1,\ldots,r
\label{20-8-22}
\end{equation}
where $Y$ is $\sC^{2,1}$ manifold of codimension $3$. Then the classically allowed zone is where after rescaling 
$\mu \mapsto \mu_\new=\mu \gamma^2$ and $h\mapsto h_\new =h/\gamma$ we have 
$\mu _\new h_\new \lesssim 1$ i.e. $\gamma \le \bar{\gamma}\Def 1/(\mu h)$. In the border we have $h_\new \asymp \mu h^2$, $\mu_\new \asymp h_\new^{-1}$. To be in the frame of the semiclassical theory we need to assume that $h_\new \ll 1$ i.e. $\mu \ll h^{-2}$.

\begin{Problem}\label{Problem-20-8-9}
Under assumption (\ref{20-8-22}) prove that for the Schr\"odinger operator

\begin{enumerate}[label=(\roman*), fullwidth]
\item Let $\mu  \ge C_0h^{-2}$; then $e(x,y,0)=O(\mu^{-\infty})$;

\item Let  $C_0 h^{-1}\le \mu \le C_0h^{-2}$; then the principal part of the asymptotics is given by a magnetic Weyl expression and has a magnitude of the  very strong magnetic field rescaled case multiplied by $\bar{\gamma}^{3-d}$, namely 
\begin{equation}
h_\new^{-d}\bar{\gamma}^{3-d} \asymp h^{-d} \bar{\gamma}^3\asymp \mu^{-3}h^{-3-d}
\label{20-8-23}
\end{equation}
while the error coincides with one of the very strong magnetic field case rescaled, multiplied by $\bar{\gamma}^2$:
\begin{align}
&\R^\MW \le C h_\new^{1-d}\bar{\gamma}^{3-d} \asymp C\mu^{-2}h^{-1-d}\qquad &&\text{as\ \ } q\ge 3,
\label{20-8-24}\\[3pt]
&\R^\MW \le C h_\new^{\frac{q}{3}-d}\bar{\gamma}^{3-d} \asymp C\mu^{-3+\frac{q}{3}}h^{-3+\frac{2q}{3}-d}
&& \text{as\ \ } q\le 3;
\label{20-8-25}
\end{align}
in the general case; 
\item Further, assuming that
\begin{equation}
|f_j -f_k|\asymp \gamma \qquad \forall j\ne k
\label{20-8-26}
\end{equation}
prove for $(l,\sigma)= (2,0)$ estimate (\ref{20-8-24}) as $q=2$ and 
\begin{equation}
\R^\MW \le h_\new^{\frac{1}{2}-d}\bar{\gamma}^{3-d} \asymp \mu^{-\frac{5}{2}}h^{-2-d}\qquad q=1;
\label{20-8-27}
\end{equation}
Also consider  $(l,\sigma)\prec (2,0)$ for $q=1,2$ and $(l,\sigma)\succ (2,0)$ for $q=1$;
\item Furthermore, introduce a  notions of microhyperbolicity and non-degeneracy and under these assumptions prove estimate (\ref{20-8-24}) as $q=1,2$.
\end{enumerate}
\end{Problem}

\subsection{Schr\"odinger-Pauli operator}
\label{sect-20-8-4-2} 

In the  of the Schr\"odinger-Pauli operator the classically allowed zone does not shrink and we  need neither assumption  that $\mu \le C_0h^{-2}$ nor that all $f_j$ vanish on $Y$; so we return to our previous assumption that only $f_1$ vanishes on $Y$. 

However in the inner zone $\{x:\, \gamma(x)\le C_0/(\mu h)\}$ we do not have a ``scalar operator case'' anymore and it can spoil the remainder estimate as $q=1,2$. It does not happen, however, it does not happen as $\mu \ge C_0h^{-2}$ as well. Thus we arrive to (i), (ii) of Problem~\ref{Problem-20-8-10} below.

However  we are essentially in the ``constant multiplicities of $f_j$'' case as $f_1$ is disjoint from all others and with $f_2,\ldots,f_r$ we are restricted to the lowest Landau level. Further, microhyperbolicity with respect to $z$ arguments work unless $\alpha_1 \le C_1$ and therefore we can formulate notions of the microhyperbolicity and non-degeneracy in the terms of $W|Y$.

\begin{Problem}\label{Problem-20-8-10}
Under assumptions (\ref{20-8-14})--(\ref{20-8-15})  for the Schr\"odinger-Pauli  operator
\begin{enumerate}[label=(\roman*), fullwidth]
\item 
Prove that the principal part of the asymptotics is as if there was no degeneration at all (i.e.  of magnitude $\mu^r h^{r-d})$);

\item Prove that as $q\ge 3$   the remainder estimate is as if there was no degeneration at all (i.e.   $O\bigl(\mu^r h^{1+r-d}\bigr)$);

\item Formulate notions of the microhyperbolicity and non-degeneracy in the terms of $W|Y$ and prove that under either of these assumptions the remainder estimate is also $O\bigl(\mu^r h^{1+r-d}\bigr)$;

\item In the general case of $q=1,2$ as $(l,\sigma)= (2,0)$ prove that the remainder estimate is as if there was no degeneration at all (i.e.   $O\bigl(\mu^r h^{1+r-d}\bigr)$ for $q=2$ and 
$O\bigl(\mu^r h^{\frac{1}{2}+r-d}\bigr)$ for $q=1$); also consider  $(l,\sigma)\prec (2,0)$ for $q=1,2$ and $(l,\sigma)\succ (2,0)$ for $q=1$.
\end{enumerate}
\end{Problem}

\bibliographystyle{alpha}

\begin{thebibliography}{BrIvr}




\bibitem[1]{harman}
\textsc{G. Harman}.
 \newblock{Metric Number Theory, 
 Clarendon Press, Oxford},  1998, xviii+297.


\bibitem[Ivr1]{Ivr1}
\textsc{V.~Ivrii}.
 \emph{Microlocal Analysis and Precise Spectral Asymptotics}, 
 Springer-Verlag, SMM, 1998, xv+731.


\bibitem[Ivr2]{futurebook}
\textsc{V.~Ivrii}.
 \emph{Microlocal Analysis and Sharp Spectral Asymptotics}, 
 in progress: available online at \newline
\href{http://www.math.toronto.edu/ivrii/futurebook.pdf}{http://www.math.toronto.edu/ivrii/futurebook.pdf}
 



\bibitem[Ivr3]{ivrii:IRO5}
\href{http://weyl.math.toronto.edu:8888/victor2/preprints/IRO5.pdf}{
Sharp spectral asymptotics for operators with irregular coefficients.
V. Multidimensional Schr\"odinger operator with a strong magnetic field.
Non-full-rank case},  arXiv:math/0510328
(Aug. 17, 2004), 78pp.

\bibitem[Ivr4]{ivrii:IRO8}
\href{http://weyl.math.toronto.edu:8888/victor2/preprints/IRO8.pdf}{
Sharp spectral asymptotics for four-dimensional Schr\"odinger operator
with a strong degenerating magnetic field.} , arXiv:math/0605298  (May 11, 2006), 93pp.

\bibitem[Ivr5]{ivrii:IRO9}
\href{http://weyl.math.toronto.edu:8888/victor2/preprints/IRO9.pdf}{
Sharp spectral asymptotics for four-dimensional Schr\"odinger operator with a strong magnetic field. II},  arXiv:math/0612252   (December 9, 2006), 57pp
%
 \end{thebibliography}

\providecommand{\bysame}{\leavevmode\hbox to3em{\hrulefill}\thinspace}

\vglue .06truein

\begin{tabular}{rrl}
&{\hskip 200 pt} &Department of Mathematics,\cr
&&University of Toronto,\cr
&&40, St.George Str.,\cr
&&Toronto, Ontario M5S 2E4\cr
&&Canada\cr
&&ivrii@math.toronto.edu\cr
&&Fax: (416)978-4107\cr
\end{tabular}

\end{document}